\documentclass[hidelinks,onefignum,onetabnum]{siamart251216}
\usepackage{changepage}
\usepackage{diagbox}
\usepackage[caption=false]{subfig}
\usepackage{multirow}
\usepackage[percent]{overpic}
\usepackage{enumitem}
\usepackage{algorithm}
\usepackage{algorithmic}
\usepackage{setspace}
\usepackage{amsmath, amssymb}
\usepackage{booktabs}
\usepackage{tikz}
\usepackage{pgfplots}
\pgfplotsset{compat=newest}
\usepgfplotslibrary{colormaps}

\definecolor{mygray}{gray}{0.8}

\usepackage{soul}
\usepackage[normalem]{ulem}
\allowdisplaybreaks
\patchcmd\newpage{\vfil}{}{}{}
\usepackage{lipsum}
\usepackage{amsfonts}
\usepackage{graphicx}
\usepackage{epstopdf}
\usepackage{algorithmic}
\ifpdf
  \DeclareGraphicsExtensions{.eps,.pdf,.png,.jpg}
\else
  \DeclareGraphicsExtensions{.eps}
\fi

\newsiamremark{remark}{Remark}
\newsiamremark{hypothesis}{Hypothesis}
\crefname{hypothesis}{Hypothesis}{Hypotheses}
\newsiamthm{claim}{Claim}

\usepackage{amsopn}

\headers{Xiaobin Li, Leevan Ling, Yizhong Sun}
        {Local Regularity Estimation}
 
\title{Local Regularity Estimation through Sobolev-Scale Norm Profile
\thanks{ Submitted to
\funding{This work was supported by the General Research Fund (GRF No.~12301021, 12300922, 12301824) of the Hong Kong Research Grants Council.}}}
\author{ Xiaobin Li
\thanks{Department of Mathematics, Hong Kong Baptist University  ({25481800@life.hkbu.edu.hk}).}
\and
Leevan Ling\thanks{Department of Mathematics, Hong Kong Baptist University  ({lling@hkbu.edu.hk}).}
\and
Yizhong Sun\thanks{Department of Mathematics, Hong Kong Baptist University  ({yzsun95@hkbu.edu.hk}).}
}

\begin{document}
\maketitle
\begin{abstract}
We develop a kernel-based approach for estimating the spatially varying Sobolev regularity~$s$ of an unknown $d$-variate function~$f$ from scattered sampling data, which quantifies the degree of local differentiability supported by the data.
Relying only on neighborhood data near the point of interest $z\in \Omega_z$, our method constructs a sequence of Sobolev-space reproducing kernel interpolants whose kernel smoothness order is specified by an index~$m > d/2$.
The native-space norms of these interpolants are evaluated over a bounded range of~$m$, producing a \emph{Sobolev-scale norm profile}.
The elbow of this profile serves as a quantitative probe of the underlying local regularity~$s(\Omega_z)$.
In particular, when $m > s(\Omega_z)$, the profile exhibits rapid, near-worst-case growth governed by the classical upper bound associated with the conditioning of the kernel matrix.
A band-limited surrogate analysis explains this transition and establishes a lower bound linking native-norm growth to the Sobolev regularity of~$f$.
Two complementary strategies are incorporated for further enhancement:
(i)~a \emph{stencil-shift} subroutine, which repositions local neighborhoods to avoid crossing discontinuities whenever possible, thereby suppressing artifacts in the norm estimates; and (ii)
a \emph{secant-based tail screening strategy} that uses two high-order norm
evaluations to identify candidate low-regularity neighborhoods at reduced
computational cost.
Numerical experiments on synthetic test functions and turbulent-flow data
demonstrate recovery of spatially varying regularity, while a surface
conservation-law example illustrates the detection of evolving low-regularity
regions in time-dependent PDE data.
\end{abstract}


\section{Introduction}
\label{sec:Introduction}

Estimating spatially varying smoothness from scattered data is a fundamental problem in numerical analysis and data-driven modeling.
We consider an unknown $d$-variate function
\[
  f : \Omega \subset \mathbb{R}^d \to \mathbb{R},
  \qquad
  (X, f(X)) = \{(\mathbf{x}_i, f(\mathbf{x}_i))\}_{i=1}^{N},
\]
where $X \subset \Omega$ is a finite set of sampling locations.
Given a point of interest $z \in \Omega$, we denote by $\Omega_z$ a small neighborhood containing $z$ and by $X_z = X \cap \Omega_z$ the corresponding subset of samples.
Our goal is to estimate the \emph{local Sobolev regularity}
\[
  s = s(\Omega_z)
  \quad \text{such that} \quad
  f \in H^{s}(\Omega_z),
\]
that is, the Sobolev order for which all weak derivatives $D^{\alpha}f$ with $|\alpha| < s(\Omega_z)$ lie in $L^2(\Omega_z)$.
The function $s(\Omega_z)$ measures the degree of differentiability supported by the data near $z$, for example, a jump discontinuity admits $s < \tfrac 12$, while a corner or slope discontinuity yields $s < \tfrac 32$.
Knowledge of $s(\Omega_z)$ is valuable for detecting nonsmooth regions, selecting reconstruction stencils, and guiding adaptive solvers. We focus here on deterministic data, including discrete solution values
generated by numerical PDE solvers. The treatment of
measurement noise would require a regularized formulation in place of exact
interpolation and is beyond the scope of the present study.

Traditional edge or singularity detection methods locate discontinuities via local peaks in interpolation coefficients or derivatives, often through polynomial or transform-based approaches such as ENO/WENO schemes, wavelet transforms, or kernel filters~\cite{canny1986computational,demarchi2020shape,de2020jumping,jung2009iterative,jung2011iterative,shrivakshan2012comparison}.
While classical wavelets~\cite{Mallat1999Wavelet} and high-order reconstructions like ENO/WENO~\cite{HartenEngquistOsherChakravarthy1987ENO,JiangShu1996WENO,HuShu1999WENOUnstructured} are highly effective, many of their standard formulations rely on structured
or mesh-based stencils, which are not directly available for arbitrary scattered
data.
Moreover, these techniques yield mostly qualitative, often binary, indicators of smoothness rather than quantitative measures of differentiability. Polynomial-annihilation techniques provide another local approach for
detecting nonsmooth features from irregular samples \cite{ArchibaldGelbYoon2005}. Our objective is different: we seek a quantitative estimate of local
Sobolev regularity rather than only the location of nonsmooth features.

Multiscale kernel constructions and samplet analysis provide complementary perspectives on localized approximation. Hierarchical and sparse-grid kernel schemes~\cite{Kempf+Wendland-Highapprwithkern:23,Opfer-Multkern:06,Wendland-MultanalSobospac:10} introduce explicit multilevel structure into RBF approximation, enabling efficient representation of multiscale features on scattered data. In parallel, samplets~\cite{Avesani2025Smoothness,Harbrecht2022Samplets} combine the vanishing-moment property of wavelets with the geometric flexibility of meshfree methods. Constructed via hierarchical cluster trees, they form localized, distributional wavelets whose coefficients decay at rates governed by local smoothness. Recent studies further demonstrate that the samplet transform can compress and sparsify generalized Vandermonde matrices arising in multiscale RBF interpolation, particularly for globally supported kernels such as the Whittle-Mat\'ern family~\cite{Avesani2024Multiscale}. 

These approaches primarily address feature detection or multiscale
representation. Here we instead seek a quantitative estimate of local Sobolev
regularity directly from scattered samples.
In this paper, we propose a different approach: directly estimating Sobolev regularity from arbitrary local stencils of scattered samples via a \emph{Sobolev-scale norm profile}.
Section~\ref{sec:kernel-detector} defines this profile and introduces an elbow-based estimator of~$s(\Omega_z)$.
Section~\ref{sec:nativenorm-bandlimited} develops the theoretical foundation linking native-norm growth to Sobolev regularity through Sobolev–native norm relations, band-limited surrogate analysis, and inverse inequalities. The main result, Theorem~3.4, gives a lower bound for $\eta(m)$ in terms of a
minimum-norm band-limited interpolant and a spectral-median factor, quantifying
the additional norm growth caused by imposing smoothness beyond a reference
Sobolev order. Section~\ref{sec:regularity_aware_algorithms} presents refinement strategies,
including the stencil-shift adjustment that reduces spatial smearing and a low-cost secant-based tail screening strategy.
Finally, Section~\ref{sec:Numerical experiments} presents synthetic and
turbulent-flow tests, together with an application to time-dependent surface
conservation-law data, where the detector is used to identify evolving
low-regularity regions.


\section{Local Regularity Estimation via Sobolev-Scale Norm Profiles}
\label{sec:kernel-detector}

In many computational settings, function values are available only at irregular locations, and the underlying field may exhibit nonsmooth features.
Standard kernel interpolation methods often exhibit oscillatory behavior near such features.
Rather than treating this instability as a numerical defect, we interpret it as diagnostic information about local loss of regularity.

This section introduces a kernel-based approach for quantifying local smoothness.
The key idea is to examine, within each interior neighborhood, how the native norm of a kernel interpolant varies across a family of Sobolev-scale kernels.
The resulting dependence on the kernel smoothness parameter~$m$ defines a local interpolant norm profile, which serves as a quantitative probe of the underlying regularity.

\subsection{Preliminaries}
\label{sec:kernel-prelim}

We consider a family of translation-invariant, symmetric, positive definite (SPD) kernels
$\{\Phi_m\}_{m > d/2}$ on $\mathbb{R}^d$, parameterized by the smoothness index $m$,
where $\Phi_m(x,z) = \phi_m(\|x-z\|_2)$ for a radial function $\phi_m:[0,\infty)\to\mathbb{R}$.
These kernels are characterized by the algebraic decay of their Fourier transforms,
\begin{equation}
c_1 (1+\|\omega\|_2^2)^{-m}
\le \widehat{\Phi}_m(\omega)
\le c_2 (1+\|\omega\|_2^2)^{-m},
\qquad \omega \in \mathbb{R}^d,
\label{eq:fourier-decay}
\end{equation}
for constants $0 < c_1 \le c_2 < \infty$ that remain uniform over bounded intervals of $m > d/2$.

A standard example is the Whittle-Mat\'ern kernel family,
\begin{equation}\label{eq:wmat-def-ex}
\Phi_m(x,z)
= \|x-z\|_2^{\,m-d/2}\,
  \mathcal{K}_{m-d/2}(\|x-z\|_2),
  \qquad m > d/2,
\end{equation}
where $\mathcal{K}_{m-d/2}$ denotes the modified Bessel function of the second kind.
Its associated reproducing kernel Hilbert space (also called native space)
$\mathcal{N}_{\Phi_m}(\mathbb{R}^d)$
is norm-equivalent to the Sobolev space $H^{m}(\mathbb{R}^d)$~\cite[Thm.~10.47]{wendland2004scattered}.
If a shape parameter $\varepsilon>0$ is introduced via
$\Phi_{m,\varepsilon}(x,z) = \Phi_m(\varepsilon\|x-z\|_2)$,
this norm equivalence remains valid.

As in the introduction, let $X_z$ denote the local neighborhood of data sites whose convex hull contains the evaluation point $z\in\Omega$.
We record a few geometric quantities used in later analysis.
The fill distance of $X$ with respect to $\Omega_z$ and the separation distance of $X$ are
\begin{equation}
h_{X,\Omega_z} := \sup_{x\in\Omega_z}\min_{x_j\in X}\|x-x_j\|_2,
\qquad
q_X := \tfrac12 \min_{j\ne k}\|x_j-x_k\|_2 ,
\label{eq:fill-sep-local}
\end{equation}
and the local mesh ratio is $\rho_{X,\Omega_z} := h_{X,\Omega_z}/q_X$.
From now on, we drop the subscript $z$ and keep in mind that all data are local and located near $z$.

Let $Y = f|_X$ denote the data values of an unknown function $f$ at the sites $X$.
The finite-dimensional kernel trial space associated with $X$ is
\begin{equation}
\mathcal{U}_{X,\Phi_m}
:= \operatorname{span}\{\Phi_m(\cdot, x_j): x_j\in X\}.
\label{eq:trial-space}
\end{equation}
The interpolant $u_m \in \mathcal{U}_{X,\Phi_m}$ satisfies
$u_m|_X = Y$ and can be written as
\[
  u_m(x) = \sum_{j=1}^n \alpha_j \Phi_m(x, x_j),
\]
where the coefficient vector $\boldsymbol{\alpha}$ solves
\[
  \Phi_m(X,X)\, \boldsymbol{\alpha} = Y,
\]
and $\Phi_m(X,X)\in\mathbb{R}^{n\times n}$ has entries
$[\Phi_m(X,X)]_{jk} = \Phi_m(x_j, x_k)$.
The squared native norm of the interpolant is then
\begin{equation}
\|u_m\|_{\mathcal{N}_{\Phi_m}}^2
  = \boldsymbol{\alpha}^T \Phi_m(X,X)\, \boldsymbol{\alpha}
  = Y^T [\Phi_m(X,X)]^{-1} Y.
\label{eq:native-norm-matrix}
\end{equation}
By the standard native-space identity for kernel expansions,
\eqref{eq:native-norm-matrix} expresses the native norm directly in terms of the data pair
$(X,Y)$ and the kernel order $m$ \cite[Ch.~10]{wendland2004scattered}. This quantity forms the basic building block of the interpolant norm profile introduced in the next subsection.

\begin{figure}[t]
\centering
\includegraphics[width=0.75\textwidth]{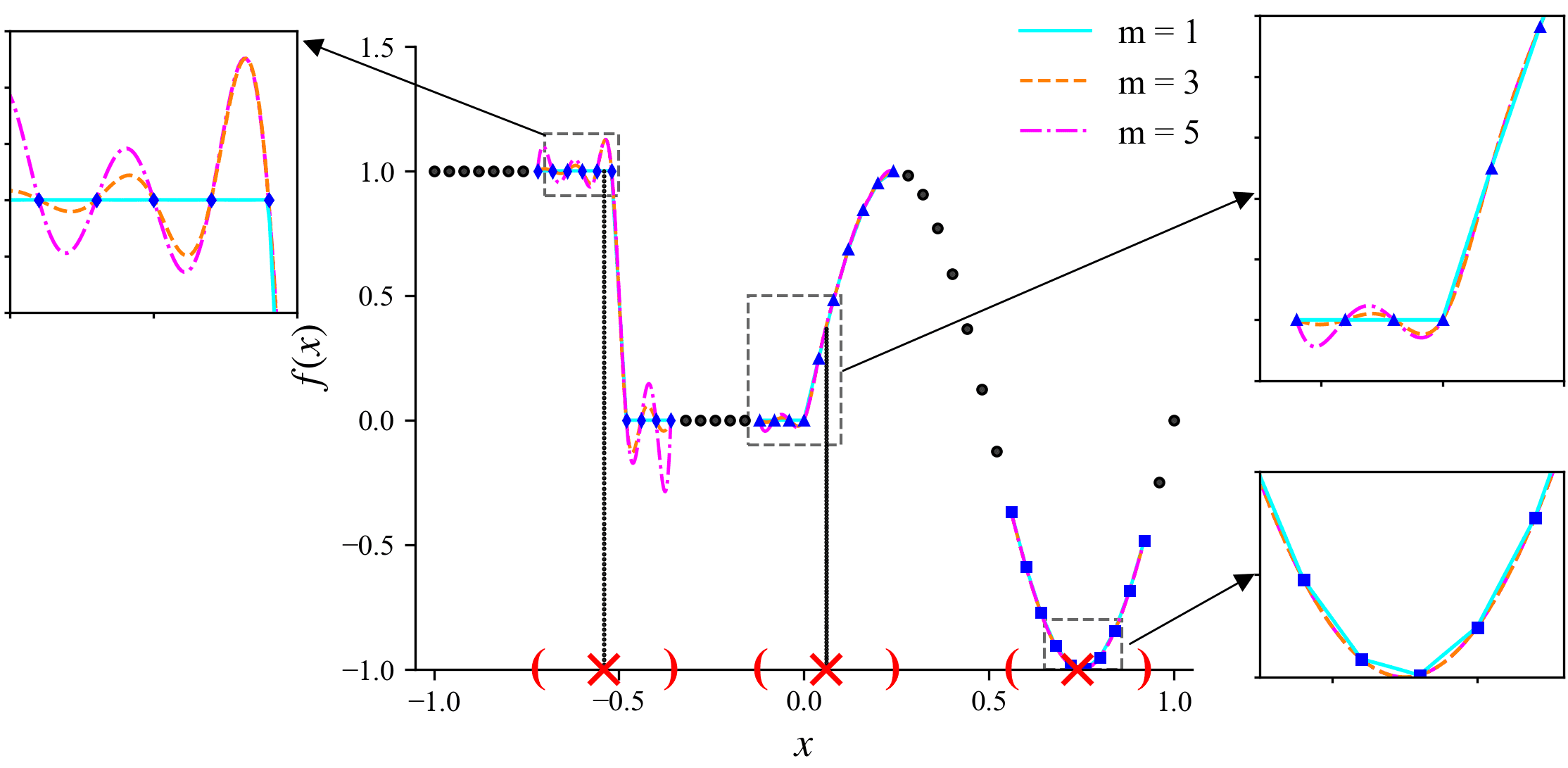}
\caption{
One-dimensional RBF interpolation on local stencils using Whittle--Mat\'ern kernels with smoothness
$m=1,3,5$.
The test function is $f(x)=\sin(2\pi x)$ for $x>0$, $f(x)=0$ for $-0.5<x\le0$, and $f(x)=1$ for $-1\le x\le-0.5$.
A total of 51 sample points are used, and each local interpolant is computed on a stencil of 10 points.
}
\label{fig:oscillation}
\end{figure}
Although the interpolation condition $u_m|_X=Y$ holds for all $m$, the off-sample behavior of $u_m$ depends strongly on the kernel smoothness.
Larger $m$ enforces higher smoothness and may produce Gibbs-type oscillations near nonsmooth regions,
while a smaller $m$ yields rougher approximants that can underfit smooth areas.
Figure~\ref{fig:oscillation} illustrates this dependence and motivates examining the variation of the native norm with respect to $m$ as an indicator of local regularity.

\subsection{Sobolev-Scale Norm Profile and Data-Driven Regularity}
\label{sec:sobolev-scale-profile}

To quantify local smoothness from scattered samples, we examine how the native norm of a kernel interpolant varies over the Sobolev scale.
For each local stencil $X$ with associated function values $Y$ in the neighborhood of an interior target $z$, we apply a local normalization so that the data vector $Y$ has unit discrete $\ell^2_h$-norm.
This normalization is performed separately for each stencil, making the resulting Sobolev-scale norm profile $\eta(m)$ (defined later in \eqref{eq:native-norm-profile}) invariant to local amplitude and directly comparable across neighboring locations.

The diagnostic role of varying kernel smoothness was already illustrated in
Figure~\ref{fig:oscillation}.
We consider the family of interpolants
\[
  u_m \in \mathcal U_{X,\Phi_m}, \qquad
  m \in [m_{\min},m_{\max}] \subset (d/2,\infty),
\]
constructed for the data $(X,Y)$ using the kernel interpolation method introduced in Section~\ref{sec:kernel-prelim}.
Their native norms are evaluated as in~\eqref{eq:native-norm-matrix}, and we define the
\emph{Sobolev-scale  norm profile} for $z$ to be the native space norm of these interpolants:
\begin{equation}
  \eta(m) := \|u_m\|_{\mathcal{N}_{\Phi_m}}
=\sqrt{Y^T[\Phi_m(X,X)]^{-1}Y},
  \label{eq:native-norm-profile}
\end{equation}
which depends only on the normalized data and on the kernel smoothness order~$m$.
The normalization removes amplitude effects so that $\eta(m)$ primarily reflects how the
shape of the local samples interacts with the assumed smoothness order.

The maximal admissible growth of $\eta(m)$ can be bounded using the smallest eigenvalue of
$\Phi_m(X,X)$ and the local separation distance $q_X$:
\begin{equation}
  \eta(m) \le \|Y\|_2 \lambda_{\min}^{-1/2}(\Phi_m(X,X))
  \le c_{m,d}\,\|Y\|_2\, q_X^{\,d/2-m},
  \label{eq:upper-bound}
\end{equation}
where $c_{m,d}>0$ depends only on the kernel family and the dimension.
This classical stability bound~\cite{narcowich2006sobolev} provides a
worst-case reference rate describing the steep regime of $\eta(m)$.

\begin{figure}
\centering
\subfloat[$q_{X} = 0.1$]{%
\begin{overpic}[width=0.32\textwidth]{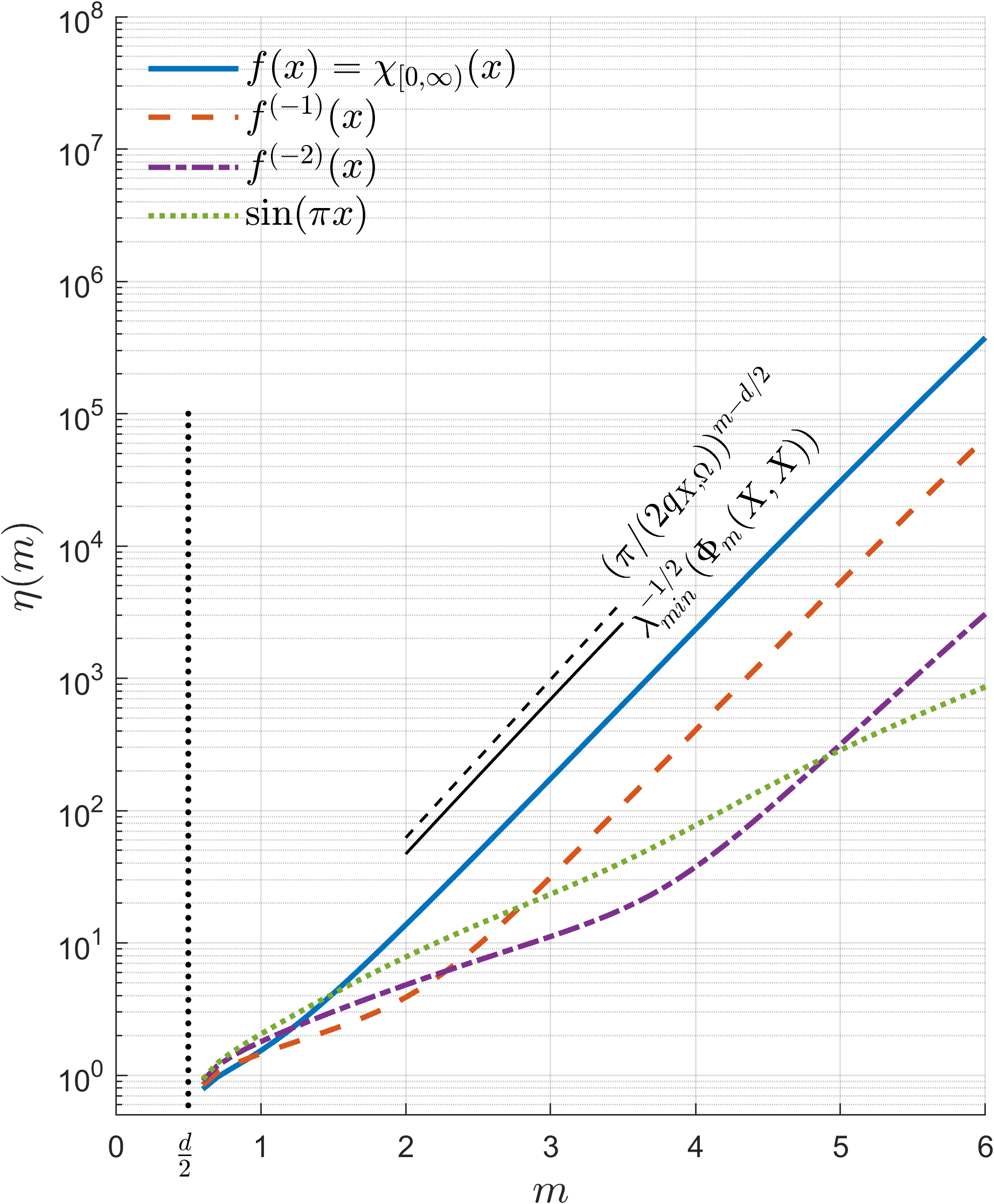}
\end{overpic}}
\hspace{0.1mm}
\subfloat[$q_{X} = 0.05$]{%
\begin{overpic}[width=0.32\textwidth]{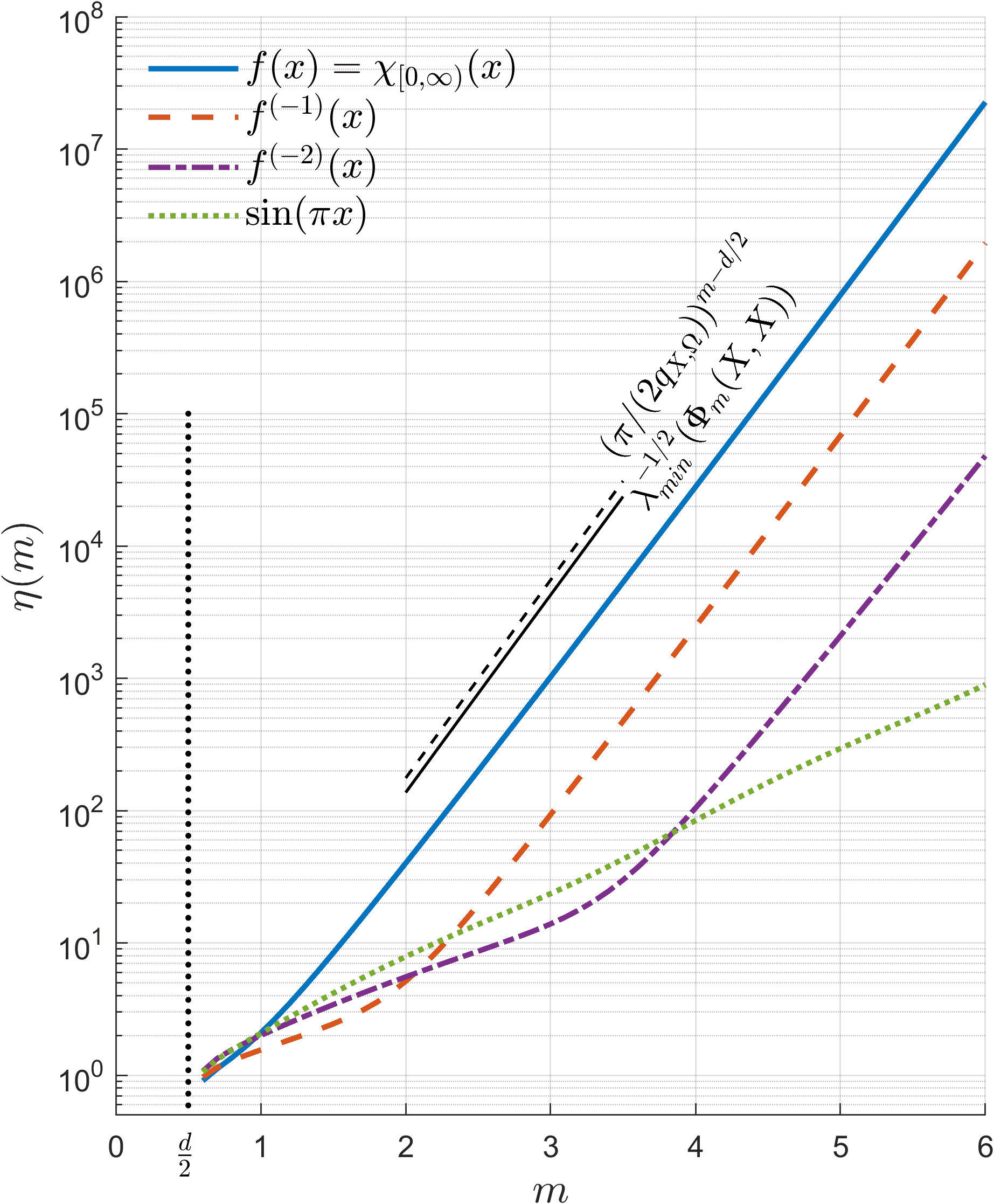}
\end{overpic}}
\hspace{0.1mm}
\subfloat[$q_{X} = 0.025$]{%
\begin{overpic}[width=0.32\textwidth]{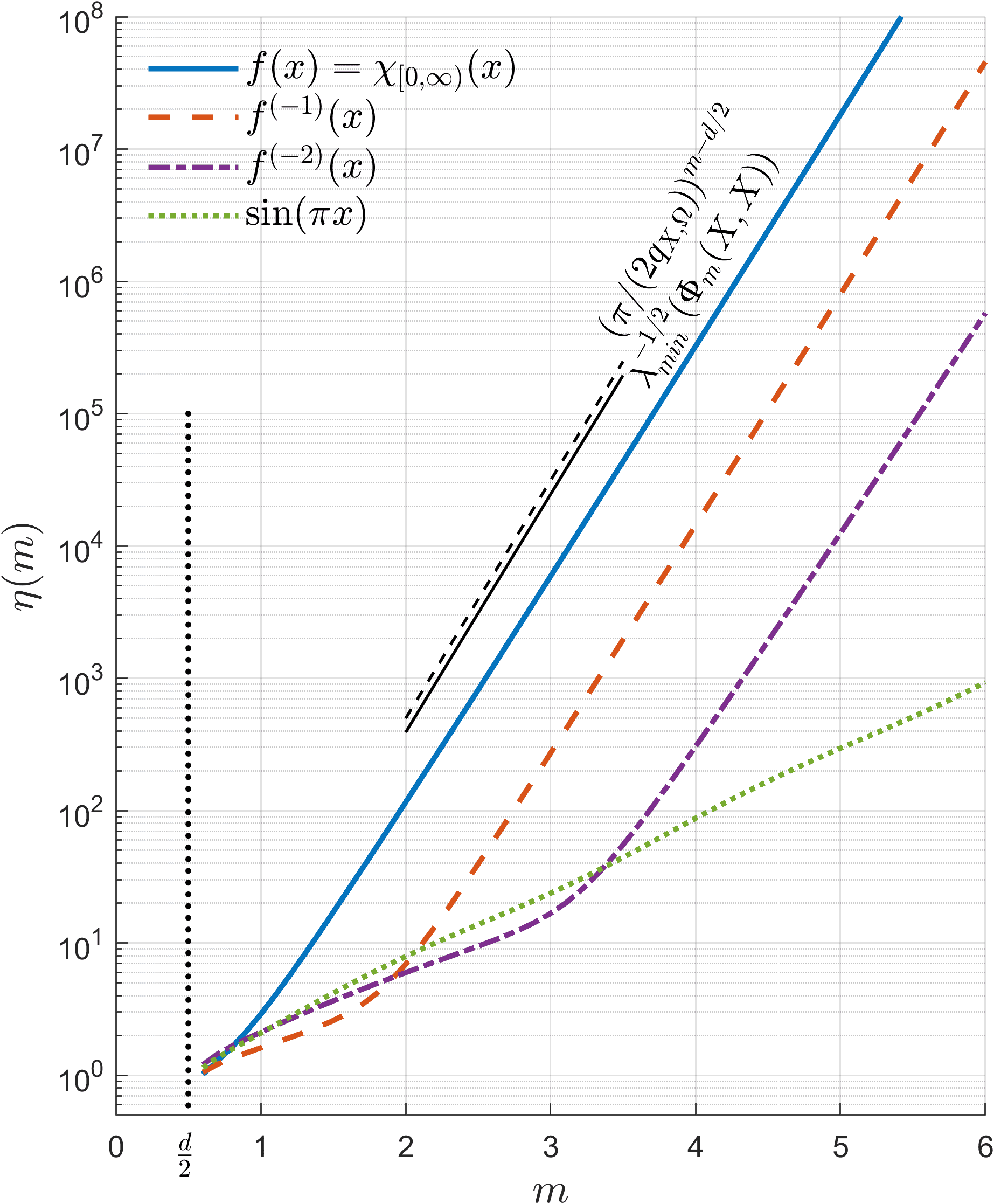}
\end{overpic}}

\caption{
Native‑norm profiles $\eta(m)$ versus kernel smoothness $m$ for RBF interpolants centered at $z=0$, evaluated on local neighborhoods $\Omega_z=[-1,1]$ constructed from point sets with separation distances $q_X=0.1,0.05,$ and $0.025$.
Colored curves correspond to a step function $f(x)=\chi_{(0,\infty)}(x)$, first- and second-order kink
(given by the first and second antiderivatives $f^{(-1)}$ and $f^{(-2)}$ of $f$), and a smooth reference $\sin(\pi x)$.
Dashed and solid black lines show $\lambda_{\min}^{-1/2}(\Phi_m(X,X))$ and its approximation
$(\pi/(2q_{X}))^{\,m-d/2}$, respectively, illustrating the worst-case rate in \eqref{eq:upper-bound}.
}

\label{fig:norm-growth}
\end{figure}

Figure~\ref{fig:norm-growth} illustrates the upper bound~\eqref{eq:upper-bound} for one-dimensional functions
with local regularity~$s=0.5, 1.5, 2.5,$ and~$\infty$
on uniformly spaced samples ($h=2q_X$), where the asymptotic rate $q_X^{\,d/2-m}$ is well approximated by $(\pi/(2q_X))^{\,m-d/2}$.
This follows from a simple one‑dimensional heuristic: for evenly spaced nodes with spacing $h=2q_X$, the largest representable frequency is $\omega_{\max}=\pi/(2q_X)$, so the proxy growth $\omega_{\max}^{\,m-d/2}$ yields exactly this approximation.

For lower regularity, where the true local order satisfies $s<m_{\min}\approx d/2$,
the profile $\eta(m)$ grows immediately at the worst-case rate for all tested $m$,
resulting in a nearly linear log-profile with no distinct elbow.

If $s \in [m_{\min},m_{\max}]$, the profile $\eta(m)$ typically evolves through two regimes:
for small $m$, the values vary moderately, whereas beyond a certain point they increase rapidly.
Intuitively, if the underlying function is locally in $H^{s}$, then for $m<s$ the data remain compatible with the kernel space and the native norm stays nearly stable.
Once $m>s$, the interpolant enforces excessive smoothness, causing $\eta(m)$ to grow quickly.
On a semi-logarithmic scale $(m,\log\eta(m))$, this behavior forms an L-shaped curve:
a flat region for $m\lesssim s$ followed by a steep, near power-law rise for
$m\gtrsim s$, close to the worst-case rate~\eqref{eq:upper-bound}.
We identify the corner of this transition as a surrogate for the largest
Sobolev order detectable from the samples and define the corresponding
\emph{L-curve elbow} via the curvature of the log-norm curve:
\begin{equation}
  m^*
  = \arg\max_m
  \frac{|(\log\eta(m))''|}
         {(1 + ((\log\eta(m))')^2)^{3/2}}.
  \label{eq:m-eta-star-def}
\end{equation}

For smooth data with $s>m_{\max}$, no elbow is expected within the tested
range. The tail may appear approximately linear on the log scale, but its
slope depends on the target function, the local stencil, and the kernel family.
We therefore interpret it only relative to the worst-case reference rate
in~\eqref{eq:upper-bound}.

\begin{definition}[Data-driven regularity]
\label{def:data-driven-regularity}
Let $\eta(m)$ denote the Sobolev-scale norm  profile evaluated over a
prescribed  interval of Sobolev orders $[m_{\min},m_{\max}]$, and $m^*$ be the
elbow defined in \eqref{eq:m-eta-star-def}.
Then, the \emph{data-driven regularity} $\tilde{s}(\Omega_z)$ is defined by
\begin{enumerate}
\item $\tilde{s}(\Omega_z)=m^* $, if a distinct elbow\footnote{%
A \emph{distinct elbow} refers to a point $m^*$ at which the log--norm profile
$(m,\log\eta(m))$ exhibits a clear change of slope, in the sense that the
difference between the local pre- and post-elbow slopes exceeds a chosen
contrast threshold (typically several times the numerical slope variation
within each region). In practice, this corresponds to a visually discernible
transition between the mild- and rapid-growth regimes of $\eta(m)$.%
} is present;
\item $\tilde{s}(\Omega_z)=m_{\min}$, if the profile is nearly linear with slope comparable to the worst case rate \eqref{eq:upper-bound};
\item $\tilde{s}(\Omega_z)=m_{\max}$, if the profile has no distinct elbow and its tail
slope is significantly smaller \footnote{``Significantly smaller'' is assessed by comparing the
observed tail slope of $\log\eta(m)$ with the slope of the worst-case reference
curve in \eqref{eq:upper-bound}. In our numerical tests, we use a user-chosen threshold $\tau_{\rm sm}\in(0,1)$ for this
ratio; values around $0.3$--$0.5$ were found stable in the tests below.}
than a prescribed fraction of the worst-case slope in
\eqref{eq:upper-bound}. 
\end{enumerate}
\end{definition}
A theoretical justification for this procedure is developed later in
Section~\ref{sec:nativenorm-bandlimited} using Sobolev--native norm relations and
band-limited surrogate estimates.

\begin{figure}[t]
  \centering
  \begin{overpic}[width=0.8\textwidth]{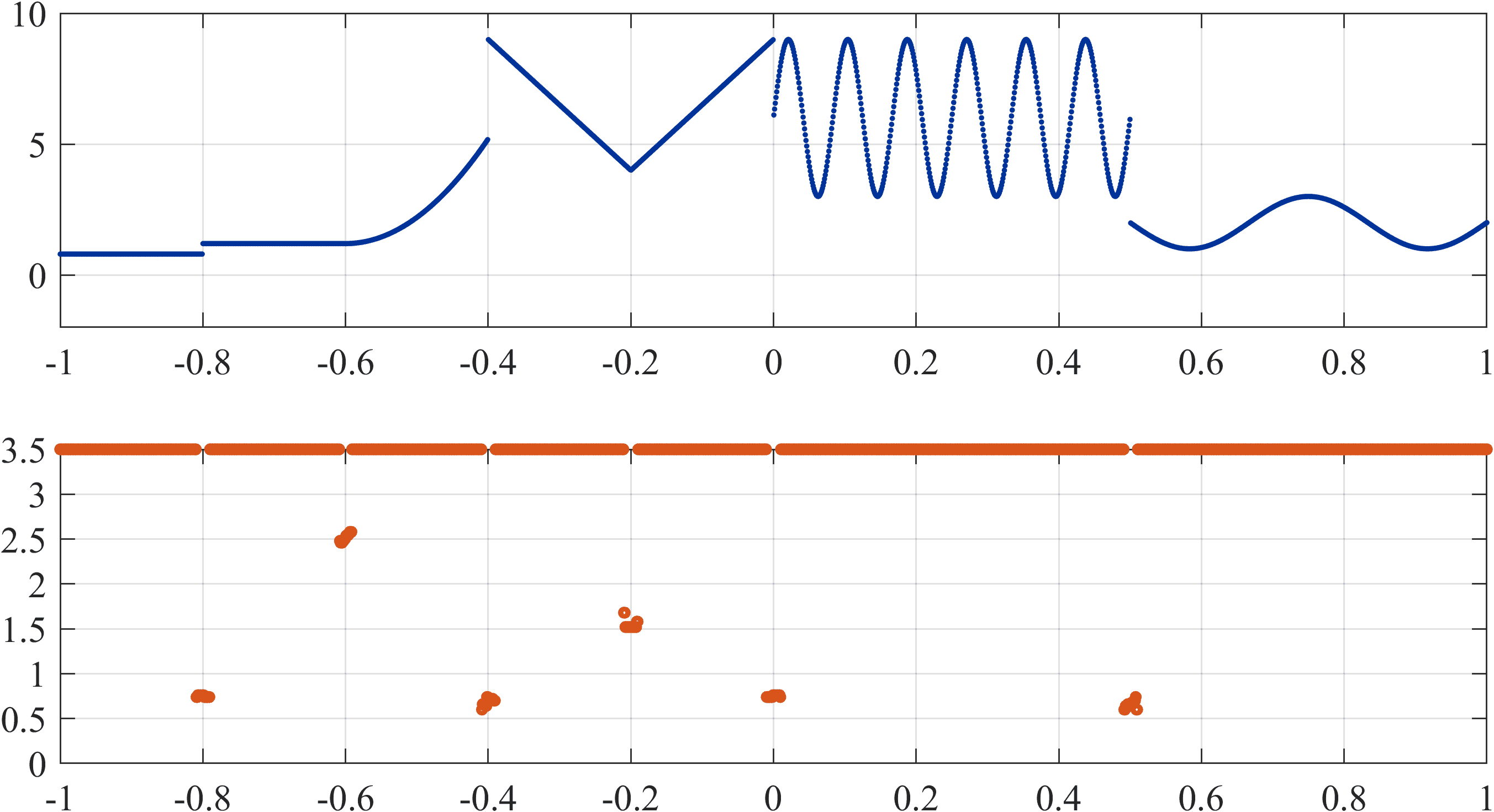}
  \put(-2,40){\scriptsize\rotatebox{90}{$f(z)$}}
  \put(-3,10){\scriptsize\rotatebox{90}{$\tilde{s}(\Omega_z)$}}
  \put(55,0){\scriptsize{$z$}}
   \put(55,29){\scriptsize{$z$}}
  \put(25,-3){\scriptsize{(b) Data-driven regularity $\tilde{s}(\Omega_z)$ }}
  \put(25,26.5){\scriptsize{(a)Piecewise-smooth test function (\Cref{eq:testfun1d}) }}
  \end{overpic}
  \caption{
  Native-norm profiles $\eta(m)$ are evaluated for
  $m\in[0.6,3.5]$ with step size $\Delta m=0.025$
  on local neighborhoods of size~20 (total 2000 samples).
  The elbow location $m^*(z)$ is extracted from the discrete
  curvature of the log-profile~\eqref{eq:m-eta-star-def}; we report $\tilde{s}(\Omega_z)$ accordingly.
  The resulting spatial map shows pronounced drops near
  nonsmooth regions, indicating reduced local Sobolev regularity.}
  \label{fig:fx_order}
\end{figure}

Figure~\ref{fig:fx_order} illustrates this procedure for a piecewise-smooth
function \eqref{eq:testfun1d} defined in Appendix~\ref{appendix:testfun1d} containing several singular features.
The recovered spatial map $m^\ast(z)$, interpreted as the data-driven regularity $\tilde s(\Omega_z)$, provides a sharp quantitative indicator of spatially varying smoothness.
In this example, the data-driven regularity $\tilde{s}(\Omega_z)$ distinguishes
different types of singular behavior:
near jump discontinuities $\tilde{s}(\Omega_z)\gtrsim 0.5$;
around corners (integrals of jumps) $\tilde{s}(\Omega_z)\gtrsim 1.5$;
and after a second integration $\tilde{s}(\Omega_z)\gtrsim 2.5$.
Although the last case ($s=2.5$) may appear visually smooth, the native-norm profile $\eta(m)$
reveals its limited local regularity.
A slight overestimation of $\tilde{s}(\Omega_z)$ is expected, since the analysis is based on
discrete samples.
A more detailed discussion of this effect is provided in
Section~\ref{sec:nativenorm-bandlimited}.


\section{Band-limited Analysis of the Sobolev-scale Norm Profile}
\label{sec:nativenorm-bandlimited}

The construction in the previous section, where the data-driven regularity
$\tilde{s}(\Omega_z)$ was introduced through the curvature of the
Sobolev-scale norm profile~$\eta(m)$, was based on the empirical observation
that $\eta(m)$ approximately follows the upper envelope~\eqref{eq:upper-bound}.
To complete the theoretical picture, we now seek a complementary
\emph{lower bound} for~$\eta(m)$ that characterizes its minimal growth
for functions of a given Sobolev regularity~$s$.
Deriving such a bound is nontrivial, since it must hold for all
$f\in H^s$, including arbitrarily smooth functions.
To this end, we introduce a band-limited surrogate formulation that
captures the essential spectral content of Sobolev functions while enabling
explicit comparison between Sobolev and native norms.
This band-limited analysis provides a quantitative lower bound on~$\eta(m)$
and forms the theoretical foundation for the main theorem presented below.

\subsection{Spectral Median of Band-Limited Functions}
\label{subsec:bandlimited-spectral-median}

We work in the Fourier-defined (Bessel-potential) Sobolev spaces.
For any $m \ge 0$ and any $g \in L^2(\mathbb{R}^d)$, define
\begin{equation}\label{eq:sobolev_fourier_def}
\|g\|_{H^m(\mathbb{R}^d)}^2
:= \int_{\mathbb{R}^d} (1 + \|\omega\|_2^2)^m \, |\widehat g(\omega)|^2 \, d\omega,
\end{equation}
which is well-defined also for fractional $m$.

Given a cutoff $\sigma > 0$, define the band-limited space
$\mathcal{B}_\sigma := \{ f_\sigma \in L^2(\mathbb{R}^d) :
\operatorname{supp} \widehat f_\sigma \subseteq B(0, \sigma) \}$.
Thus, $\mathcal{B}_\sigma$ consists of all functions whose Fourier transforms vanish
outside the ball of radius~$\sigma$ centered at the origin.
To describe how the $H^m$-norm of a band-limited function is distributed across
frequencies near the spectral cutoff, we introduce a balance parameter called the
\emph{spectral median}.

\begin{definition}[Spectral median]
Let $f_\sigma \in \mathcal{B}_\sigma$ and $t > 0$.
For $\theta \in [0,1]$, define the spectral shell
$\mathcal{S}_{\theta,\sigma} := \{ \omega \in \mathbb{R}^d :
\theta \le \|\omega\|_2/\sigma \le 1 \}$.
The \emph{spectral median} of $f_\sigma$ is defined as
\begin{equation}\label{eq:beta-halfeq}
\beta_t(f_\sigma)
:= \max \Bigl\{ \theta \in [0,1] :
\int_{\mathcal{S}_{\theta,\sigma}} (1+\|\omega\|_2^2)^t
\, |\widehat f_\sigma(\omega)|^2 \, d\omega
\ge \tfrac12 \|f_\sigma\|_{H^t(\mathbb{R}^d)}^2 \Bigr\}.
\end{equation}
\end{definition}

Therefore, the spectral median identifies the largest normalized cutoff $\beta_t\in [0,1]$ such that the
highpass region $\mathcal{S}_{\beta_t(f_\sigma),\sigma}$ still contributes at least
half of the $H^t$-norm of $f_\sigma$. The map
$\theta \mapsto \int_{\mathcal S_{\theta,\sigma}} (1+\|\omega\|_2^2)^t
|\widehat f_\sigma(\omega)|^2 \, d\omega$
is continuous and nonincreasing in $\theta$ by dominated convergence.
Hence the admissible set in~\eqref{eq:beta-halfeq}  is closed and nonempty, and the maximum is well defined.
The maximization convention selects the largest admissible cutoff; no
separate uniqueness of a balancing point is required.

The next result shows that the growth of Sobolev norms with respect to the order $m$ is bounded from below by the spectral median and a lower‑order Sobolev norm.

\begin{lemma}\label{lem:inverse-ineq}
For any $f_\sigma \in \mathcal{B}_\sigma(\mathbb{R}^d)$ and $0 < s < m$, we have
\begin{equation}\label{eq:inv-ineq-integrated}
\|f_\sigma\|_{H^m(\mathbb{R}^d)}
\;\ge\;
\|f_\sigma\|_{H^s(\mathbb{R}^d)}\,
\exp\!\biggl(
  \frac{1}{4}\int_{s}^m
  \log\!\bigl(1+(\beta_t(f_\sigma)\sigma)^2\bigr)\,dt
\biggr).
\end{equation}
\end{lemma}
\begin{proof}
If $f_\sigma \equiv 0$ the claim is trivial.
For $f_\sigma \neq 0$, fix any $t\in[s,m]$.
We begin with the logarithmic derivative of the squared Sobolev norm,
\begin{equation}\label{eq:log-deriv-quotient}
\frac{d}{dt}\,\log \|f_\sigma\|_{H^t(\mathbb{R}^d)}^2
= \|f_\sigma\|_{H^t(\mathbb{R}^d)}^{-2}
\left( \frac{d}{dt}\,\|f_\sigma\|_{H^t(\mathbb{R}^d)}^2 \right).
\end{equation}
The derivative of $\|f_\sigma\|_{H^t}^2$ with respect to $t$ is
\begin{equation}\label{eq:Ht-derivative}
\frac{d}{dt}\,\|f_\sigma\|_{H^t(\mathbb{R}^d)}^2
= \int_{\|\omega\|_2\le \sigma}
(1+\|\omega\|_2^2)^t\,|\widehat f_\sigma(\omega)|^2\,
\log(1+\|\omega\|_2^2)\,d\omega.
\end{equation}
Combining \eqref{eq:log-deriv-quotient} and \eqref{eq:Ht-derivative} gives
\[
\frac{d}{dt}\,\log \|f_\sigma\|_{H^t(\mathbb{R}^d)}^2
=
\|f_\sigma\|_{H^t(\mathbb{R}^d)}^{-2}
\int_{\|\omega\|_2\le \sigma}
(1+\|\omega\|_2^2)^t\,|\widehat f_\sigma(\omega)|^2\,
\log(1+\|\omega\|_2^2)\,d\omega.
\]

For simplicity, write $\beta_t=\beta_t(f_\sigma)$.
We restrict the integral in~\eqref{eq:Ht-derivative} to the high‑pass region
$\mathcal S_{\beta_t,\sigma}$ and
apply the defining condition~\eqref{eq:beta-halfeq} for $\beta_t$ to obtain
\[
\begin{aligned}
\frac{d}{dt}\,\log \|f_\sigma\|_{H^t(\mathbb{R}^d)}^2
&\ge
\|f_\sigma\|_{H^t(\mathbb{R}^d)}^{-2}
\int_{\mathcal{S}_{\beta_t,\sigma}}
(1+\|\omega\|_2^2)^t\,|\widehat{f}_\sigma(\omega)|^2\,
\log(1+\|\omega\|_2^2)\,d\omega \\
&\ge
\log(1+(\beta_t\sigma)^2)\,
\|f_\sigma\|_{H^t(\mathbb{R}^d)}^{-2}
\int_{\mathcal{S}_{\beta_t,\sigma}}
(1+\|\omega\|_2^2)^t\,|\widehat{f}_\sigma(\omega)|^2\,d\omega, \\
\end{aligned}
\]
where the second line uses that
$\log(1+\|\omega\|_2^2)\ge \log(1+(\beta_t\sigma)^2)$
for every $\omega\in\mathcal{S}_{\beta_t,\sigma}$.
By using the half‑energy property~\eqref{eq:beta-halfeq}, we obtain
\begin{equation}\label{eq:lower-bound-calc}
\frac{d}{dt}\,\log \|f_\sigma\|_{H^t(\mathbb{R}^d)}
\ge
\frac{1}{4}\,\log\bigl(1+(\beta_t\sigma)^2\bigr).
\end{equation}
Integrating \eqref{eq:lower-bound-calc} with respect to $t$
over the interval $[s,m]$ completes the proof.
\end{proof}

\subsection{Lower Bound for the Sobolev-Scale Norm Profile}
\label{subsec:LowerBoundSobolevScale}

We begin by recalling that the existence and basic approximation
properties of band-limited interpolants are well established
\cite[Section~3]{narcowich2006sobolev}.
Let $m,t >0$ with $d/2<t \le m$, and let
$\Omega\subset\mathbb{R}^{d}$ be a bounded Lipschitz domain equipped with a
continuous extension operator
$E_{\Omega}:H^{m}(\Omega)\to H^{m}(\mathbb{R}^{d})$.
For any $u\in H^{m}(\Omega)$ and $\sigma=\kappa_{t ,d}\,q_{X,\Omega}^{-1}$ with $\kappa_{t,d}$ chosen sufficiently large, there exists a band-limited function
$g_{\sigma,u,t,\Omega}\in \mathcal{B}_{\sigma}$
that interpolates the data on $X$, i.e.
$u|_{X}=g_{\sigma,u,t ,\Omega}|_{X}$,
and satisfies the uniform stability and approximation bounds
\begin{subequations}
\begin{align}
 \|g_{\sigma,u,t ,\Omega}\|_{H^{t }(\mathbb{R}^{d})}   
&\le C\,
   \|u\|_{H^{t }(\Omega)},
   \label{eq:BandlimitedProp_2}\\[2mm]
 \|u-g_{\sigma,u,t ,\Omega}\|_{H^{t }(\Omega)}  
&\le C'\,
   q_{X,\Omega}^{\,m-t }\,
   \|u\|_{H^{m}(\Omega)}.
   \label{eq:BandlimitedProp_2b}
\end{align}
\end{subequations}
For each $t\in[s,m]$, we introduce a canonical representative
chosen to minimize the $H^{t}$ norm subject to the interpolation
constraints:
\begin{equation}\label{eq:fstar}
  f^*_{\sigma,t,X,Y}
  := \operatorname*{arg\,min}
     \Bigl\{
        \|v\|_{H^t(\mathbb R^d)} :
        v \in \mathcal B_\sigma(\mathbb R^d),\
        v|_X = Y
     \Bigr\}.
\end{equation}
The admissible set $ \{v\in  \mathcal B_\sigma(\mathbb R^d): v|_X=Y\}$
is nonempty by the band-limited interpolation result above and is a closed
affine subspace, since point evaluation is continuous on
$\mathcal B_\sigma(\mathbb R^d)$. Hence the Hilbert projection theorem gives
a unique minimizer $f^*_{\sigma,t,X,Y}$.
As $t$ varies, the interpolation constraints remain fixed, while the Sobolev
weight $(1+\|\omega\|_2^2)^t$ changes. Figure~4 illustrates the corresponding
minimum-$H^t$ interpolants $f^*_{\sigma,t,X,Y}$.

\begin{figure}[t]
  \centering
  \subfloat[$\sigma=\pi/h$]{
    \begin{overpic}[width=60mm,height=43mm]{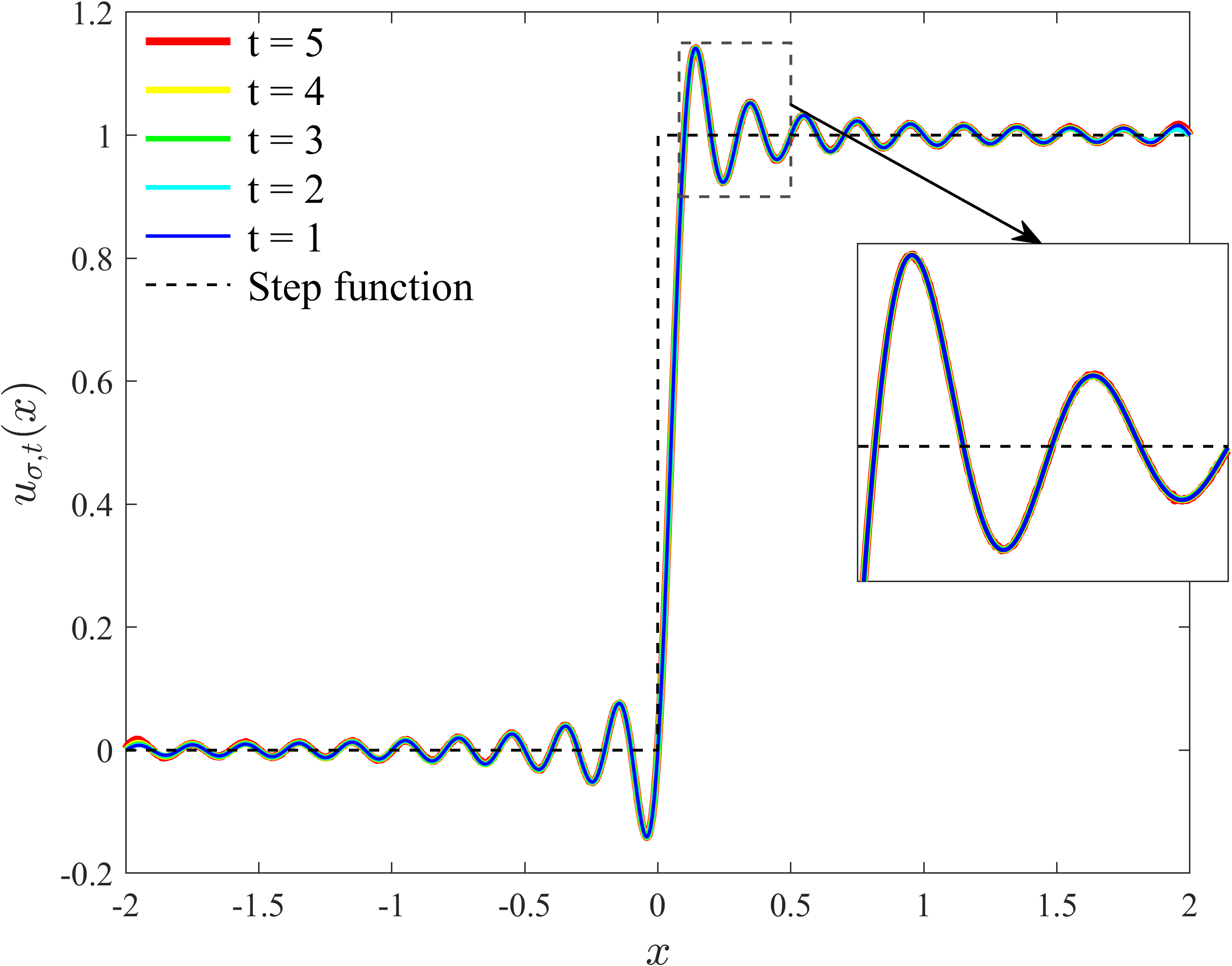}\label{fig1}\end{overpic}}
    \hspace{2mm}
    \subfloat[$\sigma=2\pi/h$]{
    \begin{overpic}[width=60mm,height=43mm]{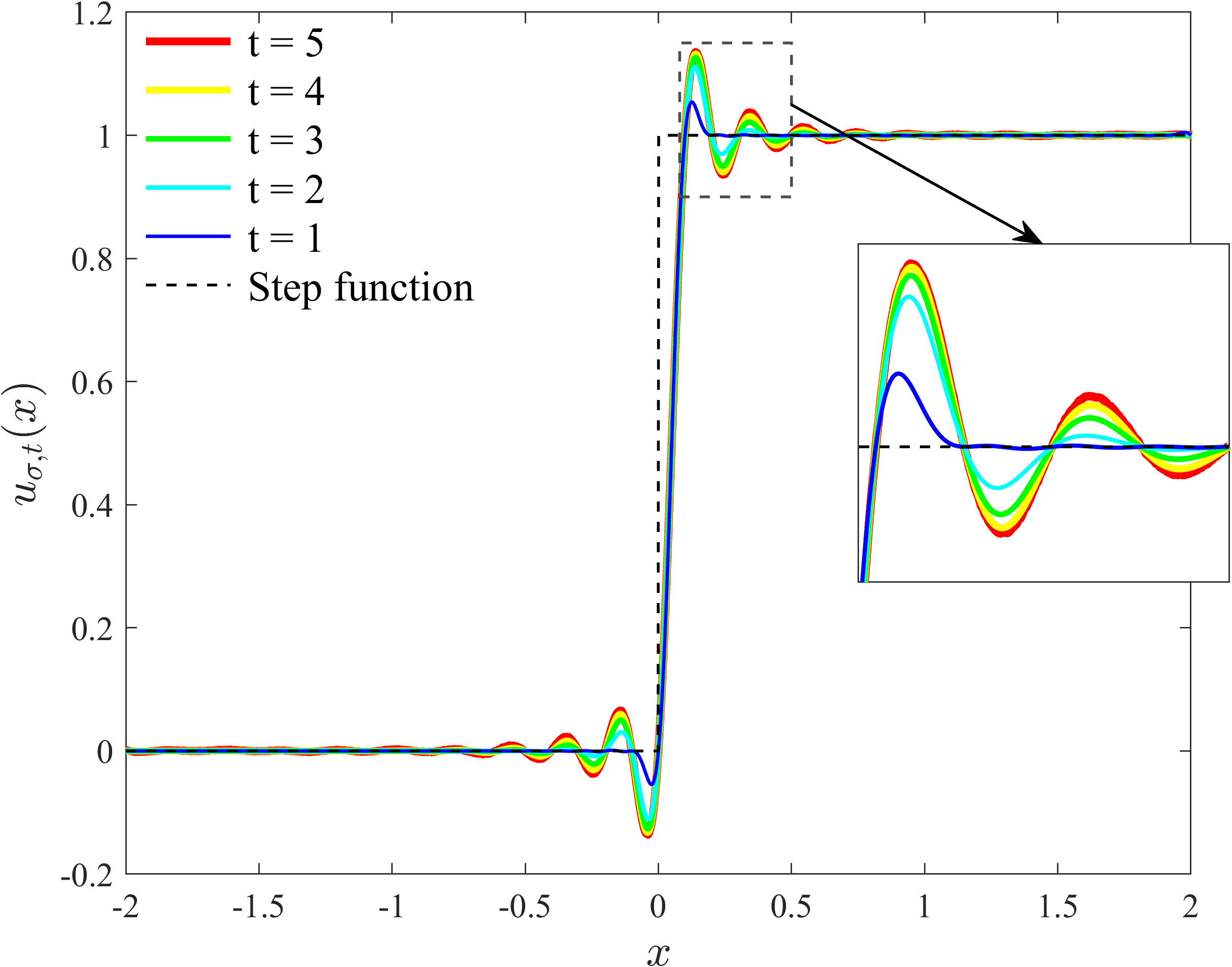}\label{fig2}
    \end{overpic}
    }
\caption{
  Minimum-$H^t$ band-limited interpolants $f^*_{\sigma,t,X,Y}$ defined in \eqref{eq:fstar} of the one-dimensional step function
  $f(x)=\chi_{(0,\infty)}(x)$ from uniformly spaced samples.
    They are realized here through the
  Fourier-domain kernels
  $\widehat{\Phi}_{\sigma,t}(\omega)=(1+\omega^2)^{-t}\,\mathbf{1}_{|\omega|\le\sigma}$,
  with $t=1,\ldots,5$.
  The Fourier-truncated kernel is used here solely to illustrate the
band-limited surrogate viewpoint and the effect of the $H^t$ weighting on the
minimal-norm interpolant.}
  \label{fig:fourier-step}
\end{figure}


We now relate the kernel interpolant $u_m\in \mathcal U_{X,\Phi_m}$ to the band-limited
interpolants introduced above.
 For each fixed
$m\in[m_{\min},m_{\max}]$, applying \eqref{eq:BandlimitedProp_2} with
$u=u_m$ and $t=m$, and choosing
$\sigma=\kappa_{m,d}q_{X}^{-1}$ sufficiently large, we obtain a
band-limited function
$g_{\sigma,u_m,m,\Omega}\in\mathcal B_\sigma(\mathbb R^d)$ such that 
    $g_{\sigma,u_m,m,\Omega}|_X=u_m|_X=Y$
and
\begin{equation}\label{eq:BL-prop-reuse}
    \|g_{\sigma,u_m,m,\Omega}\|_{H^m(\mathbb R^d)}
    \le C_m\|u_m\|_{H^m(\Omega)} .
\end{equation}
Over the fixed compact range $[m_{\min},m_{\max}]$, we assume the constants
in the band-limited interpolation estimate can be chosen locally uniformly in
$m$. We therefore take finite constants
\[
\kappa_d:=\sup_{m\in[m_{\min},m_{\max}]}\kappa_{m,d}<\infty,\qquad
C:=\sup_{m\in[m_{\min},m_{\max}]} C_m<\infty,
\]
such that \eqref{eq:BL-prop-reuse} holds uniformly for all $m$ with the common bandwidth $\sigma=\kappa_d\,q_{X,\Omega}^{-1}$.
Next, we note that the equivalence constant between the native and Sobolev norms for
Sobolev reproducing kernels is independent of~$m$.

\begin{lemma}\label{lem:kernel-native-sobolev-equiv}
Assume $\Phi_m$ ($m>d/2$) satisfies the Fourier decay condition~\eqref{eq:fourier-decay} with upper-bound constant~$c_2$.
Then there exists a constant $C_{d,c_2}>0$, depending only on~$d$, such that
\begin{equation}\label{eq:native-dominates-sobolev}
  \|v\|_{\mathcal N_{\Phi_m}}
  \;\ge\;
  C_{d,c_2}\,\|v\|_{H^m(\mathbb R^d)},
  \qquad
  \text{for all }v\in\mathcal N_{\Phi_m}(\mathbb R^d).
\end{equation}
\end{lemma}
\begin{proof}
The claim follows directly from the Fourier decay assumption~\eqref{eq:fourier-decay}
and the standard Fourier characterization of native spaces
(see, e.g.,~\cite[Cor.~10.48]{wendland2004scattered}).
\end{proof}

Combining~\eqref{eq:BL-prop-reuse} and
\eqref{eq:native-dominates-sobolev} gives
\begin{equation}\label{eq:eta-to-fsigma-ut}
  \eta(m)
  \;:=\;
  \|u_m\|_{\mathcal N_{\Phi_m}}
  \;\gtrsim\;
  \|u_m\|_{H^{m}(\Omega)}
  \;\gtrsim\;
  \|g_{\sigma,u_m,m,\Omega}\|_{H^{m}(\mathbb R^d)}.
\end{equation}
By definition of the optimal band-limited interpolant
$f^*_{\sigma,m,X,Y}$ in~\eqref{eq:fstar}, its $H^m$-norm is minimal
among all band-limited functions matching the data; hence, there exists a constant $C_0$ such that the following holds:
\begin{equation}\label{eq:eta-controls-fstar-short}
  \eta(m)
  \;\ge\;
  C_0\,\|f^*_{\sigma,m,X,Y}\|_{H^{m}(\mathbb R^d)},
  \quad\text{for all $m\in[m_{\min},m_{\max}]$}.
\end{equation}
Thus, the kernel-based norm profile $\eta(m)$ is bounded below by the
family of minimal-norm band-limited interpolants
$f^*_{\sigma,m,X,Y}$.
We are now ready to apply
Lemma~\ref{lem:inverse-ineq} to $f^*_{\sigma,m,X,Y}$
and establish the quantitative lower bound leading to the main theorem.

\begin{theorem}
\label{thm:native-lower}
Let $\{\Phi_m\}_{m\in[m_{\min},m_{\max}]}$ be a family of admissible kernels on $\mathbb{R}^d$
such that, for each $m>d/2$, the associated native space $\mathcal N_{\Phi_m}$
is norm-equivalent to the Sobolev space $H^m(\mathbb{R}^d)$.
For each $m\in[m_{\min},m_{\max}]$, let $u_m\in\mathcal U_{X,\Phi_m}$
be the kernel interpolant satisfying $u_m|_X=Y$, and define
$\eta(m) := \|u_m\|_{\mathcal N_{\Phi_m}}$.
For any reference Sobolev order $s>0$, set $s_\ast := \min\{s,m\}$.
For $\sigma = \kappa\,q_X^{-1}$ with some sufficiently large fixed $\kappa>0$,
let $f^*_{\sigma,t,X,Y}$ denote the minimal-norm band-limited interpolants
defined in~\eqref{eq:fstar}.
Then there exists a constant $C>0$, independent of $X$, $Y$, and
$m\in[m_{\min},m_{\max}]$, such that
\begin{equation}\label{eq:native-norm-lower-final}
   \eta(m)
   \;\ge\;
   C\,
   \|f^*_{\sigma,s_\ast,X,Y}\|_{H^{s_\ast}(\mathbb{R}^d)}\,
   \exp\!\Biggl(
     \frac{1}{4}\int_{s_\ast}^{m}
     \log\!\Bigl(1+\bigl(\beta_t(f^*_{\sigma,m,X,Y})\,\sigma\bigr)^2\Bigr)\,dt
   \Biggr),
\end{equation}
where, for each $t\in[s_\ast,m]$, $\beta_t(\cdot)$ denotes the spectral median
as defined in~\eqref{eq:beta-halfeq}.
\end{theorem}
\begin{remark}

The bound in Theorem~\ref{thm:native-lower} is parametrized by the reference order $s$. It should therefore be understood as a smoothness-misspecification estimate for the norm profile, not as an exact recovery statement for $\tilde s(\Omega_z)$. Suppose that the local data have Sobolev regularity $s(\Omega_z)$ and that $s=s(\Omega_z)$ is used only for interpretation. Then the exponential factor in \eqref{eq:native-norm-lower-final} is absent for
$m\le s(\Omega_z)$ and becomes active for $m>s(\Omega_z)$.
 Moreover, since $\sigma=\kappa q_X^{-1}$, the post-transition growth is tied to the local sampling scale. For quasi-uniform stencils $q_X\asymp h$. If there exist an interval $I\subset (s(\Omega_z),m_{\max}]$ and a constant $\beta_0>0$ such that
    $\beta_t(f^*_{\sigma,m,X,Y})\ge \beta_0,
    \ t\in I,$
then
\[
    \int_I
    \log\!\left(
        1+\bigl(\beta_t(f^*_{\sigma,m,X,Y})\sigma\bigr)^2
    \right)\,dt
    \ge
    |I|\log\!\left(1+\beta_0^2\kappa^2 q_X^{-2}\right)
    \asymp |I|\log(h^{-1})
    \quad \text{as } h\to0 .
\]
Thus, under this condition, the post-threshold lower bound becomes steeper
as the stencil is refined. The refinement behavior of the numerical estimator
is examined in Section~5. The numerical elbow criterion in Definition~2.1 uses this increasing contrast to identify the transition region of the computed profile. A similar rate viewpoint appears in smoothness misspecification for Sobolev
kernel Gaussian processes, where oversmoothing causes a scale quantity to grow
at a rate determined by the gap between the imposed and true smoothness orders
\cite{karvonen2025scale}.
\end{remark}



\begin{proof}
If $m\le s$, then $s_\ast=m$, and the exponential term in
\eqref{eq:native-norm-lower-final} reduces to~$1$.
In this case, the estimate follows directly from
the basic inequality~\eqref{eq:eta-controls-fstar-short}, up to a constant factor.
If $m>s$, we must compare the kernel norm at order~$m$ with the
$H^{s}$-norm of the band-limited interpolant.
Applying Lemma~\ref{lem:inverse-ineq} to
$f_\sigma=f^*_{\sigma,m,X,Y}$ with lower order $s_\ast=s$ gives a lower bound
involving $\|f^*_{\sigma,m,X,Y}\|_{H^s(\mathbb R^d)}$.
Since $f^*_{\sigma,s,X,Y}$ minimizes the $H^s(\mathbb R^d)$-norm over the
same admissible set of band-limited interpolants satisfying $v|_X=Y$, 
$\|f^*_{\sigma,m,X,Y}\|_{H^s(\mathbb R^d)}
\ge
\|f^*_{\sigma,s,X,Y}\|_{H^s(\mathbb R^d)} $.
Combining this inequality with~\eqref{eq:eta-controls-fstar-short}
gives~\eqref{eq:native-norm-lower-final}.
Thus, $s_\ast=\min\{s,m\}$ simply ensures that the bound stops
at the smaller of the available regularities:
if the data-generating function is smoother ($s>m$), the estimate
uses its $H^m$-norm; otherwise, it uses its $H^s$-norm.
\end{proof}


\section{Practical Refinements: Stencil Shifting and Secant-Based Screening}
\label{sec:regularity_aware_algorithms}

This section presents two independent, practical extensions of the framework:
(i) a \emph{stencil-shift refinement} that tests a small set of nearby stencil
positions to improve spatial localization of low-regularity features, and
(ii) a local-and-global accelerated norm-sweep that estimates $\eta(m)$ on large datasets using a secant-based tail for each neighborhood, optionally combined with a one-point global screen for further speedup.

\subsection{Stencil-Shift for Sharper Localization}
\label{sec:stencil_shift_refinement}

The data-driven regularity $\tilde{s}(\Omega_z)$ in Definition~\ref{def:data-driven-regularity}
is computed on a neighborhood $\Omega_z$ centered at an evaluation point $z$.
When $\Omega_z$ crosses an interface or a steep transition,
its samples may contain both smooth and nonsmooth components.
In such cases, the data-driven regularity $\tilde{s}(\Omega_z)<m_{\max}$
tends to reflect the rougher side,
even when $z$ itself lies within a smoother subregion.
This can lead to spatial smearing, where low-regularity estimates
extend across nearby points whose neighborhoods intersect the same feature.

\begin{figure}[t]
  \centering
  \subfloat[One-dimensional example of stencil shifting.]{
    \begin{overpic}[width=60mm,height=43mm]{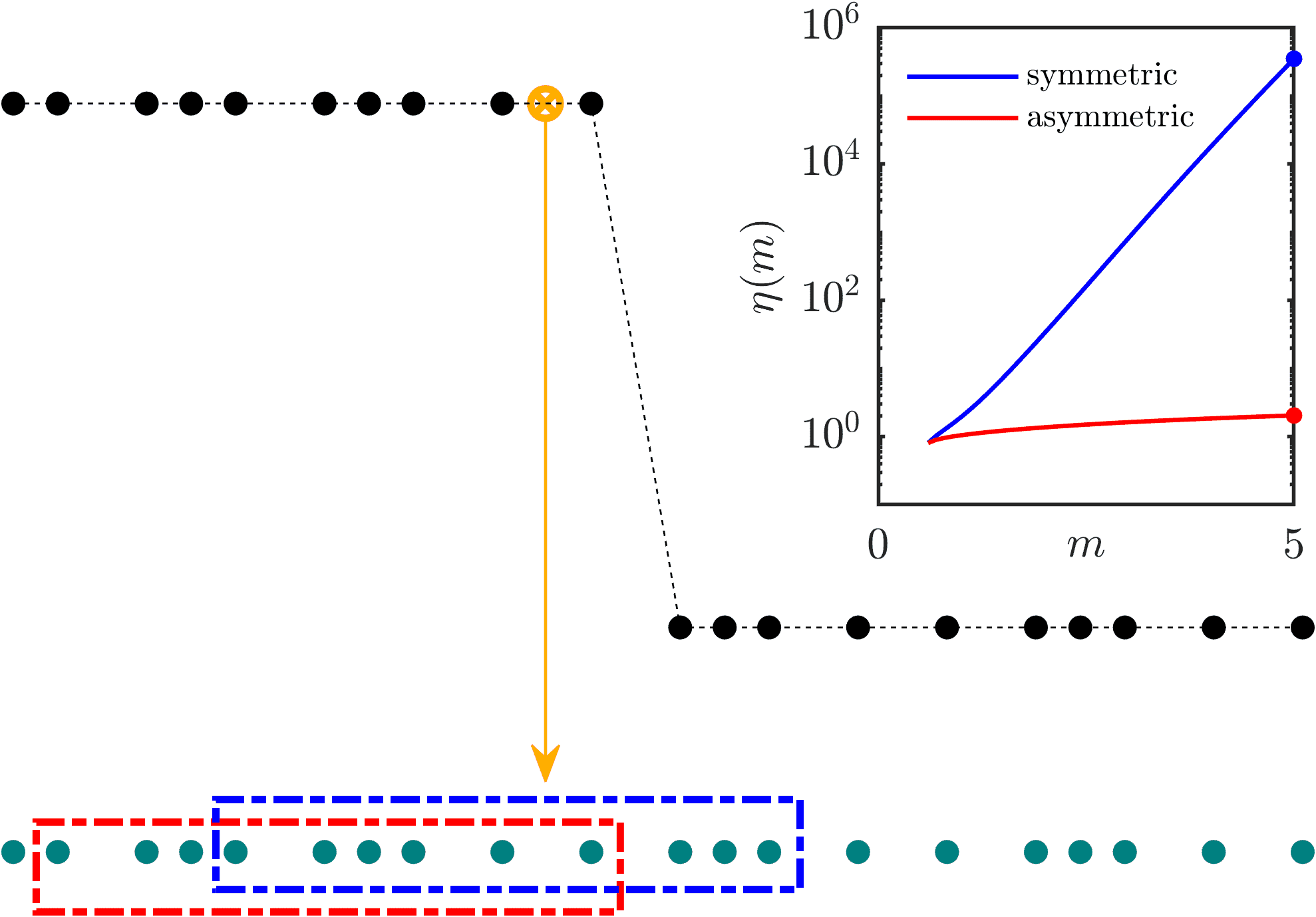}\label{fig:singularity_aware_stencil_selection_1d_a}\end{overpic}}
    \hspace{2mm}
    \subfloat[Stencil shifting near a $C^1$-but-not-$C^2$.]{
    \begin{overpic}[width=60mm,height=43mm]{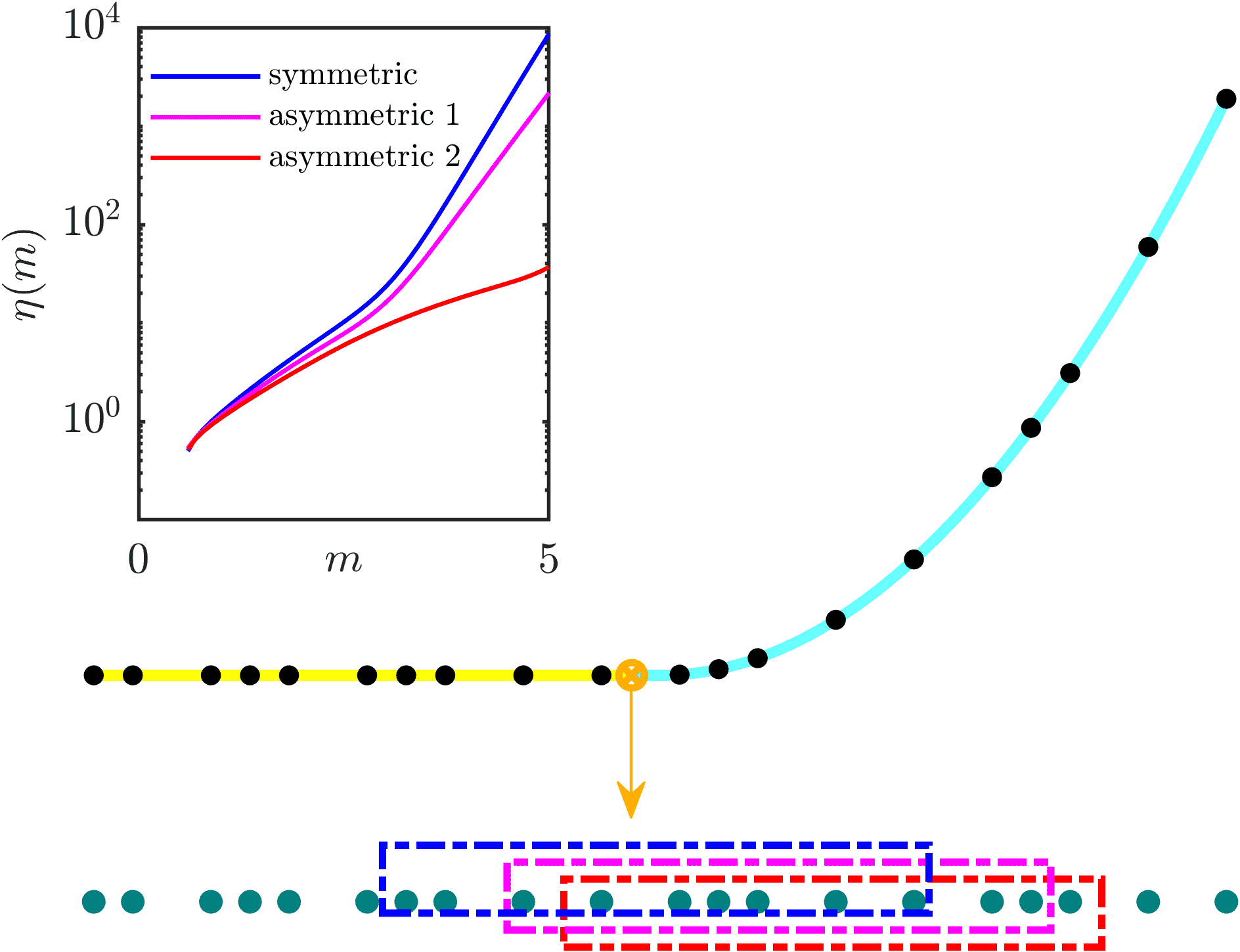}\label{fig:singularity_aware_stencil_selection_1d_b}
     \put(60,15){\scriptsize{$f(x)=\tfrac{1}{2}x^{2}\mathbf{1}_{{x\ge 0}}$}}
     \put(45,25){\scriptsize{$z = 0$}}
    \end{overpic}
    }
    \caption{(a)
A symmetric stencil spanning a discontinuity produces an underestimated data-driven regularity $\tilde{s}(\Omega_z)$.
Shifting the stencil slightly away from the jump avoids mixed samples,
leading to a higher $\tilde{s}(\Omega_z)$
and sharper localization of the smooth region. (b) The test function is $f(x)=\tfrac{1}{2}x^{2}\mathbf{1}_{{x\ge 0}}$ and the target location is $z = 0$, taken as an off-sample point ($z\notin X$).
Different shifted neighborhoods yield markedly different norm profiles $\eta(m)$. In the extreme case where the shifted stencil becomes highly one-sided and the center approaches the stencil boundary (red), the profile can appear artificially mild, leading to an overestimation of the local regularity $\tilde{s}$ under Definition~\ref{def:data-driven-regularity}. This motivates the interior condition used in the stencil-shift procedure.
}
  \label{fig:singularity_aware_stencil_selection_1d}
\end{figure}

To alleviate this effect, we test several nearby stencil positions
and retain the one that yields a higher regularity estimate.
Whenever $\tilde{s}(\Omega_z)<m_{\max}$,
the neighborhood is considered contaminated by low-regularity samples.
A small collection of alternative neighborhoods
$\{\Omega'_z\}$ is then created by locally shifting the stencil center
while preserving the interpolatory geometry.
For each candidate $\Omega'_z$,
the norm profile $\eta(m)$ is recomputed, and the corresponding
data-driven regularity $\tilde{s}(\Omega'_z)$ is evaluated.
The final regularity value at $z$ is updated as
\[
\tilde{s}(\Omega_z)
\;:=\;
\max_{\Omega'_z\in\mathfrak {N}(z)}
\tilde{s}(\Omega'_z),
\]
where $\mathfrak {N}(z)$ denotes the set of admissible shifted neighborhoods.
Figure~\ref{fig:singularity_aware_stencil_selection_1d}(a) illustrates the idea in one dimension:
a symmetric stencil spanning a discontinuity corresponds to case 2 in Definition \ref{def:data-driven-regularity}, whereas a slight lateral shift avoids mixed samples and produces a flatter $\log\eta(m)$ profile with a smaller slope, corresponding to case 3.

\begin{figure}[t]
  \centering
  \includegraphics[width=0.27\textwidth]{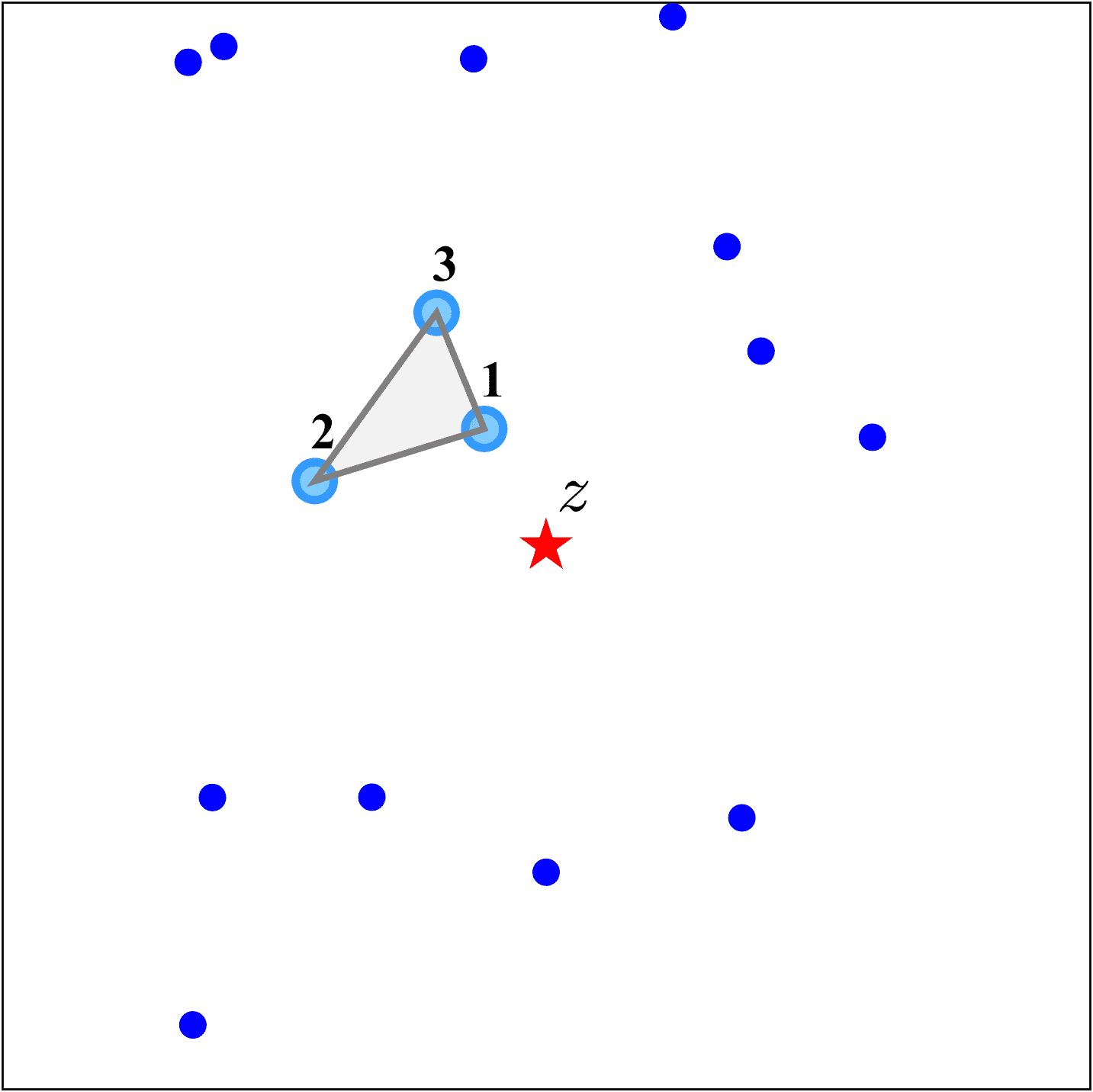}\hspace{3mm}%
  \includegraphics[width=0.27\textwidth]{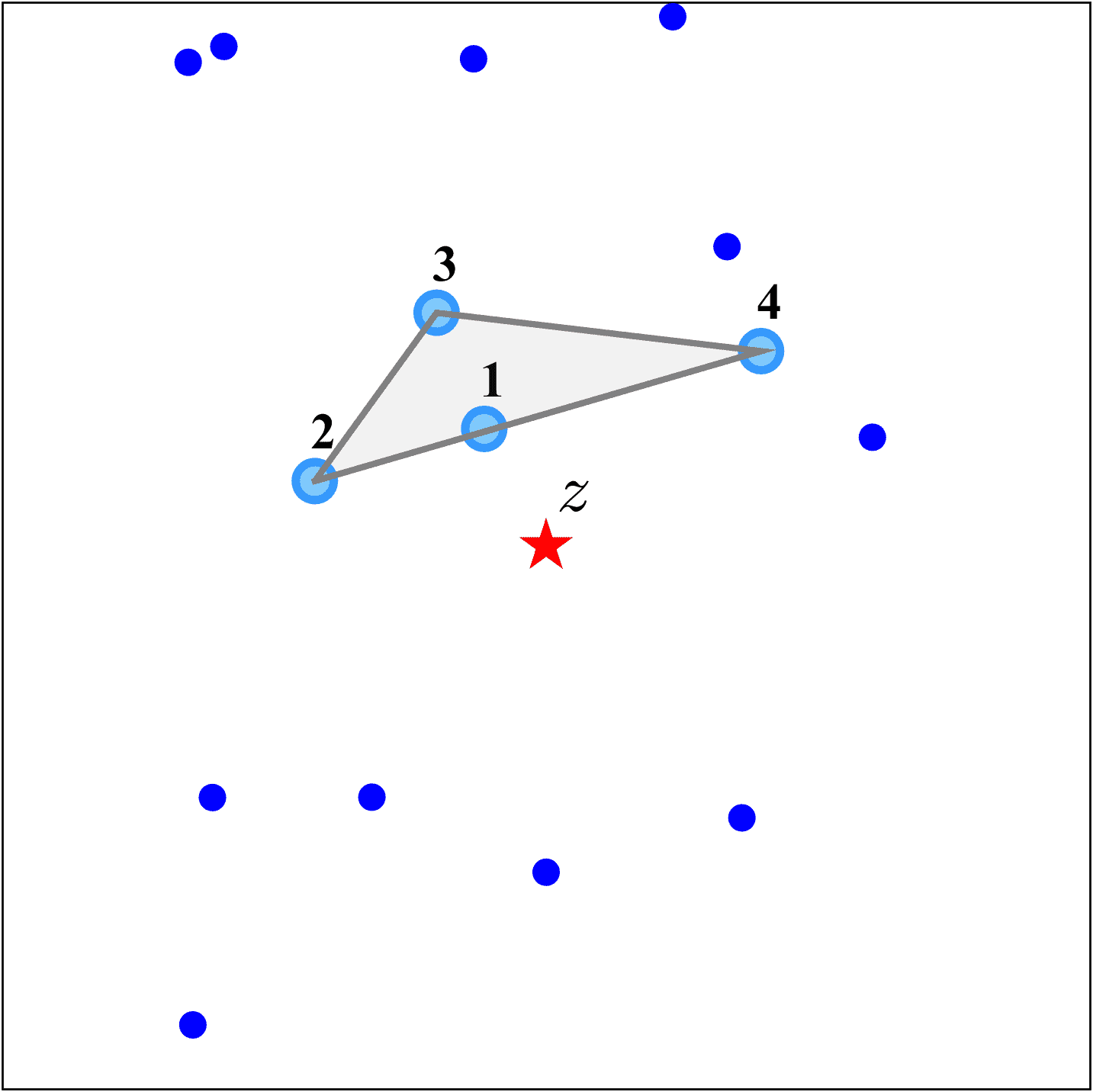}\hspace{3mm}%
  \includegraphics[width=0.27\textwidth]{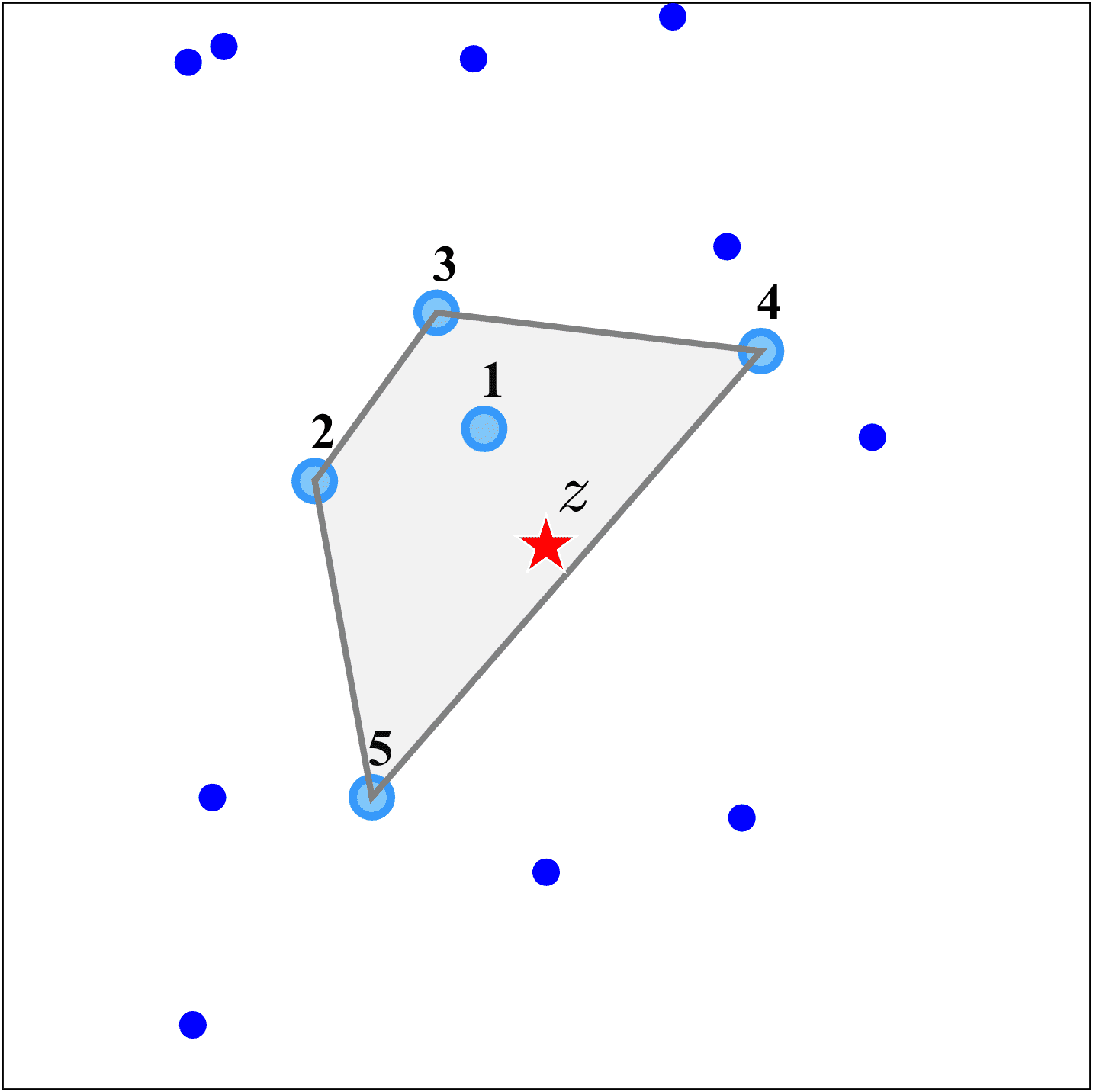}
  \caption{
Convex-hull construction around the target $z$. The target is shown by a
distinct marker, and the polygon is $\operatorname{conv}(C_z)$. The remaining
points are nearby stencil candidates and do not belong to $C_z$.
The admissibility condition
$\operatorname{dist}(z,\partial\operatorname{conv}(\mathcal C_z))
\ge c_z q_{X_z}$
prevents shifted centers from approaching the hull boundary.
  }
  \label{fig:scatter_hull}
\end{figure}

\begin{figure}[t]
  \centering
\subfloat[Candidate shift directions]{
\begin{overpic}[width=0.27\textwidth]{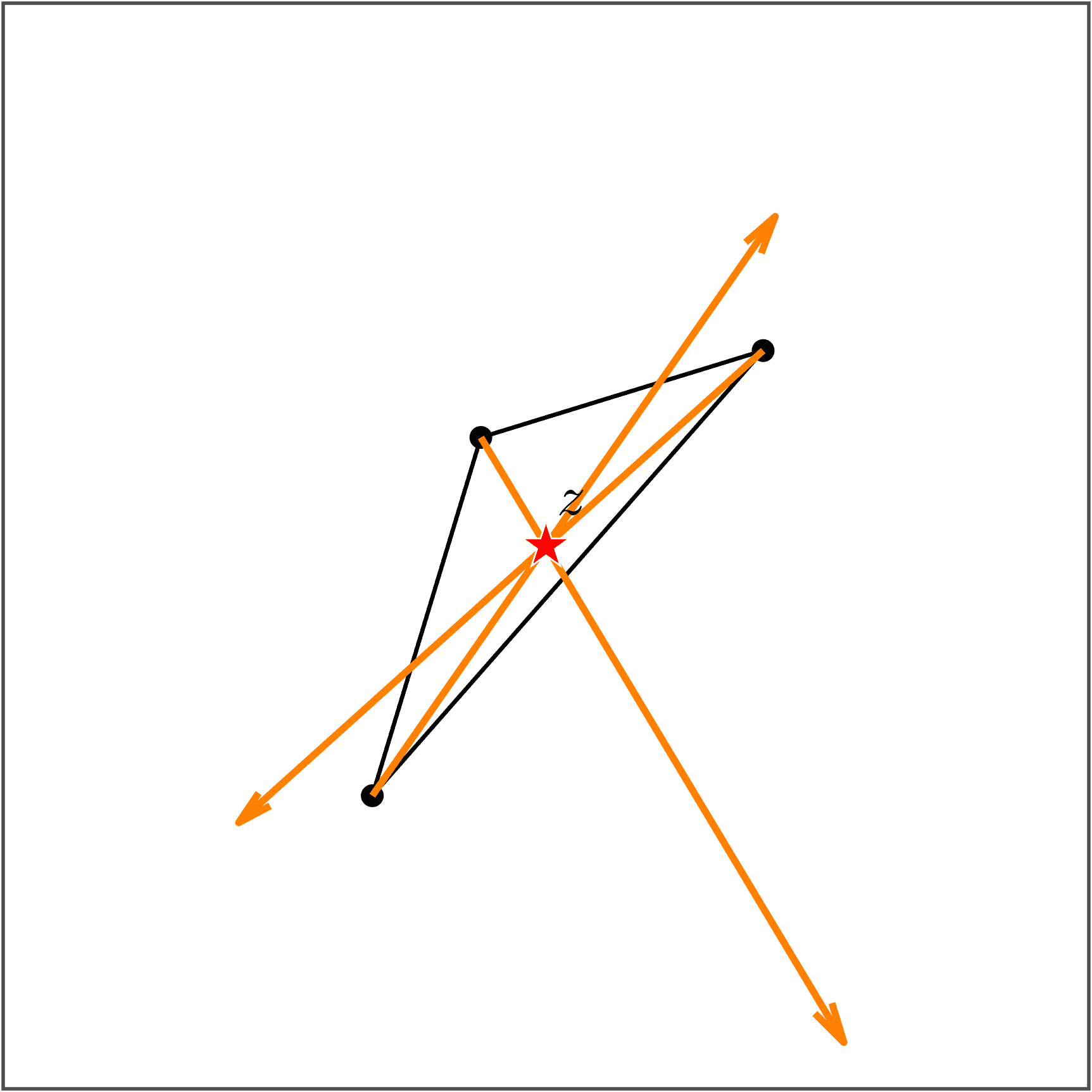}
\end{overpic}}
\hspace{1mm}
\subfloat[Shifted stencils  along \\ rays from $z$ to hull vertices]{
\begin{overpic}[width=0.27\textwidth]{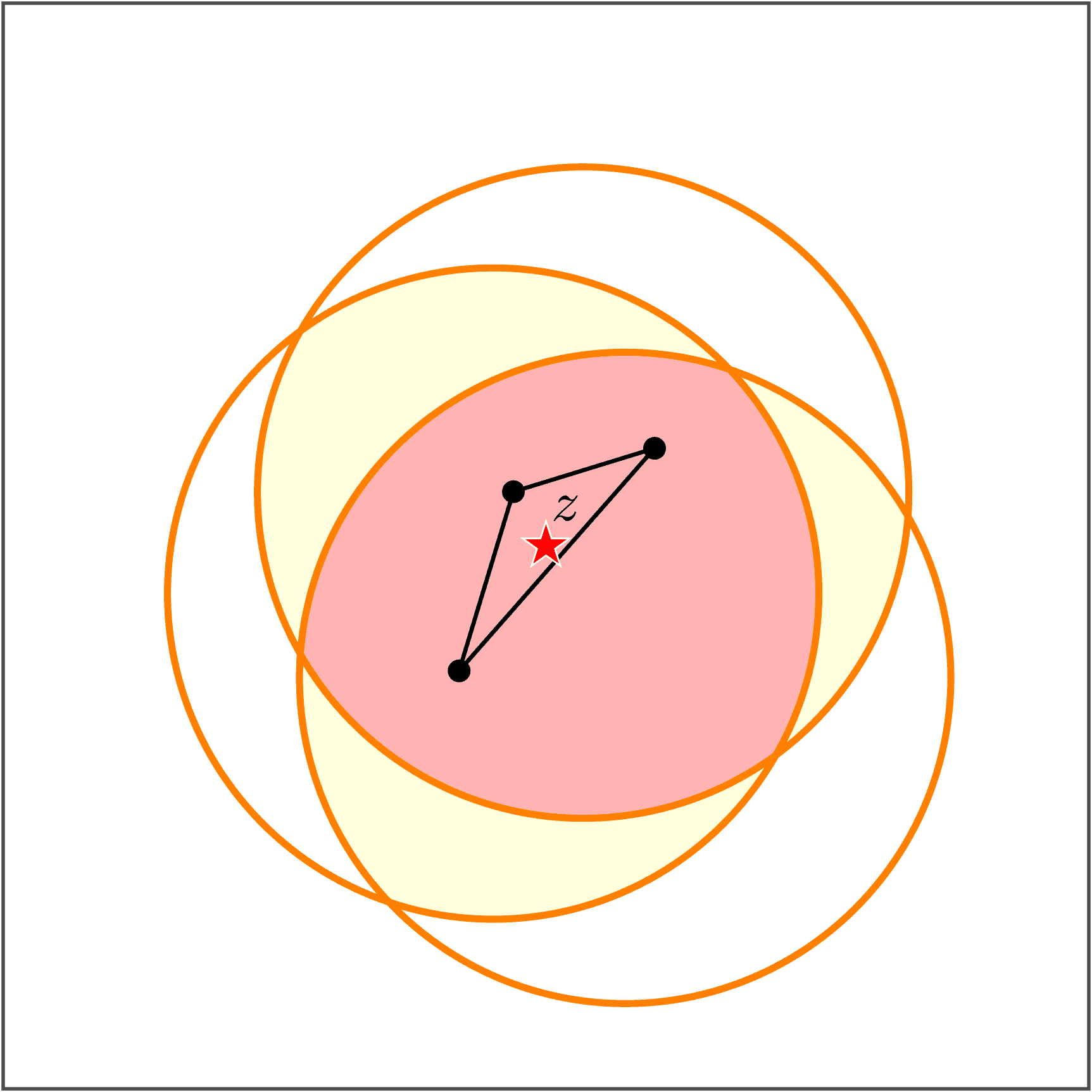}
\end{overpic}}
\hspace{1mm}
\subfloat[Shifted stencils along \\ outward hull-edge normals]{
\begin{overpic}[width=0.27\textwidth]{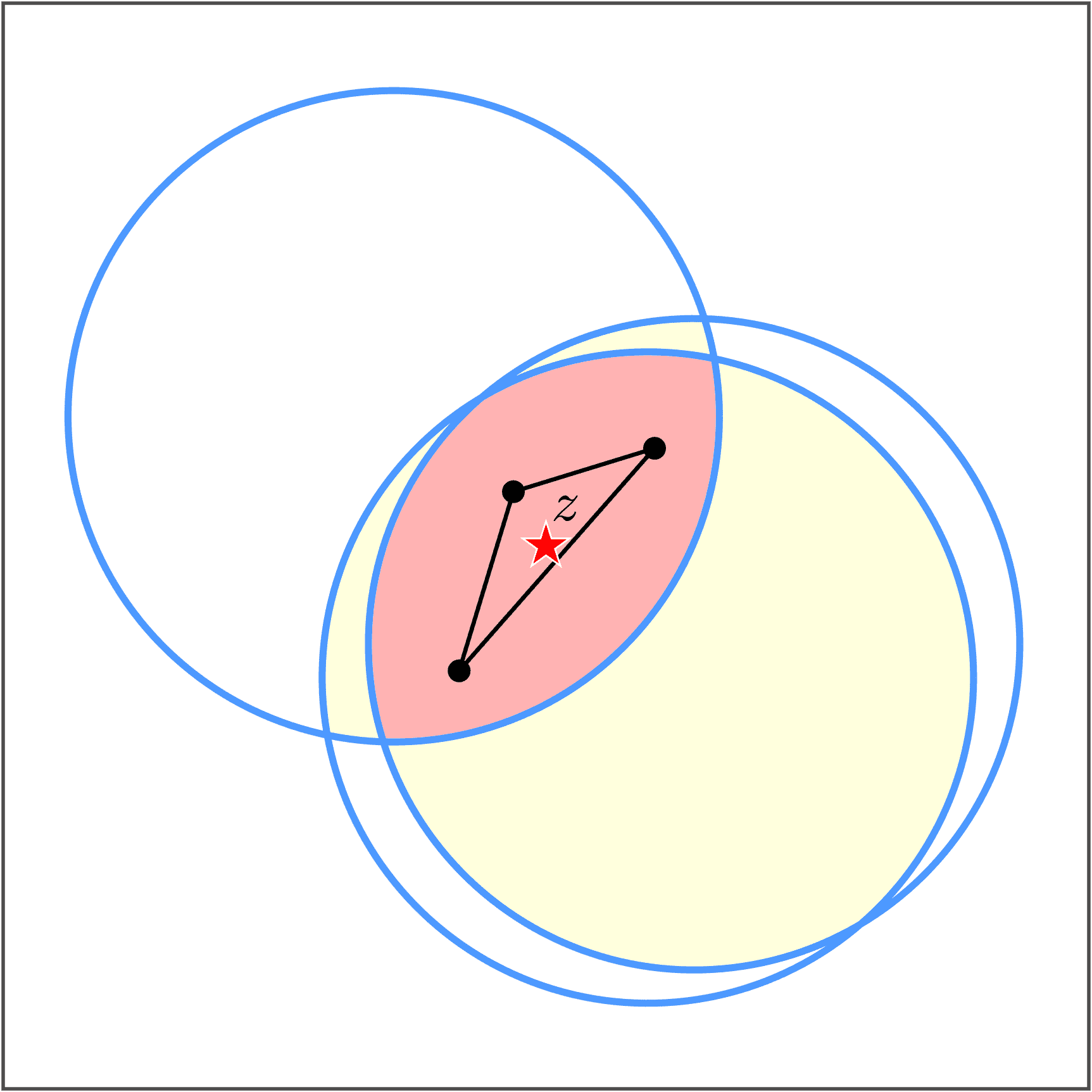}
\end{overpic}}
  \caption{
Candidate shift directions in two dimensions,
generated from the convex-hull geometry.
Each trial offset defines a circular region used to form a shifted neighborhood.
  }
  \label{fig:stencil_directions}
\end{figure}

In higher dimensions, the candidate shifts are constructed geometrically.
Starting from the symmetric stencil $X_z^{\mathrm{sym}}$,
we identify a stable core $\mathcal{C}_z \subseteq X_z^{\mathrm{sym}}$
and its convex hull
$\operatorname{conv}(\mathcal{C}_z)$ around $z$
(Figure~\ref{fig:scatter_hull}), which defines the region where interpolation remains well conditioned.
Candidate shift directions are drawn from
(i) rays connecting $z$ to the hull vertices and
(ii) outward normals to hull edges in two dimensions
or to hull faces in three dimensions
(Figure~\ref{fig:stencil_directions}).

For each direction $v_k$, we construct a few trial offsets
$z+r\,v_k$, $r\in[0,r_{\max}(z)]$.
Each offset defines a circular (in 2D) or spherical (in 3D)
region $\mathcal{B}(z+r\,v_k,\,\rho_z)$ of radius $\rho_z$,
matching the local stencil radius of $\Omega_z$.
New neighborhoods $\Omega'_z(r)$ are obtained by selecting the $n$
nearest samples to the shifted center within $\mathcal{B}(z+r\,v_k,\,\rho_z)$.
A key caveat is that overly aggressive shifts may produce stencils that become
geometrically one-sided (or otherwise under-informative) with respect to a nearby
limited-regularity feature. In this situation, the shifted neighborhood may contain
too few samples across the feature to reliably reflect the local loss of regularity.
Consequently, the tail of the $\log\eta(m)$ profile can appear artificially mild,
leading to an overestimation of the local regularity under
Definition~\ref{def:data-driven-regularity}.
To mitigate this effect, each candidate should satisfy the interior condition
\[
\operatorname{dist}\!\bigl(z,\partial\operatorname{conv}(\mathcal{C}_z)\bigr)\ge c_z\,q_{X_z},
\]
where $q_{X_z}$ is the local spacing and $c_z$ is an adaptive safety
factor. The role of $c_z$ is geometric: it keeps the evaluation point at least
$c_z$ local spacings away from the boundary of the stable core. 
In the numerical experiments,
we take $c_z=1$ for first-order-type features and $c_z=2$ for second-order-type
features, with $c_z$ capped at $2$ to avoid rejecting too many admissible
shifts.  Figure~\ref{fig:singularity_aware_stencil_selection_1d}(b) provides an example of this failure mode for a function that is locally $C^1$ but not $C^2$.
Here the target location is $z=0$, which is an off-sample point (i.e., $z\notin X$).
Depending on the shift, the norm profile may either indicate limited regularity or lead to an overestimated regularity under
Definition~\ref{def:data-driven-regularity}.
If the target point were included as a data site (i.e., $z\in X$), this overestimation is not observed in this demo.

The stencil-shift procedure uses only geometric information
from the local sample distribution and operates independently at each $z$.
It suppresses spurious low-regularity detections near interfaces,
narrows smeared transition zones,
and improves the spatial localization of the data-driven regularity $\tilde{s}(\Omega_z)$.
In addition, it can also serve as a geometric prefilter for constructing
RBF-FD or other meshfree stencils that preserve interpolation accuracy
while avoiding cross-feature sampling.

\subsection{Secant-Based Low-Regularity Screening for Large Datasets}
\label{sec:norm_sweep_acceleration}
A dense evaluation of the Sobolev-scale native-norm profile $\eta(m)$ at many $m$-values for every point $z$
quickly becomes intractable on large datasets.  We therefore use the high-order secant of $\log\eta(m)$ as a
low-cost screening indicator. This quantity is not intended to replace the
curvature-based elbow estimator in Definition~2.1.

\begin{figure}[t]
\centering
\begin{overpic}[width=0.40\textwidth]{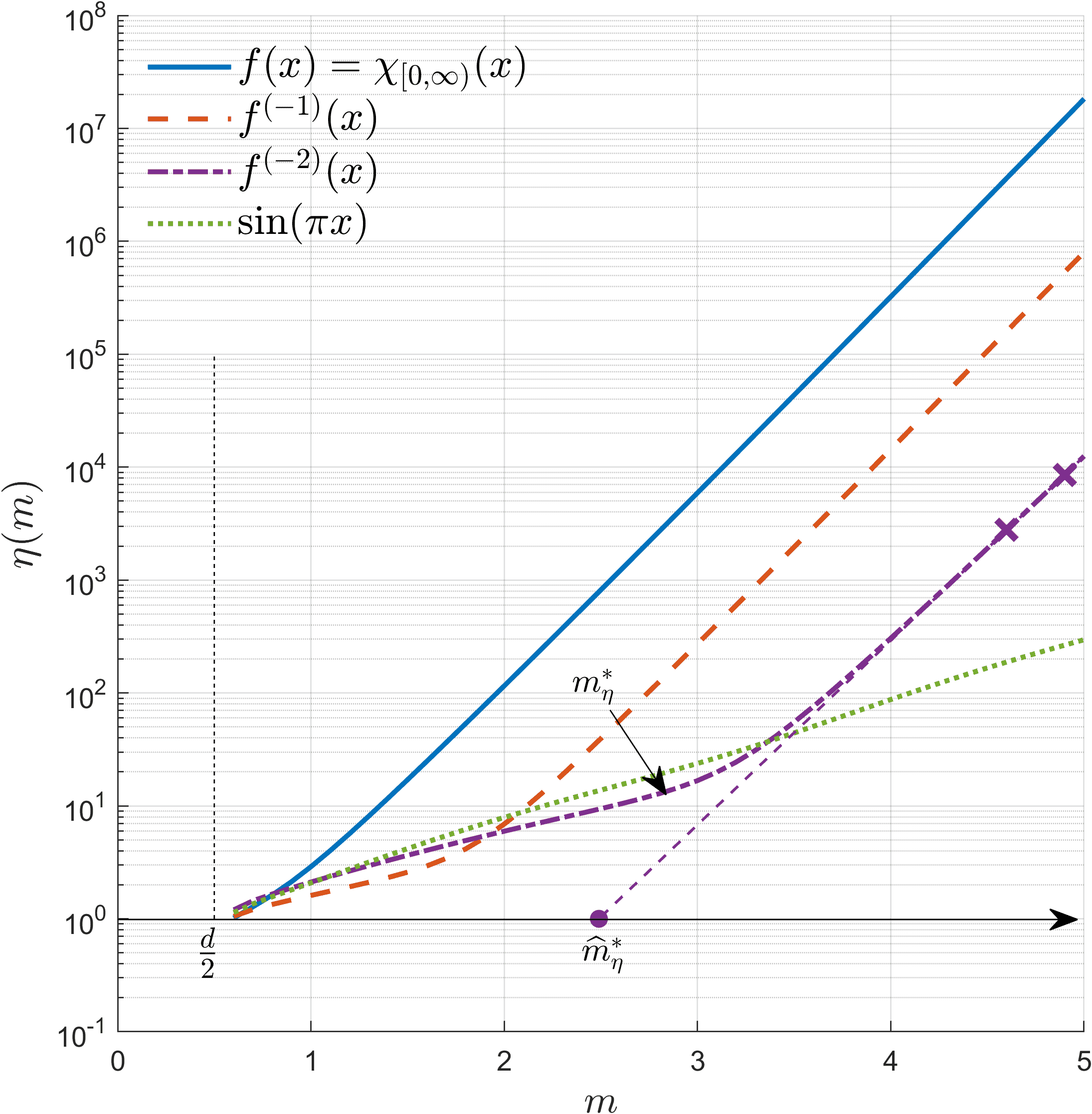}\end{overpic}
\caption{
Secant-based tail screening for candidate low-regularity neighborhoods.
Each log-profile $\log\eta(m)$ is sampled only at
$m_1=m_{\max}-\delta$ and $m_2=m_{\max}$.
If the slope at $m_{\max}$ grows significantly slower than the worst-case trend,
the neighborhood is considered locally smooth.
Otherwise, the secant between these two points approximates the high-$m$ tail,
and its intersection with $\log\eta=0$ gives an elbow-free transition estimate $\hat{m}$.
An optional one-point global screening step can further restrict such evaluations to outlier neighborhoods.
}
\label{critical_m}
\end{figure}

For each $z$, we evaluate $\eta(m)$ only at the two highest available orders,
\[
m_1 = m_{\max}-\delta, \qquad
m_2 = m_{\max},
\]
where $m_{\max}$ is the upper boundary of the tested kernel orders and $\delta>0$ is a small fixed step within
the stable range.
The pair $\{\eta(m_1),\eta(m_2)\}$ provides a short secant segment of the log-profile $\log\eta(m)$
near its high-$m$ end.
If the local slope at $m_{\max}$ is substantially smaller than that of the steep, worst-case growth in
\eqref{eq:upper-bound}, the neighborhood is screened as locally smooth and no further local sweep is triggered. 
Otherwise, the two samples define an ascending secant line that approximates the local tail of $\log\eta(m)$.
The intersection of this line with $\log\eta=0$ yields a simple tail-based transition indicator,
\begin{equation}
\label{eq:m_hat}
\hat{m}(z) = m_2 - \frac{\log\eta(m_2)}{\text{slope}(m_1,m_2)},
\end{equation}
which is used only to rank or screen candidate low-regularity
neighborhoods, rather than as a replacement for the full elbow estimator, see Figure~\ref{critical_m}.
The resulting scheme remains entirely local and requires at most two norm
evaluations per neighborhood.  It reduces the number of dense $m$-sweeps
needed in large-scale computations, while points identified as low-regularity
candidates may still be refined geometrically through the stencil-shift process
of Section~\ref{sec:stencil_shift_refinement}. 

For even larger datasets, we can include an optional global screening stage to reduce the number of these
two-point evaluations.
In this variant, we first compute a single high-order value $\eta(m_{\max})$ at all $z\in Z$.
Points whose $\log\eta(m_{\max})$ appear as outliers in the global distribution are collected into a reduced set
$Z_o$. This detection can, for example, use a robust interquartile-range (IQR) rule. 
The local two-point secant refinement is then applied only to $Z_o$, while other points 
are treated as smooth candidates unless selected by another
local criterion.

This optional one-$m$ screening further reduces the computational cost, but its
accuracy depends on the shape of the local norm profile and should be assessed
against the full sweep on representative cases. The same screening strategy can also be combined with the stencil‑shift refinement in subsection~\ref{sec:stencil_shift_refinement}, as illustrated in Figure~\ref{fig:singularity_aware_stencil_selection_1d}(a).

This screening rule is intentionally conservative. It may mark some smooth
neighborhoods as low-regularity candidates when the high-order tail has
residual curvature or when the local stencil is poorly balanced. Such false
positives mainly increase the number of neighborhoods passed to the refinement
step or slightly widen the detected transition region. This is acceptable for
the intended use, since missing genuinely low-regularity neighborhoods would be
more damaging for the subsequent stencil-shift and PDE discretization steps.

Unless otherwise stated, the large-scale numerical examples in Section~5 use
this accelerated screening rule to identify candidate low-regularity
neighborhoods. For validation, we compare the {secant} indicator with the
maximum-curvature estimate obtained from the full norm profile on representative
one-dimensional tests. The full sweep gives a sharper estimate of the elbow,
whereas the {secant} tail indicator is cheaper and conservative.

\section{Numerical Examples}
\label{sec:Numerical experiments}

This section demonstrates the proposed Sobolev-scale norm profiling procedure on one- and two-dimensional
scattered datasets, including synthetic test functions and a turbulent-flow case.
The first two examples examine parameter effects in one dimension, followed by two-dimensional and
turbulent-flow applications in later subsections.

The local analysis estimates the data-driven regularity order $\tilde{s}(\Omega_z)$ from neighborhood samples
and investigates how $\tilde{s}(\Omega_z)$ approaches the intrinsic scale $s(\Omega_z)$ as the sampling is
refined.
In all examples, the Whittle-Mat\'ern kernel family defined in~\eqref{eq:wmat-def-ex} is used, with the
shape parameter~$\varepsilon$ specified individually and smoothness index
$m \in (d/2,5]$ swept to form the discrete norm profile~$\eta(m)$. 

Unless otherwise stated, the large-scale numerical examples use the
secant-based screening indicator of Section~\ref{sec:norm_sweep_acceleration}.
We denote this quantity by $\tilde s_{\rm tail}(\Omega_z)$.
The full curvature-based estimator defined by
\eqref{eq:m-eta-star-def} is denoted by
$\tilde s_{\rm full}(\Omega_z)$ and is used only where stated explicitly.

\subsection{Example 1: Fill-Distance Refinement in One Dimension}
\label{subsec:1d-stencil}

We first study how the data resolution affects the L-curve elbow of $\eta(m)$ and the resulting
data-driven regularity estimate $\tilde{s}(\Omega_z)$.
Three reference functions on $\Omega=[-0.5,0.5]$ are considered: a step function ($H^{\frac{1}{2}-}$), a first-order kink ($H^{\frac{3}{2}-}$), and a second-order kink ($H^{\frac{5}{2}-}$), similar to the definition of Figure \ref{fig:norm-growth}.
For each target point \(z\) near the low-regularity location, a symmetric local stencil \(X_z\) is formed from \(n\) uniformly spaced points in the fixed domain \(\Omega = [-0.5, 0.5]\), with \(n \in \{5, \dots, 50\}\) and fixed shape parameter \(\varepsilon = 1\).

\begin{figure}[t]
\centering
\subfloat[Step function]{
\begin{overpic}[width=40mm,height=43mm]{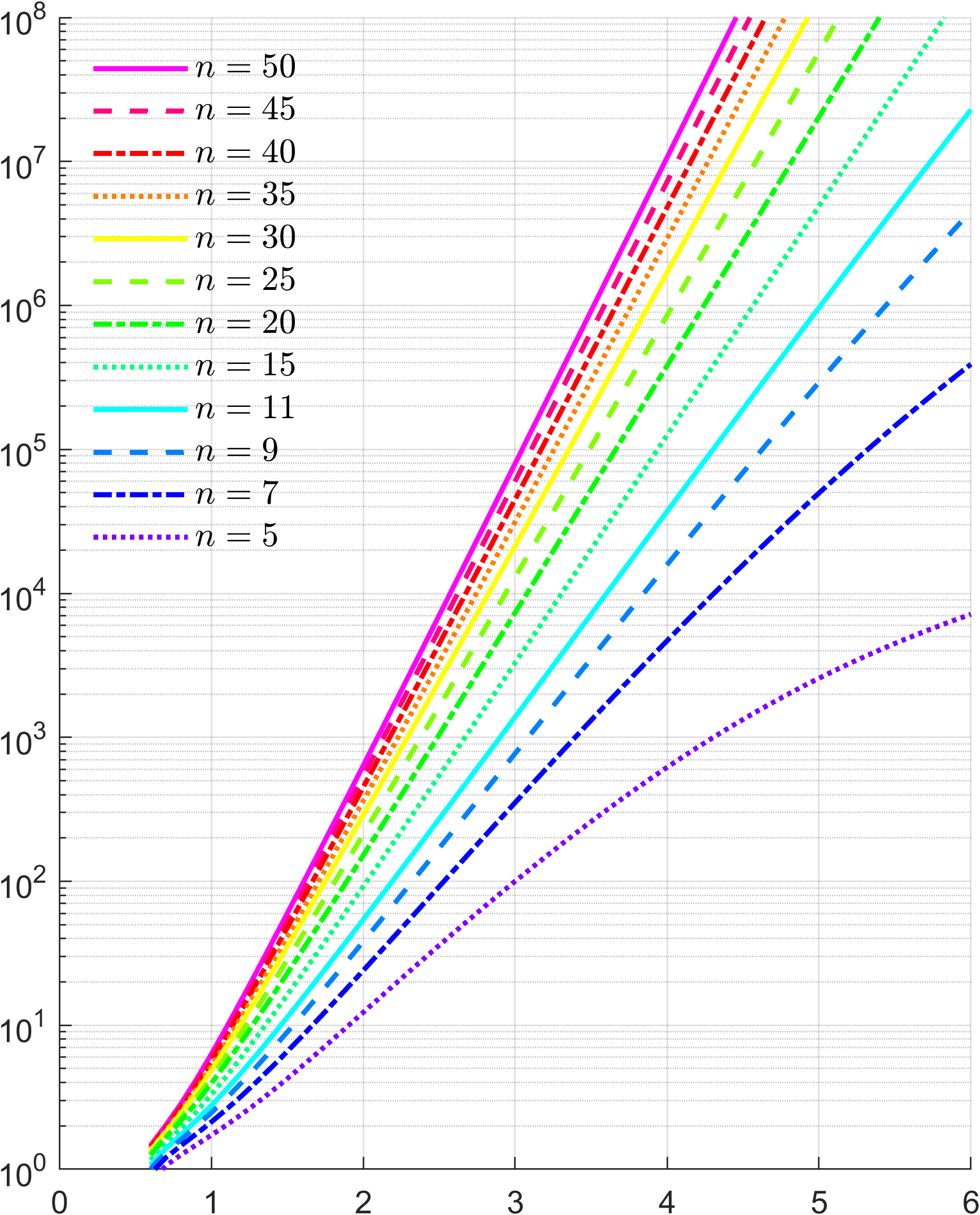}
\put(-5,45){\scriptsize\rotatebox{90}{$\eta(m)$}}
\put(45,-2){\scriptsize{$m$}}
\end{overpic}}
\hspace{2mm}
\subfloat[First-order kink]{
\begin{overpic}[width=40mm,height=43mm]{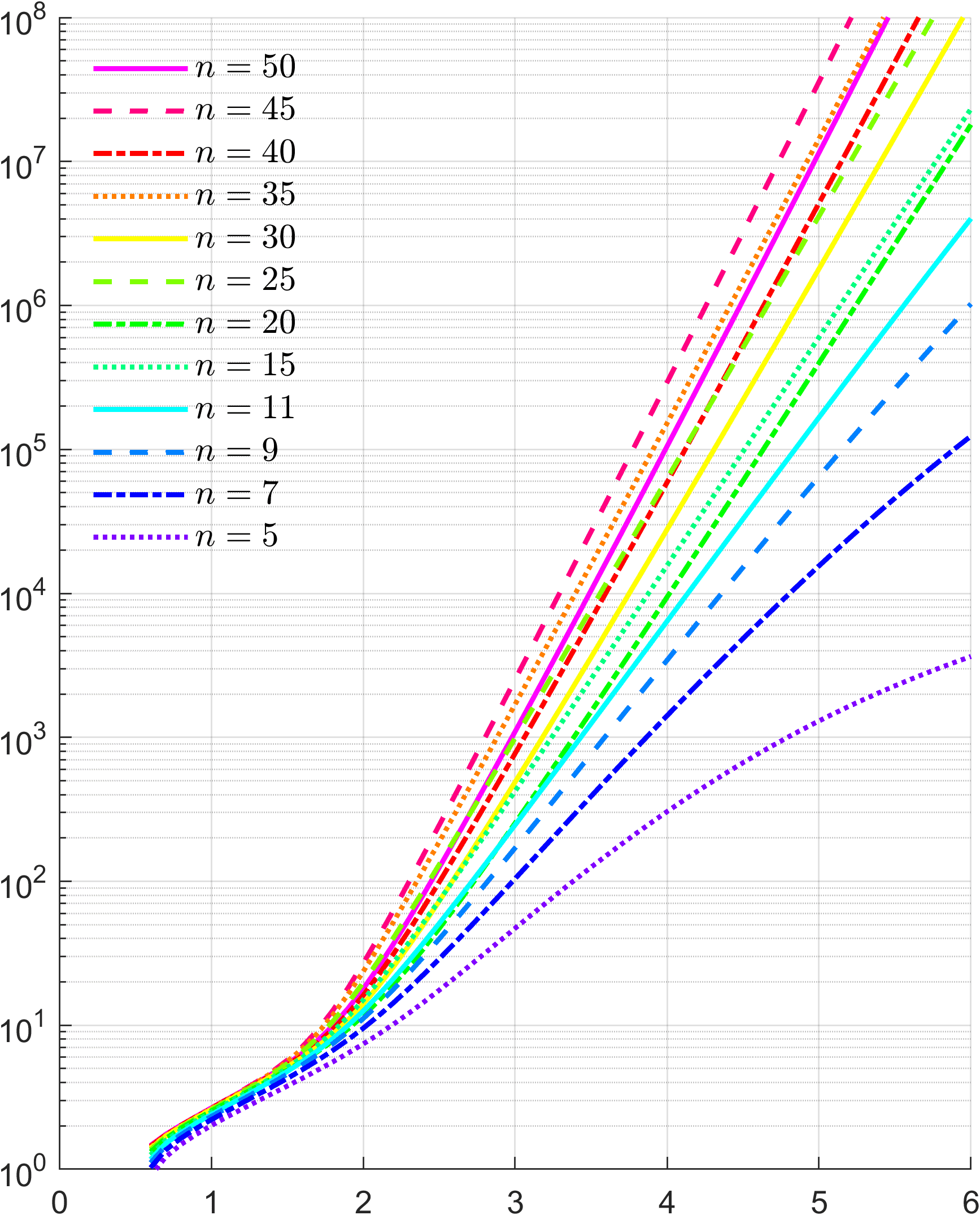}
\put(-5,45){\scriptsize\rotatebox{90}{$\eta(m)$}}
\put(45,-2){\scriptsize{$m$}}
\end{overpic}}
\hspace{2mm}
\subfloat[Second-order kink]{
\begin{overpic}[width=40mm,height=43mm]{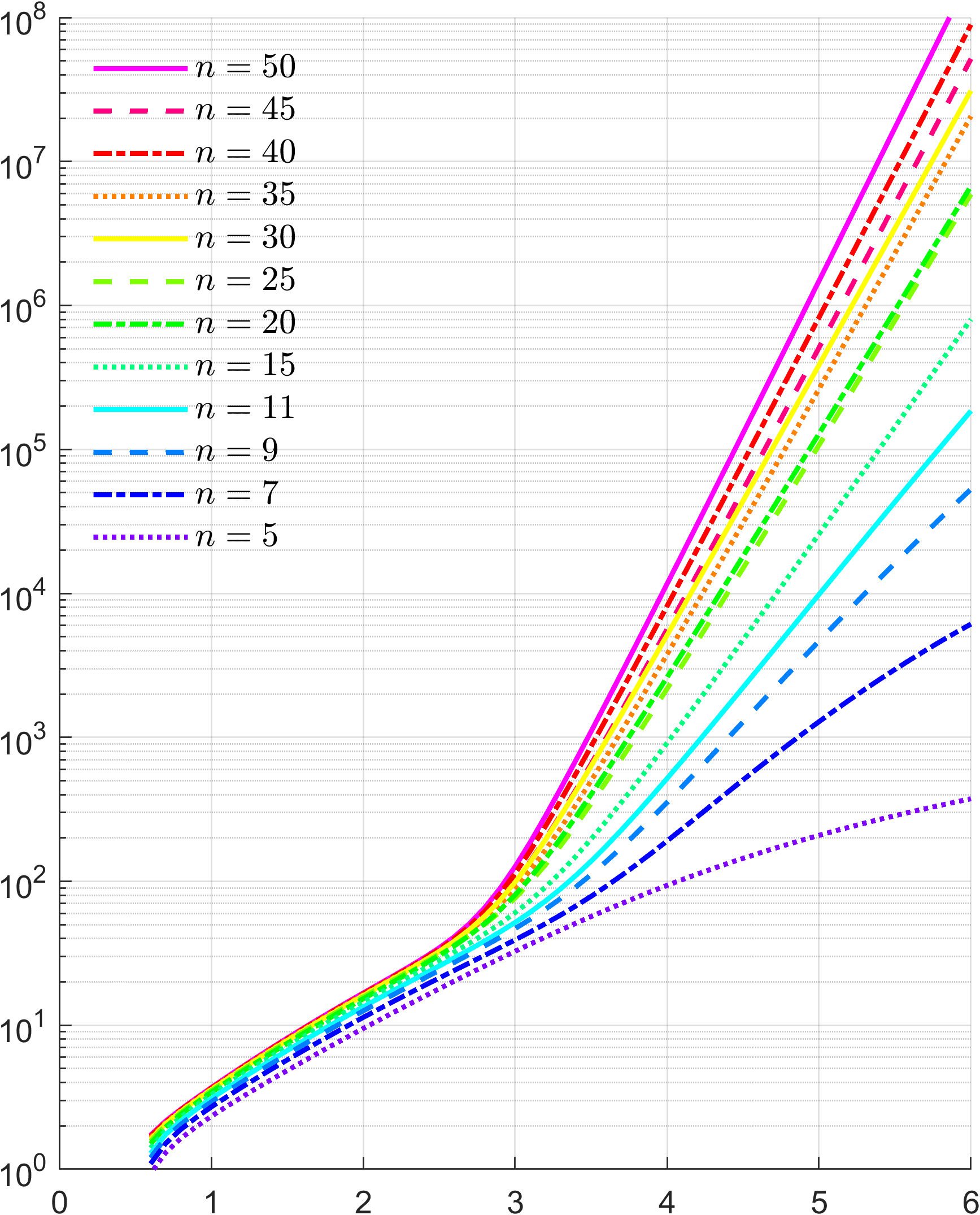}
\put(-5,45){\scriptsize\rotatebox{90}{$\eta(m)$}}
\put(45,-2){\scriptsize{$m$}}
\end{overpic}}
\caption{
\textbf{(Example 1)} Discrete native-norm profiles $\eta(m)$ for three one-dimensional test functions (a step function, a first-order kink, and a second-order kink), corresponding to Sobolev regularities $H^{\frac{1}{2}-}$, $H^{\frac{3}{2}-}$, and $H^{\frac{5}{2}-}$, respectively. A Whittle--Mat\'ern kernel with $\varepsilon=1$ is used on $\Omega=[-0.5, 0.5]$. Each curve corresponds to a different number $n$ of uniformly spaced samples in the fixed local interval. 
}
\label{fig:eta-1d-n}
\end{figure}

Figure~\ref{fig:eta-1d-n} shows the discrete profiles~$\eta(m)$ as $n$ increases.
For the step function, no clear elbow is observed; the profile enters a rapid-growth regime almost immediately across all stencil sizes, reflecting its very low regularity.
In contrast, for first- and second-order kink functions, the profiles appear smooth and lack a distinct elbow when $n$ is very small.
A distinct elbow emerges for $n\gtrsim10$ in the first-order case and $n\gtrsim15$ in the second-order case.
 This behavior illustrates the data-driven character of the norm-profile diagnostic. For small $n$, the local samples do not sufficiently resolve the singular
feature and the transition in $\log\eta(m)$ is weak. As more samples are
included, the transition becomes more pronounced near the expected Sobolev
threshold.

\begin{figure}[t]
\centering
\subfloat[Curvature]{
\begin{overpic}[width=0.3\textwidth]{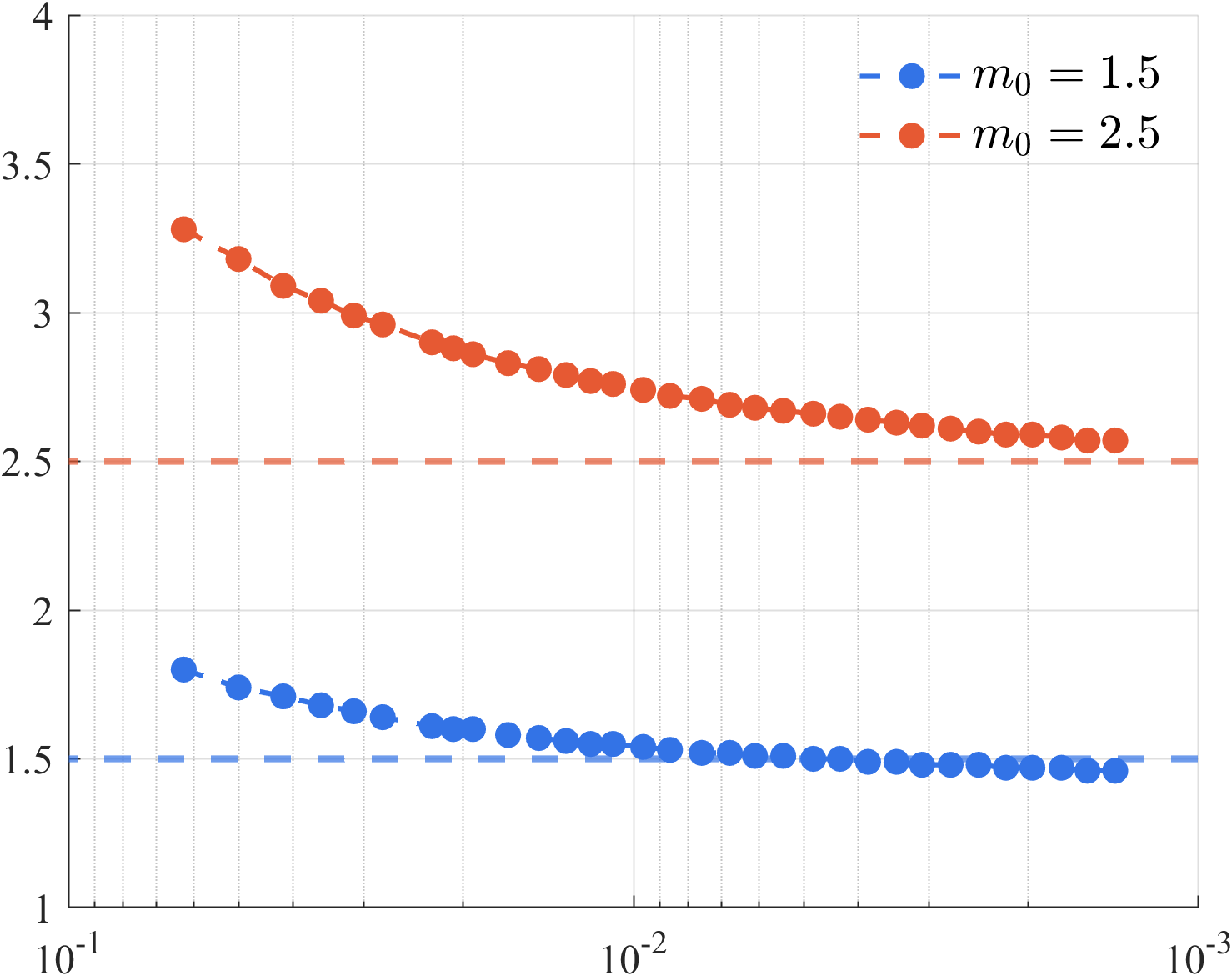}
\put(-5,30){\scriptsize\rotatebox{90}{$\widehat m_0$}}
\put(45,-2){\scriptsize{$h$}}
\end{overpic}}
\subfloat[Secant method]{
\begin{overpic}[width=0.3\textwidth]{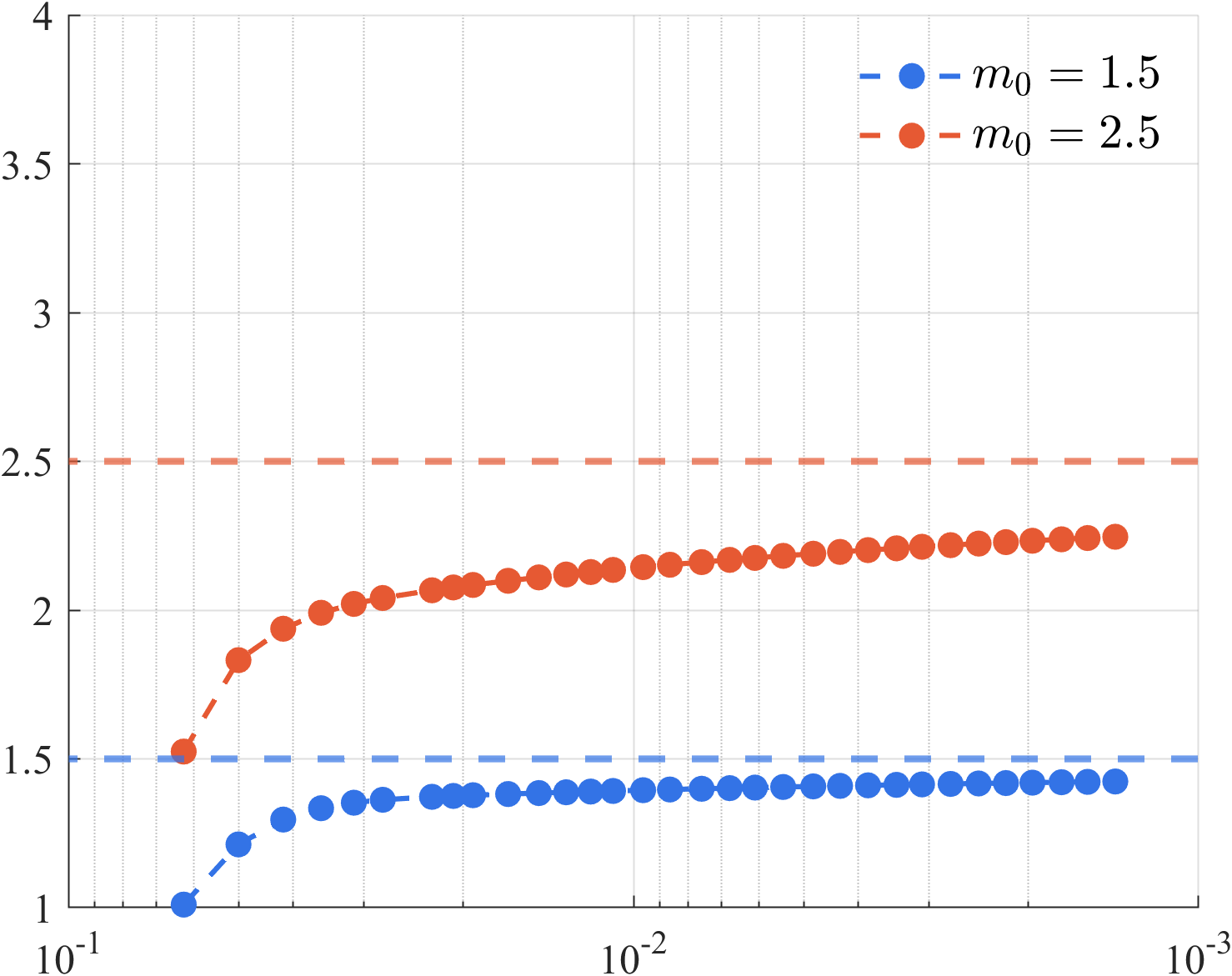}
\put(-5,30){\scriptsize\rotatebox{90}{$\widehat m_0$}}
\put(45,-2){\scriptsize{$h$}}
\end{overpic}}
\subfloat[CPU time]{
\begin{overpic}[width=0.3\textwidth]{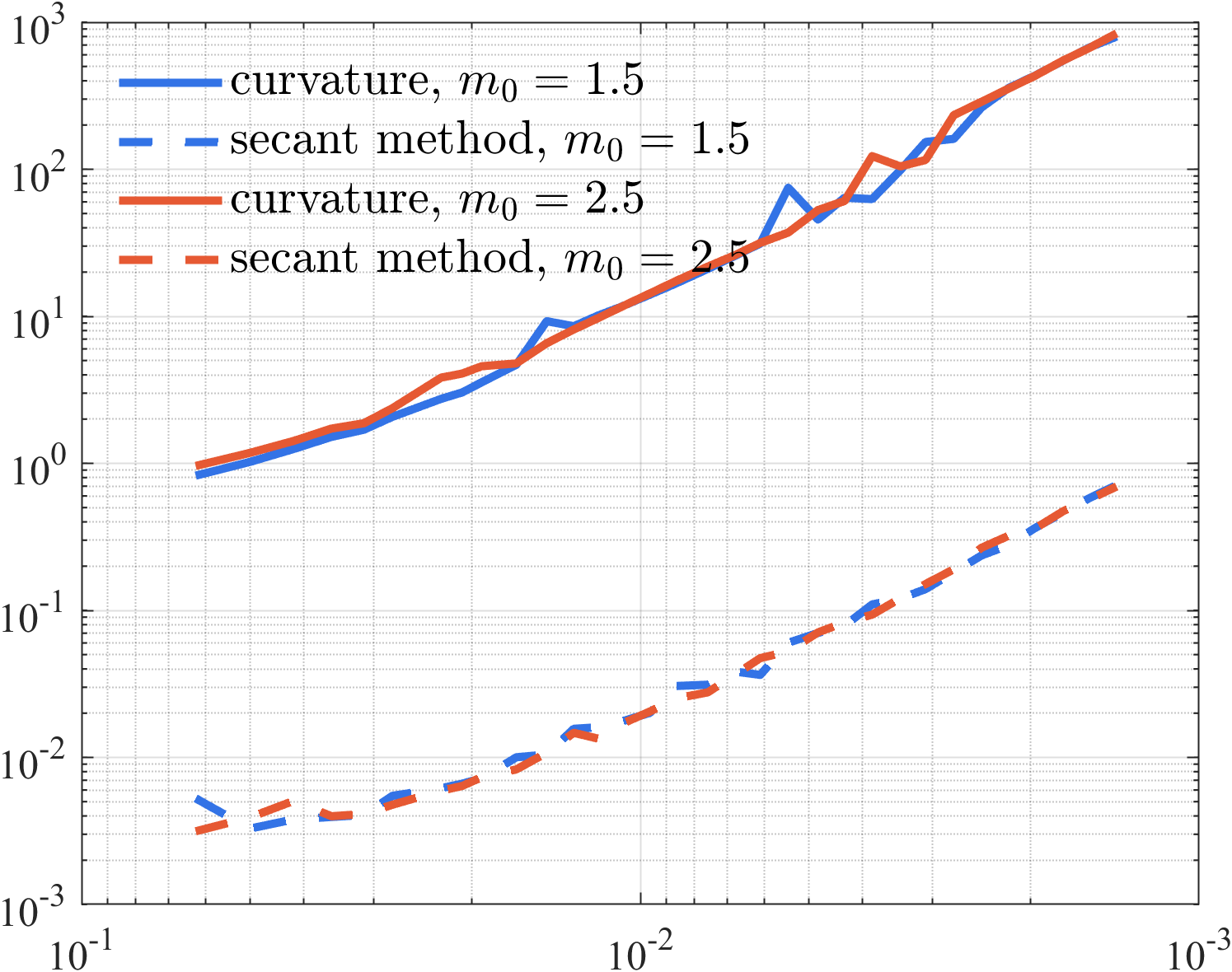}
\put(-5,30){\scriptsize\rotatebox{90}{CPU time(s)}}
\put(45,-2){\scriptsize{$h$}}
\end{overpic}}
\caption{
 \textbf{(Example 1)} Refinement behavior of regularity indicators under decreasing fill
distance $h$. The estimated threshold $\widehat m_0$ is plotted for
$f(x)=x_+$ with $m_0=1.5$ and $f(x)=x_+^2/2$ with $m_0=2.5$.
{Panels (a)} show the maximum-curvature estimate obtained from the full
Sobolev-scale norm profile by \eqref{eq:m-eta-star-def}, while panels {(b)} show the  secant
estimate used in the accelerated screening strategy of
Section~\ref{sec:norm_sweep_acceleration} with $m_{\max}=5.5$ and $\delta=0.5$. Panels {(c)} show the corresponding CPU times. The dashed horizontal
lines indicate the exact Sobolev thresholds. 
}
\label{fig:compare_estimators_1d}
\end{figure}

We next examine the estimated Sobolev threshold under fill-distance refinement. The interval $[-0.5,0.5]$ is sampled uniformly with increasing point density. We consider the first-order kink $f(x)=x_+$ with threshold $m_0=3/2$ and the second-order kink $f(x)=x_+^2/2$ with threshold $m_0=5/2$. The jump case is omitted because its threshold $1/2$ lies below the tested Sobolev range. Figure~\ref{fig:compare_estimators_1d}(a) shows the maximum-curvature estimate. The profile $\eta(m)$ is evaluated on a uniform grid $m\in[0.6,5.5]$ with $\Delta m=0.01$, and the derivatives of $\log\eta(m)$ are approximated by centered finite differences. The estimates approach the exact thresholds as $h$ decreases. The small oscillations are caused by whether the singular point $x=0$ is sampled directly or lies between two neighboring grid points. Figure~\ref{fig:compare_estimators_1d}(b) shows the secant-based indicator with $m_1=5$ and $m_2=5.5$, and Figure~\ref{fig:compare_estimators_1d}(c) reports the corresponding CPU times. For the $H^{3/2-}$ kink, the secant indicator approaches the threshold from below. For the $H^{5/2-}$ kink, it increases under refinement but remains below $m_0=2.5$ over the tested range. The secant-based criterion is therefore used as a low-cost screening indicator, while the full curvature criterion provides a sharper estimate of the elbow.

\subsection{Example 2: Effect of Kernel Shape Parameter in 1D}
\label{subsec:1d-shape}

We next investigate how the kernel shape parameter $\varepsilon$ influences the regularity estimate for the
piecewise-smooth function $f$ defined in Appendix~\ref{appendix:testfun1d}.
The function \eqref{eq:testfun1d} contains multiple singularities on $\Omega=[-1,1]$, including jump discontinuities, a cusp, and a
$C^1$ but not $C^2$ point, as already illustrated in Figure \ref{fig:fx_order}.
Experiments use $2000$ uniformly spaced samples, $n=20$ nearest neighbors, and
Whittle-Mat\'ern kernels with $\varepsilon\in\{0.6,0.8,1.0\}$.
The smoothness index $m$ is swept over $(d/2,5]$.
We also test two variable-shape cases $\varepsilon(m)=m$ and $\varepsilon(m)=1/m$ to illustrate
non-constant parameterization along the norm-sweep.

\begin{figure}[t]
\centering
\subfloat[ $z=-0.2$]{
\begin{overpic}[width=40mm,height=43mm]{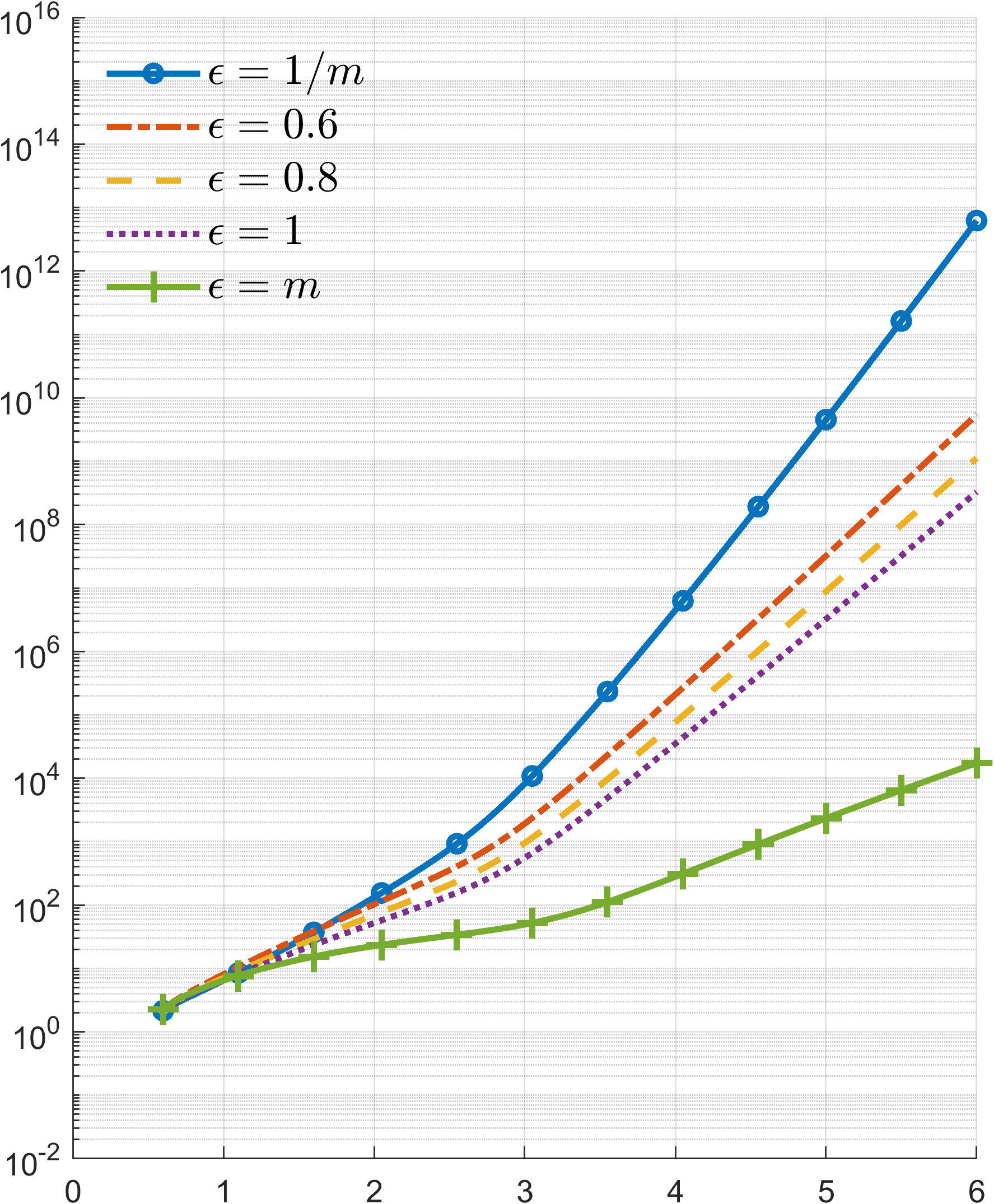}
\put(-5,45){\scriptsize\rotatebox{90}{$\eta(m)$}}
\put(52,-2){\scriptsize{$m$}}
\end{overpic}}
\hspace{2mm}
\subfloat[ $z=-0.6$]{
\begin{overpic}[width=40mm,height=43mm]{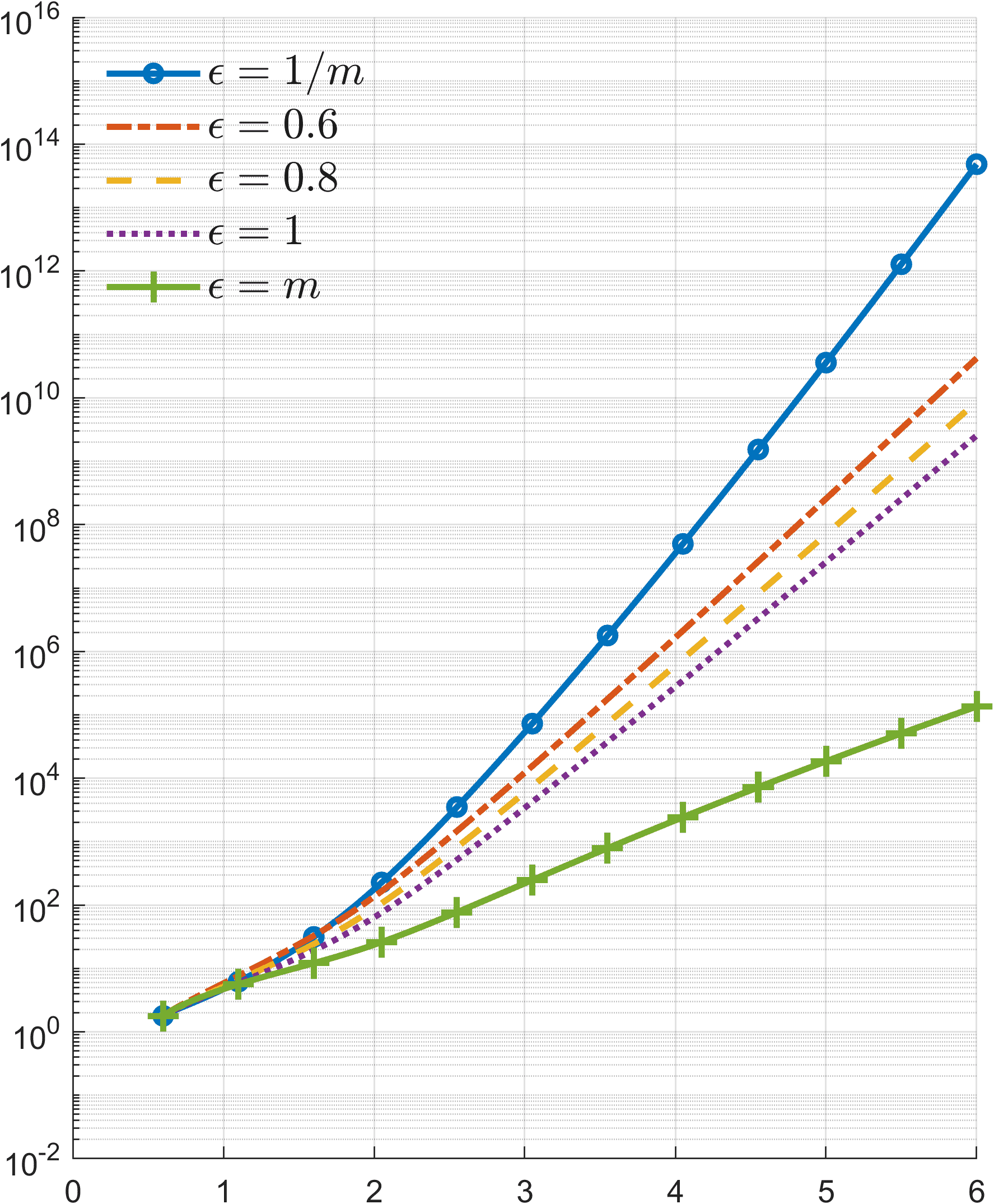}
\put(-5,45){\scriptsize\rotatebox{90}{$\eta(m)$}}
\put(52,-2){\scriptsize{$m$}}
\end{overpic}}
\hspace{2mm}
\subfloat[Theoretical upper bound]{
\begin{overpic}[width=40mm,height=43mm]{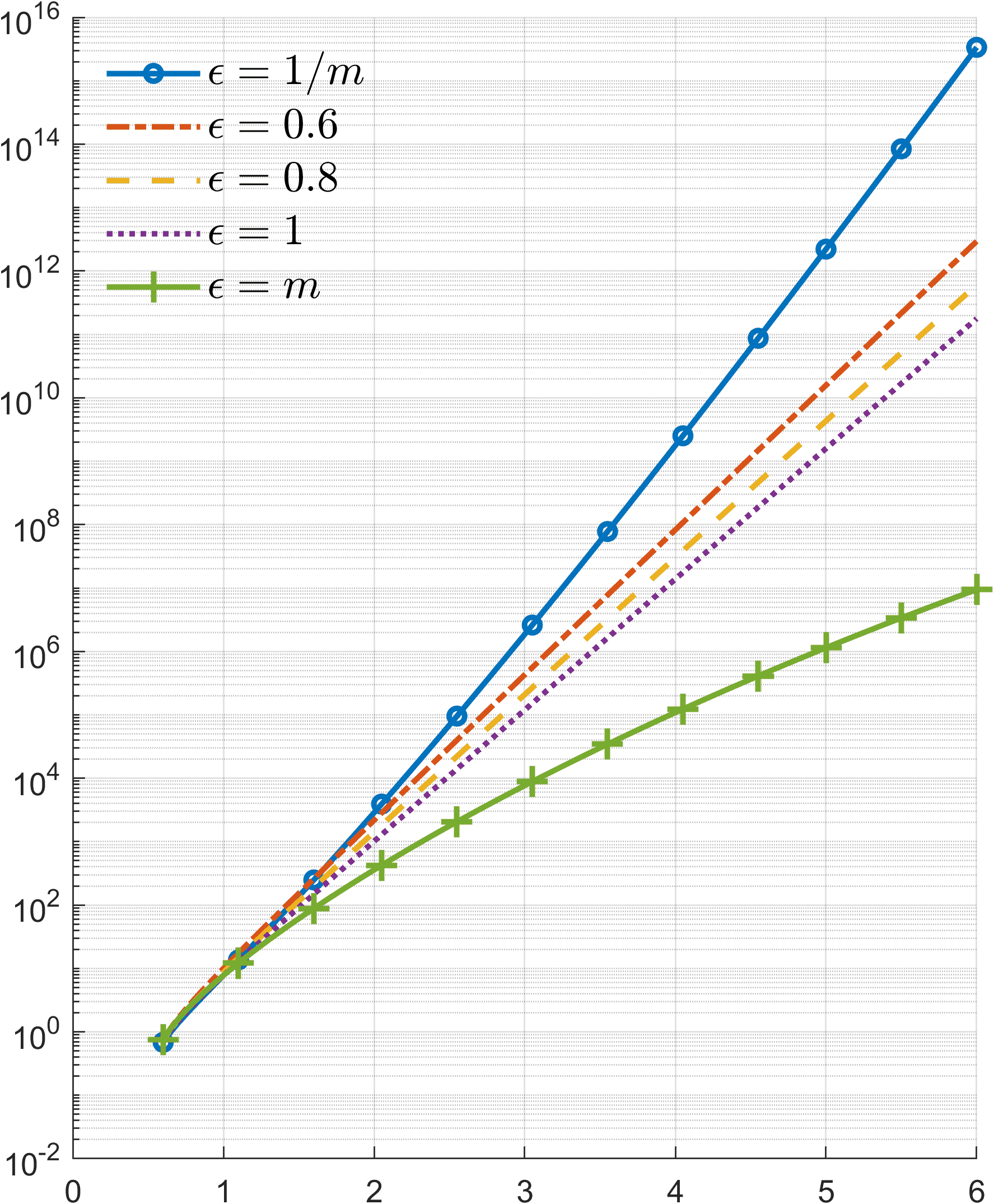}
\put(-8,30){\tiny\rotatebox{90}{$\lambda_{\min}^{-1/2}(\Phi_m(X,X))$}}
\put(52,-2){\scriptsize{$m$}}
\end{overpic}}
\caption{
\textbf{(Example 2)} Effect of the shape parameter $\varepsilon$ on the norm profiles $\eta(m)$ and the theoretical upper bound $\lambda_{\min}^{-1/2}(\Phi_m(X,X))$ in \eqref{eq:upper-bound}, evaluated for the function $f$ in \eqref{eq:testfun1d} at the local neighborhoods centered at $z=-0.6$ and $z=-0.2$, using $n = 20$ nearest neighbors.
}
\label{fig:eta-1d-eps}
\end{figure}

Figure~\ref{fig:eta-1d-eps} illustrates the influence of the kernel shape parameter~$\varepsilon$ on the
discrete norm profiles~$\eta(m)$. We focus on two representative locations, $z=-0.6$ and $z=-0.2$, corresponding to local regularities in the sense of $H^{2.5-}$ and $H^{1.5-}$ in $f$.
Smaller~$\varepsilon$ values produce flatter kernels with broader support, which sharpen the elbow in
$\eta(m)$ and improve the detection of local non-smoothness.
Larger~$\varepsilon$ values yield more peaked kernels; their weaker overlap across nodes reduces coupling and
makes $\tilde{s}(\Omega_z)$ less responsive to discontinuities.
However, excessively small~$\varepsilon$ also increases the condition number of the kernel matrix
$\Phi_m(X,X)$, as indicated by the growth of
$\lambda_{\min}^{-1/2}(\Phi_m(X,X))$ in the figure.
Among the variable-shape cases, $\varepsilon(m)=m$ preserves good conditioning as~$m$ increases,
whereas $\varepsilon(m)=1/m$ emphasizes smoothing at large~$m$.
These results show that~$\varepsilon$ need not remain fixed during the~$m$-sweep:
moderate adaptation of~$\varepsilon$ with~$m$ can help balance sensitivity to small-scale irregularities
against numerical stability.

\begin{table}[t]
\centering
\caption{
 \textbf{(Example 2)} Data-driven regularity $\widetilde s(\Omega_z)$
computed from the curvature criterion~\eqref{eq:m-eta-star-def} at selected
locations of the one-dimensional piecewise-smooth test function, for
$\varepsilon\in\{0.6,0.8,1.0\}$ and stencil sizes $n=20$ and $n=30$.
}
\label{tab:1d-epsilon-by-location}
\setlength{\tabcolsep}{4pt}
\renewcommand{\arraystretch}{1.15}
\begin{tabular}{c *{6}{cc}}
\toprule
location & \multicolumn{2}{c}{$z=-0.8$}
& \multicolumn{2}{c}{$z=-0.6$}
& \multicolumn{2}{c}{$z=-0.4$}
& \multicolumn{2}{c}{$z=-0.2$}
& \multicolumn{2}{c}{$z=0$}
& \multicolumn{2}{c}{$z=0.5$} \\
\multicolumn{1}{c}{target $s$} &
\multicolumn{2}{c}{0.5} & \multicolumn{2}{c}{2.5} & \multicolumn{2}{c}{0.5} &
\multicolumn{2}{c}{1.5} & \multicolumn{2}{c}{0.5} & \multicolumn{2}{c}{0.5} \\
\cmidrule(lr){2-3}\cmidrule(lr){4-5}\cmidrule(lr){6-7}\cmidrule(lr){8-9}\cmidrule(lr){10-11}\cmidrule(lr){12-13}
$n$& \multicolumn{1}{c}{$20$} & \multicolumn{1}{c}{$30$}
& \multicolumn{1}{c}{$20$} & \multicolumn{1}{c}{$30$}
& \multicolumn{1}{c}{$20$} & \multicolumn{1}{c}{$30$}
& \multicolumn{1}{c}{$20$} & \multicolumn{1}{c}{$30$}
& \multicolumn{1}{c}{$20$} & \multicolumn{1}{c}{$30$}
& \multicolumn{1}{c}{$20$} & \multicolumn{1}{c}{$30$} \\
\midrule
$\varepsilon=1.0$ & 0.84 & 0.78 & 5.00 & 2.82 & 0.80 & 0.74 & 2.00 & 1.76 & 0.86 & 0.80 & 0.74 & 0.64 \\
$\varepsilon=0.8$ & 0.80 & 0.76 & 2.92 & 2.70 & 0.76 & 0.72 & 1.82 & 1.66 & 0.82 & 0.78 & 0.66 & 0.64 \\
$\varepsilon=0.6$ & 0.78 & 0.74 & 2.76 & 2.54 & 0.74 & 0.72 & 1.70 & 1.56 & 0.78 & 0.74 & 0.66 & 0.64 \\
\bottomrule
\end{tabular}
\end{table}

Table~\ref{tab:1d-epsilon-by-location} reports the data-driven regularity $\widetilde s(\Omega_z)$ at selected locations for different shape parameters $\varepsilon$ and stencil sizes $n$. The reported values are computed using the full curvature criterion~\eqref{eq:m-eta-star-def}. For each case, $\eta(m)$ is evaluated on a uniform grid $m\in[0.6,5]$ with $\Delta m=0.02$, and the first and second derivatives of $\log\eta(m)$ are approximated by centered finite differences. For the scaling $\Phi_{m,\varepsilon}(x,z) =\Phi_m(\varepsilon\|x-z\|_2)$ used here, smaller $\varepsilon$ gives a flatter kernel. In this test, decreasing $\varepsilon$ generally moves the detected values closer to the nominal Sobolev thresholds, particularly at the weak $C^1$-but-not-$C^2$ transition near $z=-0.6$. The comparison between $n=20$ and $n=30$ measures the effect of the local stencil size rather than data-resolution refinement, since both cases use the same uniformly spaced 2000-point dataset.


\subsection{Example 3: Regularity Mapping and Regularity-Limited Differentiation on a 2D Synthetic Function}
\label{subsec:example3}

\begin{figure}[t]
\centering
\subfloat[top view]{
\begin{overpic}[width=0.3\textwidth]{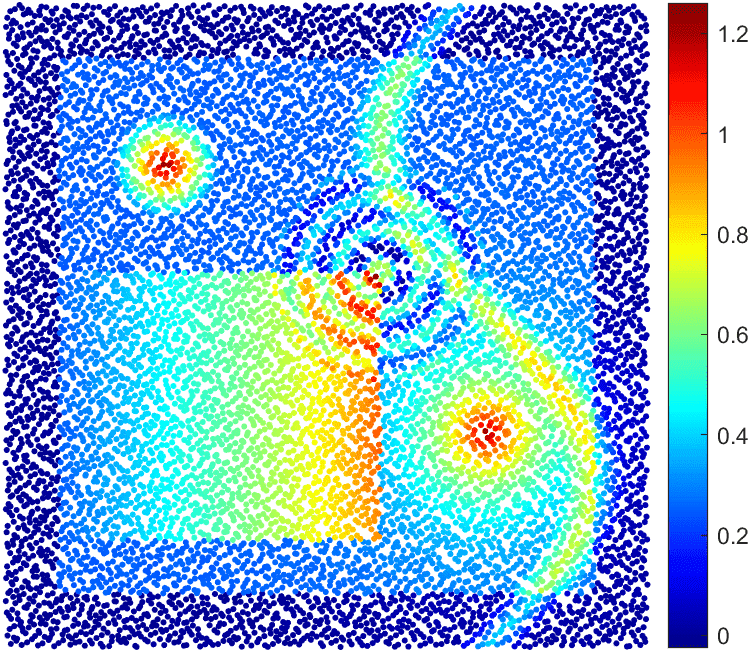}\end{overpic}}
\hspace{10mm}
\subfloat[3D view]{
\begin{overpic}[width=0.34\textwidth]{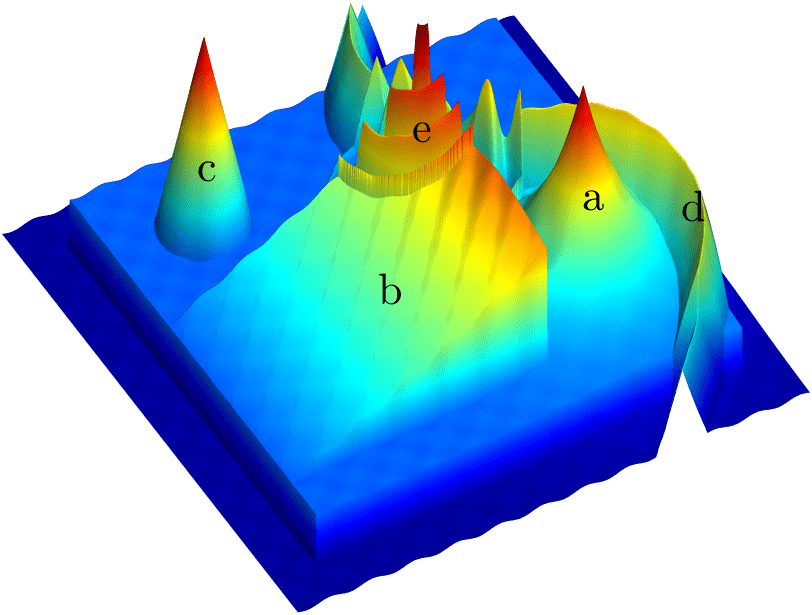}\end{overpic}}
\caption{
\textbf{(Example 3)} Visualization of the synthetic test function $f(x,y)$ on $\Omega=[0,6]^2$.
Distinct features include cusps, sharp ridges, a cone, curved ridges,
and a small oscillatory patch.
}
\label{fig:testfun_visualization}
\end{figure}

We apply the proposed Sobolev-scale norm profiling framework to the two-dimensional synthetic test function
$f(x,y)$ defined on $\Omega=[0,6]^2$, see Appendix~\ref{appendix:testfun} for the definition.
The function \eqref{eq:testfun2d} combines a smooth oscillatory background with non-smooth structures such as cusps, linear and curved ridges, a cone, and a localized high-frequency disk, see Figure~\ref{fig:testfun_visualization}.

A total of $N=160{,}000$ quasi-random Halton nodes are used, with the same number of uniform evaluation
points.  At each evaluation point $z$, a local stencil of $n=50$ nearest neighbors is selected.
For each stencil, the high-order native norms are evaluated using
Whittle--Mat\'ern kernels with $\varepsilon=1$, and the accelerated
tail-based indicator $\tilde s_{\rm tail}(\Omega_z)$ is computed as in
Section~\ref{sec:norm_sweep_acceleration}. 

\begin{figure}[t]
\centering
\subfloat[Native norm $\eta(m_{\max})$]{
\begin{overpic}[width=0.32\textwidth]{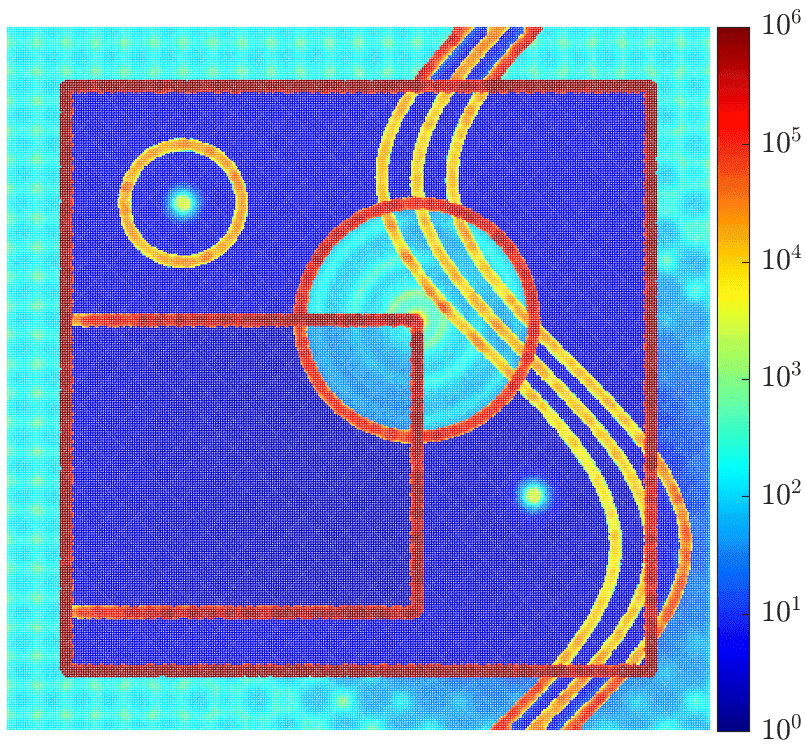}\label{fig:outlier_screening_a}\end{overpic}}
\subfloat[$\tilde{s}(\Omega_z)$ without stencil-shift]{
\begin{overpic}[width=0.32\textwidth]{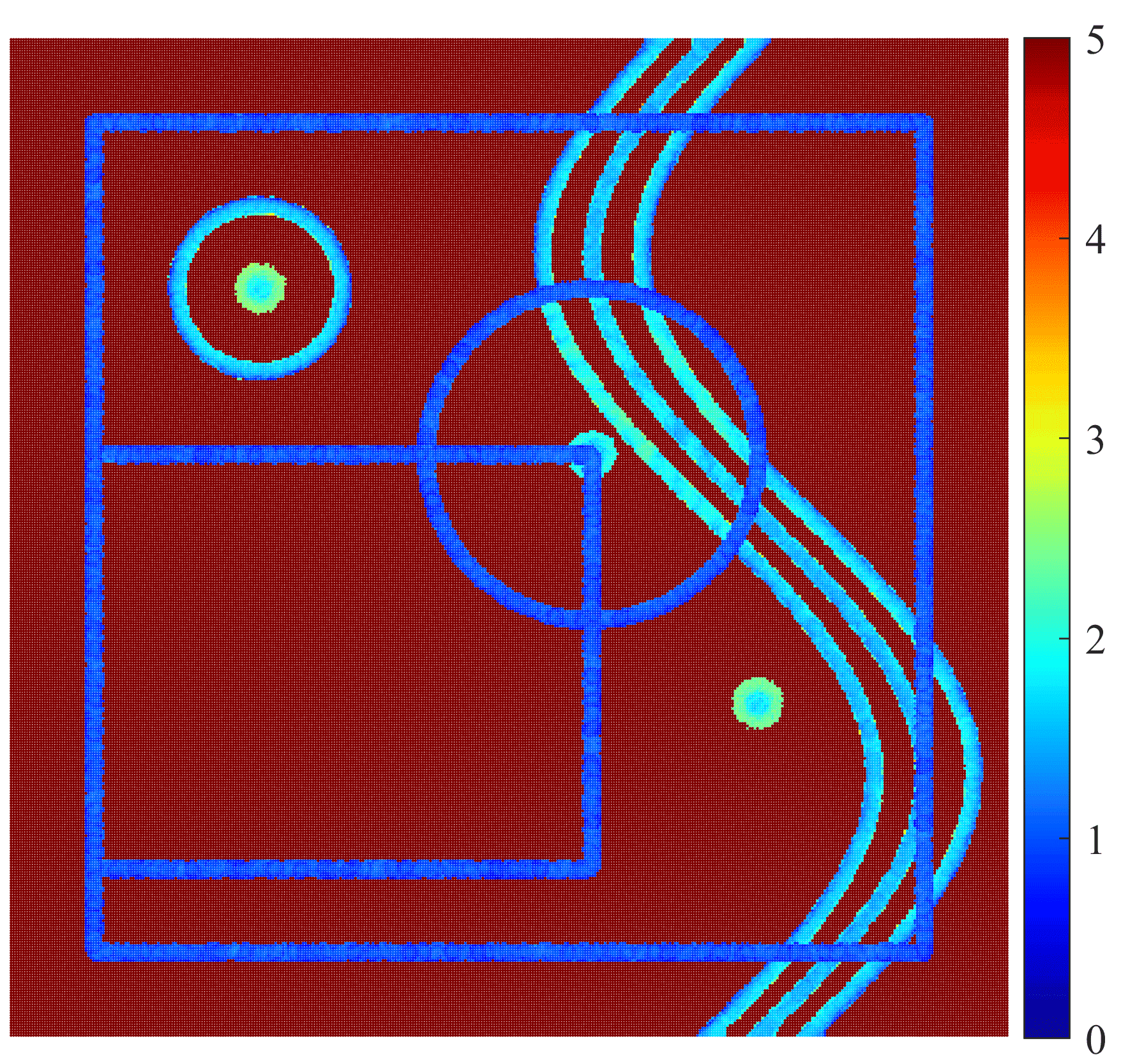}\label{fig:outlier_screening_b}\end{overpic}}
\subfloat[ $\tilde{s}(\Omega_z)$ with stencil-shift]{
\begin{overpic}[width=0.32\textwidth]{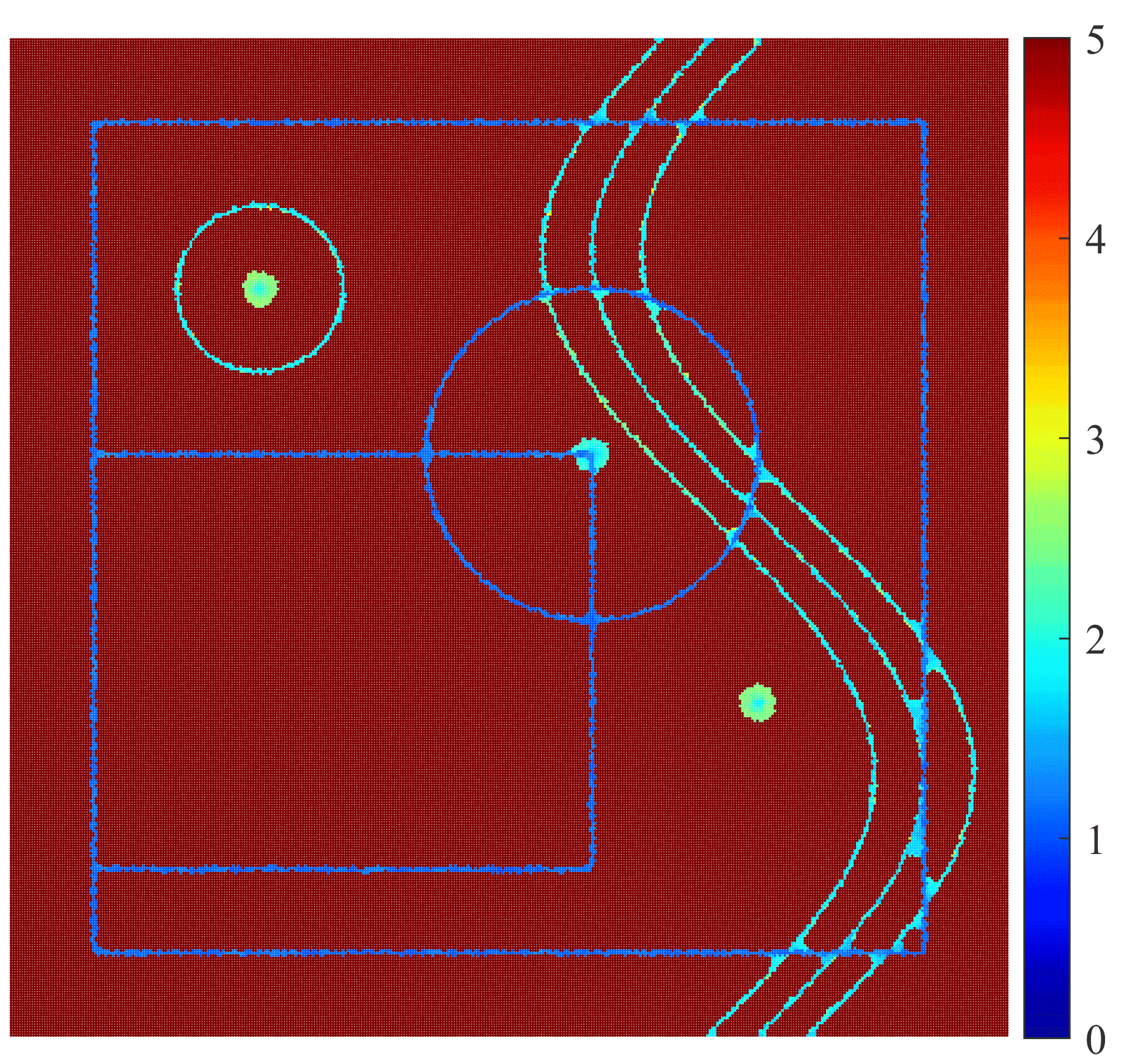}\label{fig:outlier_screening_c}\end{overpic}}
\caption{
\textbf{(Example 3)} Spatial maps of the native‑norm magnitude and data‑driven regularity $\tilde{s}(\Omega_z)$ for $f(x,y)$ on $\Omega=[0,6]^2$ with $n=50$ and $\varepsilon=1$.
}
\label{fig:outlier_screening}
\end{figure}

Figure~\ref{fig:outlier_screening} illustrates the local screening and refinement results.
Figure~\ref{fig:outlier_screening}(a) shows the native-norm values~$\eta(m_{\max})$, where large magnitudes concentrate near sharp features.
Figure~\ref{fig:outlier_screening}(b) displays the initial regularity map without stencil shifting, which already identifies the primary nonsmooth regions and their approximate types through the estimated~$\tilde{s}(\Omega_z)$.
Figure~\ref{fig:outlier_screening}(c) presents the refined map after stencil shifting, in which the low-regularity zones become more compact.

To illustrate the practical use of the data-driven regularity, we now view
$\tilde{s}(\Omega_z)$ from the perspective of stencil design.
In a PDE context, the target differential order~$k=|\alpha|$ is known in advance.
Regularity screening then identifies locations where the required Sobolev condition
$s(\Omega_z) > k + d/2$
is violated.
In practice, $\tilde{s}(\Omega_z)$ acts as a spatial reliability map:
regions satisfying the condition can support classical $k$-th order differentiation,
while regions where it fails mark zones that require a fallback or regularity-aware stencil strategy.

\begin{figure}[t]
\centering
\subfloat[$|\alpha|=2$]{
\begin{overpic}[width=0.29\textwidth]{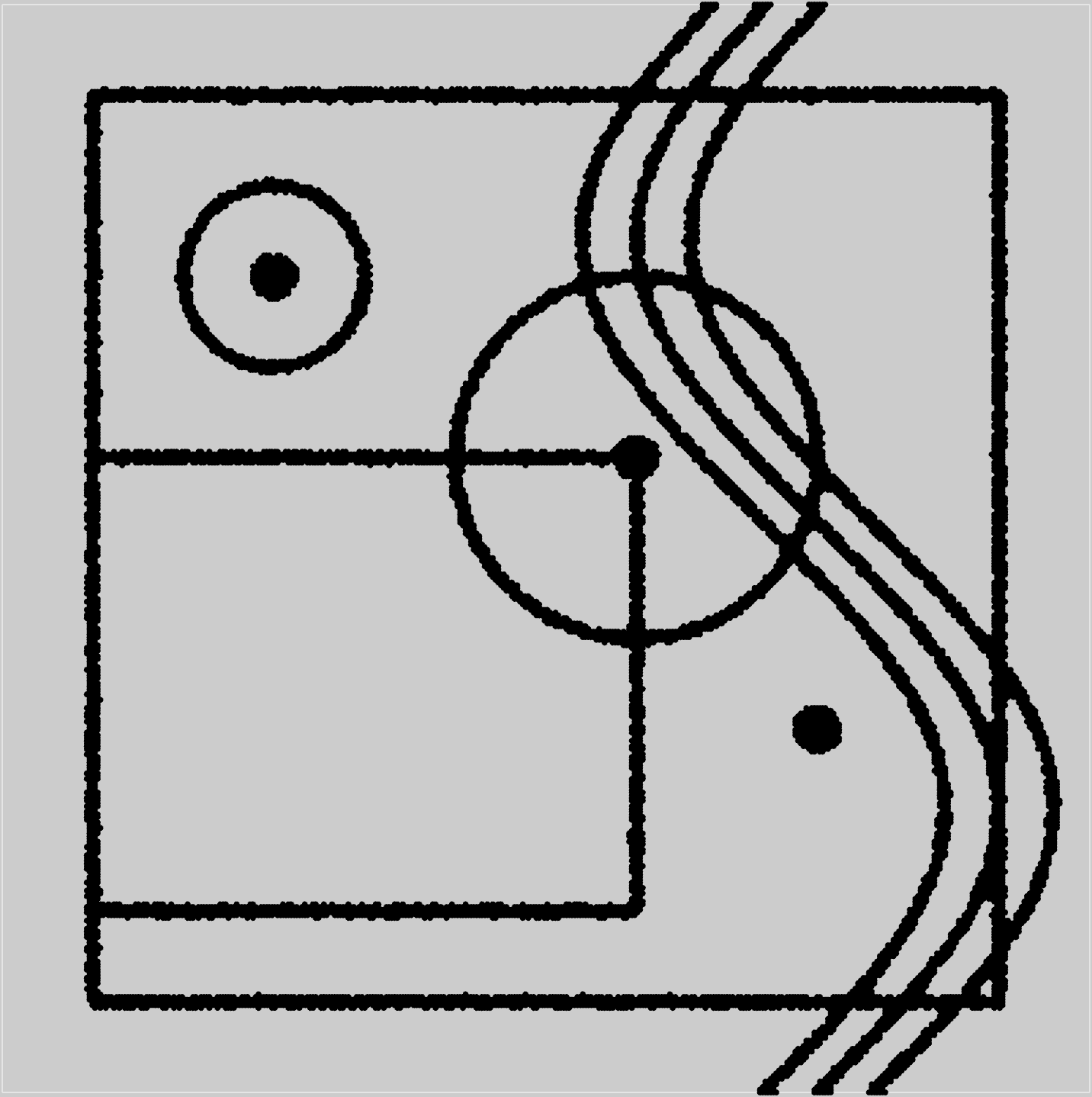}\end{overpic}}
\hspace{2mm}
\subfloat[$|\alpha|=1$]{
\begin{overpic}[width=0.29\textwidth]{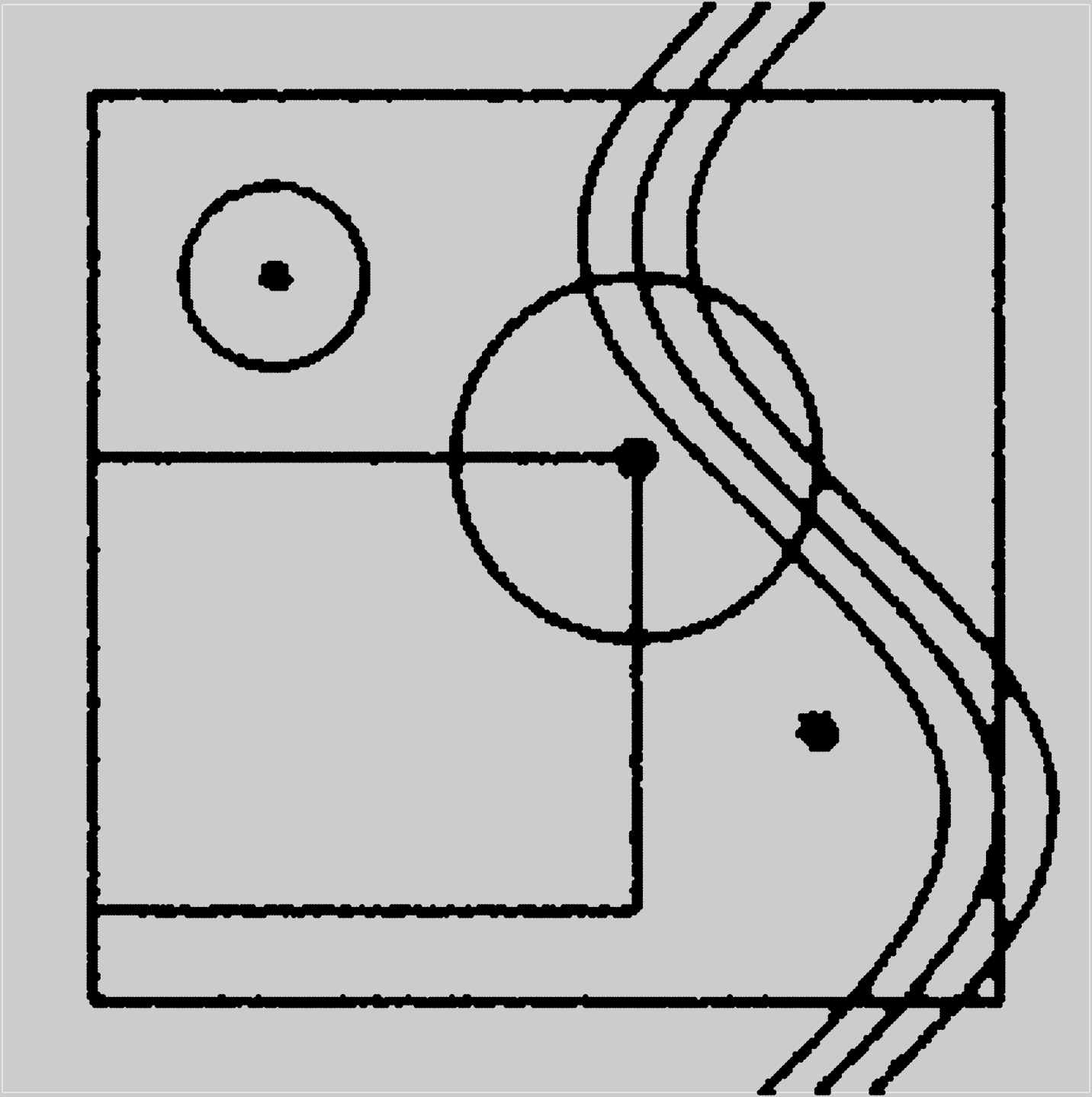}\end{overpic}}
\hspace{2mm}
\subfloat[$|\alpha|=0$]{
\begin{overpic}[width=0.29\textwidth]{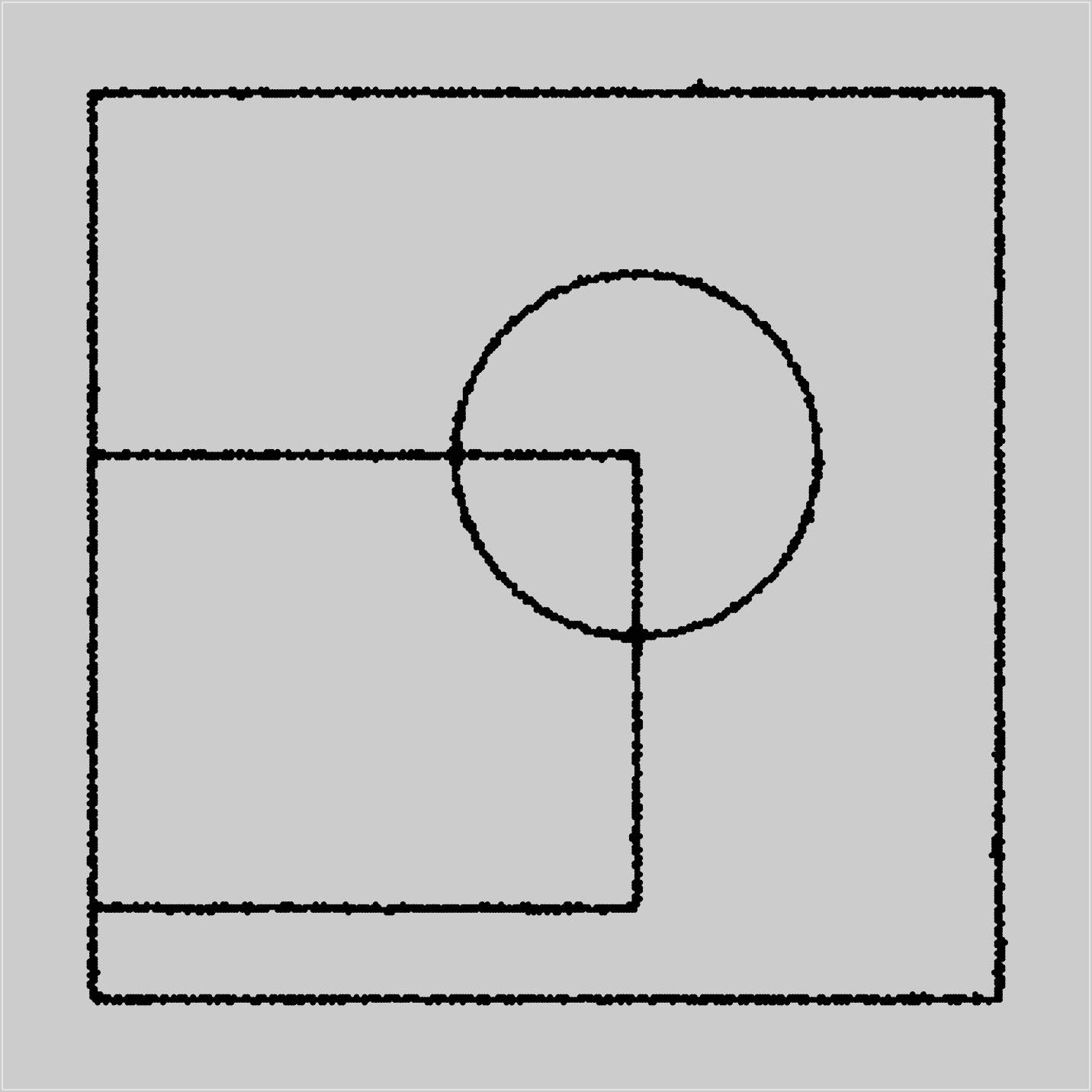}\end{overpic}}
\caption{ \textbf{(Example 3)}
Regions identified by the regularity‑limited criterion $\tilde{s}(\Omega_z)\le|\alpha|+d/2$ for nominal derivative orders
(a) $|\alpha|=2$, (b) $|\alpha|=1$, and (c) $|\alpha|=0$, based on adaptive stencil shifting.
}
\label{fig:singularity_detection_asymmetric}
\end{figure}

In the stencil-shift step, we employ an order-dependent interior admissibility margin
$\operatorname{dist}\!\bigl(z,\partial\operatorname{conv}(\mathcal{C}_z)\bigr)\ge c_{|\alpha|}\,q_{X_z},$
taking $c_{0}=0$ (interpolation only requires $z\in\operatorname{conv}(\mathcal{C}_z)$) and
increasing $c_{|\alpha|}$ with $|\alpha|$ (e.g., one- and two-spacing buffers for first- and second-order
operators), and discard candidates that violate this condition.
Figure~\ref{fig:singularity_detection_asymmetric} illustrates the regularity-based stencil reliability for different derivative orders and shows the regions flagged for $|\alpha|=2$, $1$, and $0$, respectively.
Black areas with $\tilde{s}(\Omega_z)\le|\alpha|+d/2$ indicate locally insufficient regularity to support accurate RBF-FD differentiation.
For $|\alpha|=2$, most smooth interior points satisfy the condition, whereas for $|\alpha|\le1$, only narrow zones near singular features remain excluded.
This interpretation recasts the regularity map as a data-driven filter that identifies where finite-difference operators of a chosen order are reliably admissible.
After adaptive stencil shifting, the number of flagged points decreases, respectively,
from $28{,}556$ to $16{,}667$ for $|\alpha|=2$,
from $28{,}376$ to $11{,}551$ for $|\alpha|=1$, and
from $20{,}062$ to $4{,}508$ for $|\alpha|=0$,
corresponding to reductions of roughly $40\%$--$80\%$ across these cases and yielding a substantially more localized exclusion of low-regularity neighborhoods.

\begin{figure}[t]
\centering
\begin{overpic}[width=0.32\textwidth]{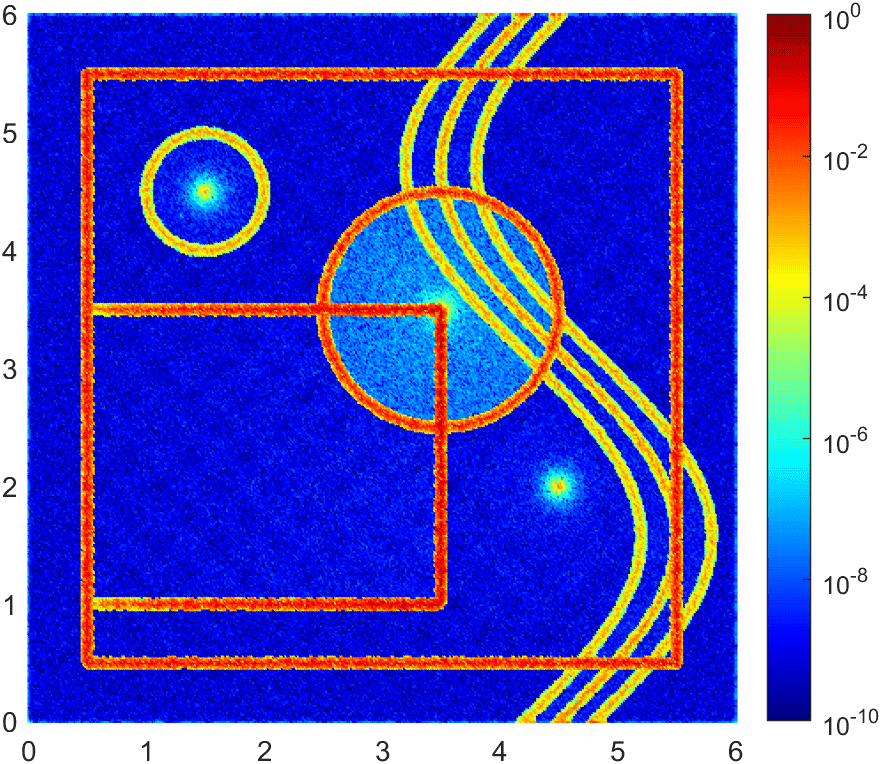}
\put(-10,30){\tiny\rotatebox{90}{direct RBF-FD}}
\put(40,-8){\scriptsize(a) $f$}
\end{overpic}
\begin{overpic}[width=0.32\textwidth]{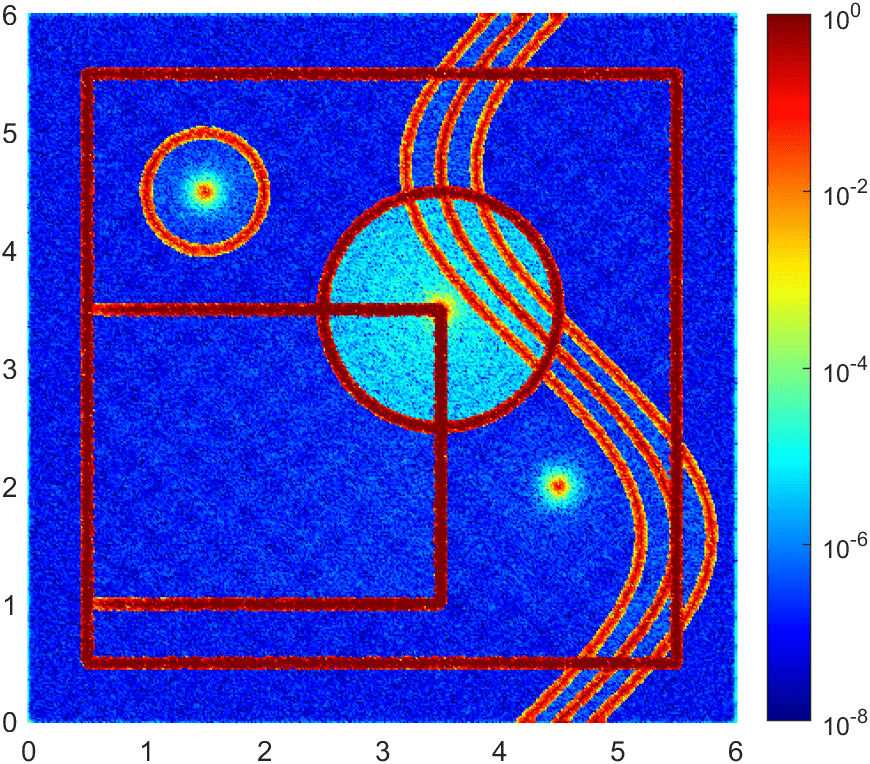}
\put(40,-8){\scriptsize(b) $D_x f$}
\end{overpic}
\begin{overpic}[width=0.32\textwidth]{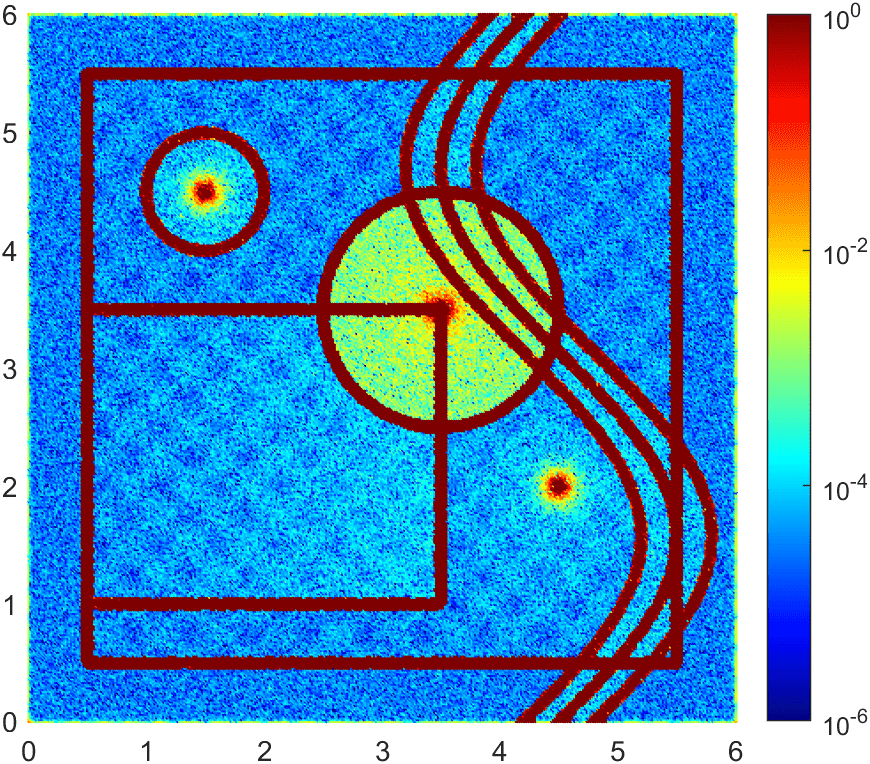}
\put(40,-8){\scriptsize(c) $\Delta f$}
\end{overpic}
\\
\begin{overpic}[width=0.32\textwidth]{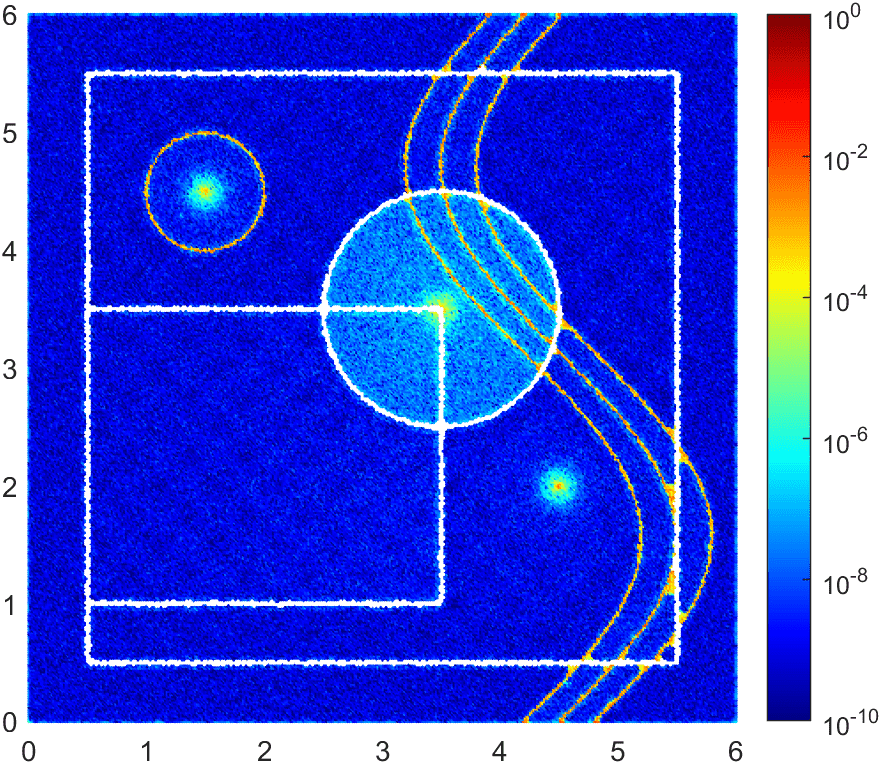}
\put(-10,10){\tiny\rotatebox{90}{regularity-aware RBF-FD}}
\put(40,-8){\scriptsize(a) $f$}
\end{overpic}
\begin{overpic}[width=0.32\textwidth]{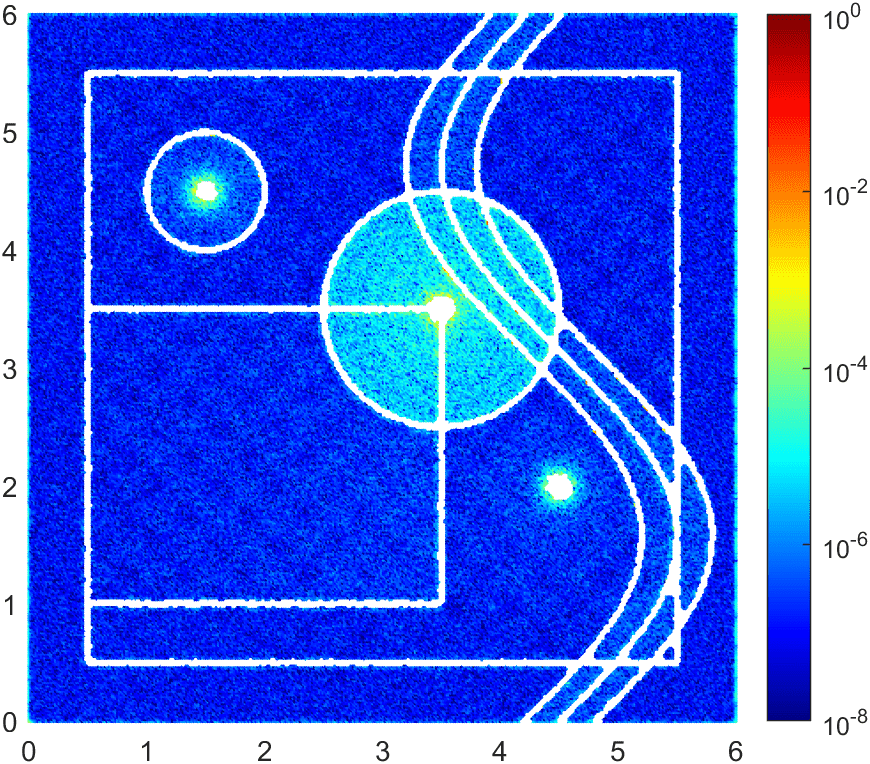}
\put(40,-8){\scriptsize(b) $D_x f$}
\end{overpic}
\begin{overpic}[width=0.32\textwidth]{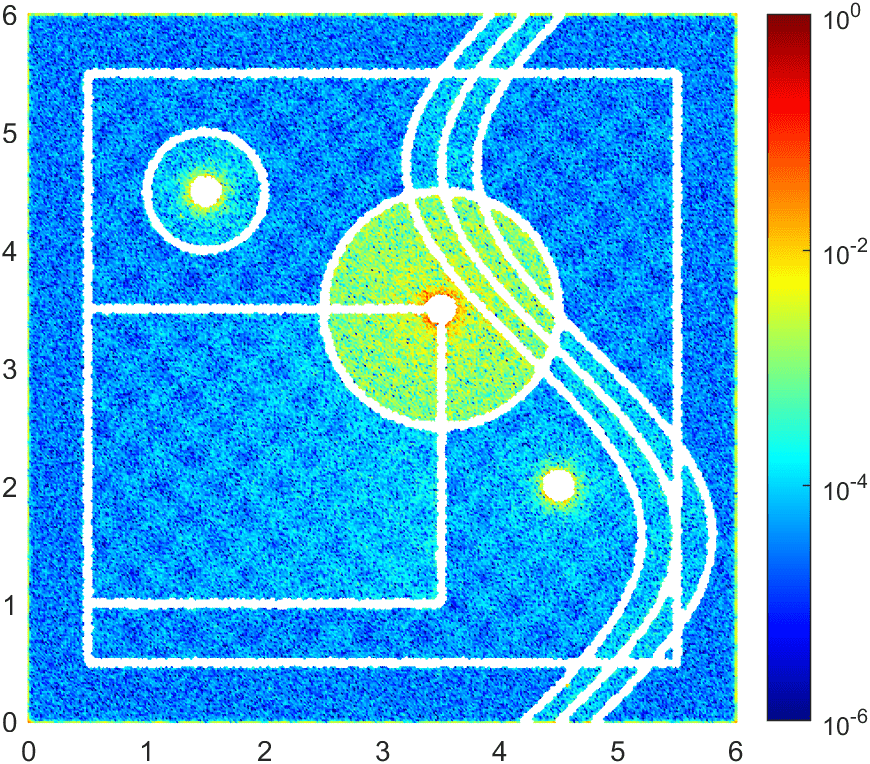}
\put(40,-8){\scriptsize(c) $\Delta f$}
\end{overpic}
\caption{
\textbf{(Example 3)} Pointwise absolute error maps for RBF-FD approximation on
\(\Omega=[0,6]^2\) at \(160{,}000\) evaluation points, computed with \(m=5\) and \(n_X=50\).
The three columns correspond to interpolation of \(f\), approximation of \(D_x f\),
and approximation of \(\Delta f\), respectively.
The top row uses direct RBF-FD stencils, whereas the bottom row uses
regularity-aware stencil shifting, with detected low-regularity locations masked out.
}
\label{fig:error_maps}
\end{figure}

Figure~\ref{fig:error_maps} presents the pointwise absolute error maps for the interpolation and differentiation results.
The three columns correspond to the reconstructed function \(f\), its first derivative \(D_x f\), and the Laplacian \(\Delta f\), respectively.
The top row shows the direct RBF-FD results, while the bottom row shows the results obtained with regularity-aware stencil shifting.
Errors remain small throughout the smooth regions, while larger values concentrate near the singular features of $f$.
The white contours indicate areas classified as low‑regularity, where finite‑difference evaluations are intentionally omitted.
These results confirm that the regularity‑aware refinement preserves accuracy in smooth portions of the domain and correctly isolates regions of limited differentiability.

\subsection{Example 4: Regularity Mapping and Differentiation in Turbulent Flow Data}
\label{subsec:turbulent-flow}

We demonstrate the applicability of the proposed framework to a real, non-uniform turbulent flow dataset.
The test case is taken from the incompressible-flow benchmark described in Section~2.13 of~\cite{hecht2012new},
comprising $104,\!444$ scattered nodes generated by FreeFEM++.
The sample distribution is highly irregular due to adaptive meshing in a complex geometry,
providing a stringent test of robustness.
Throughout this experiment, the local stencil size is fixed at $n=20$.

\begin{figure}[t]
\centering
\subfloat[PDE solution field]{\includegraphics[width=0.31\textwidth]{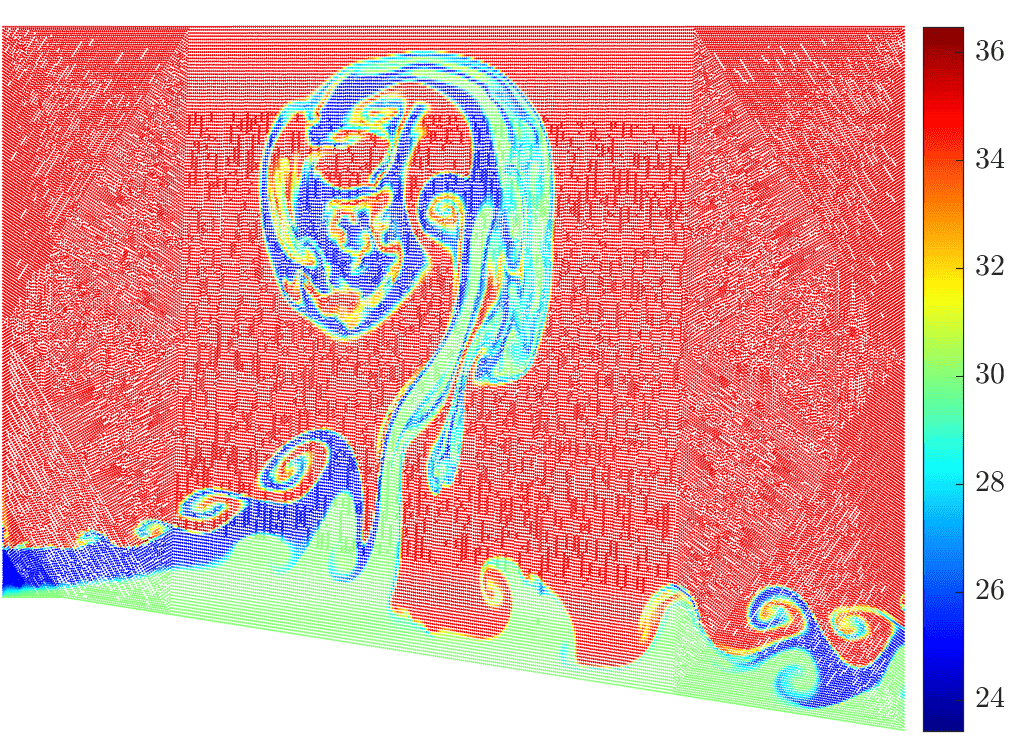}}
\hspace{1mm}
\subfloat[Native norm]{\includegraphics[width=0.31\textwidth]{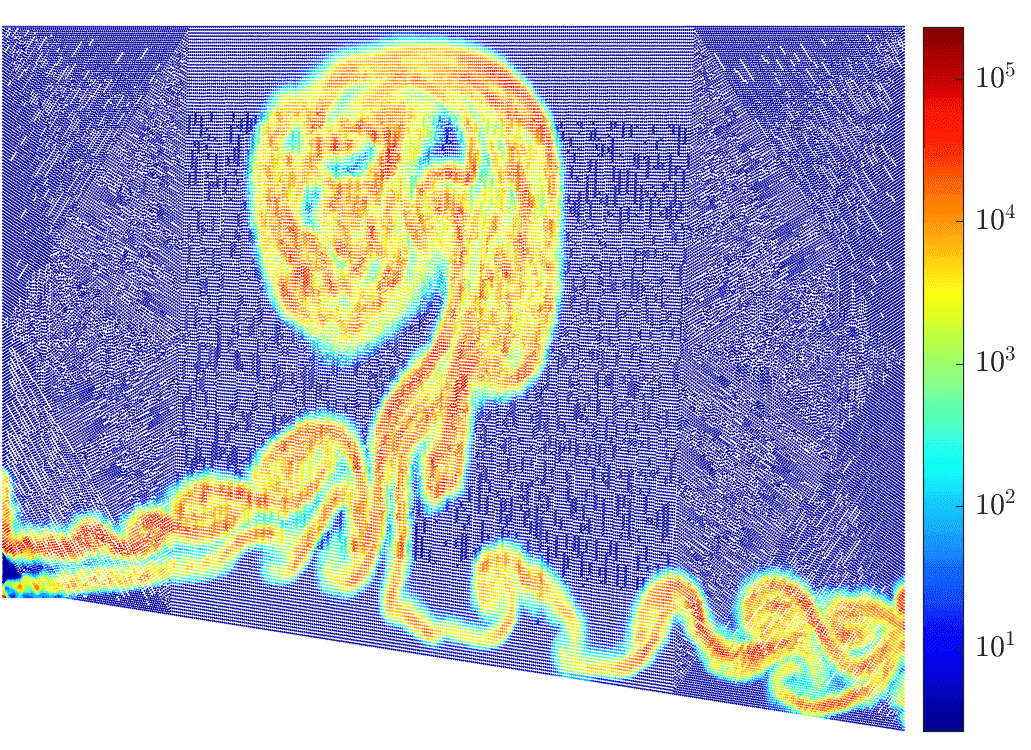}}
\hspace{1mm}
\subfloat[ $\tilde{s}(\Omega_z)$ with stencil-shift]{\includegraphics[width=0.31\textwidth]{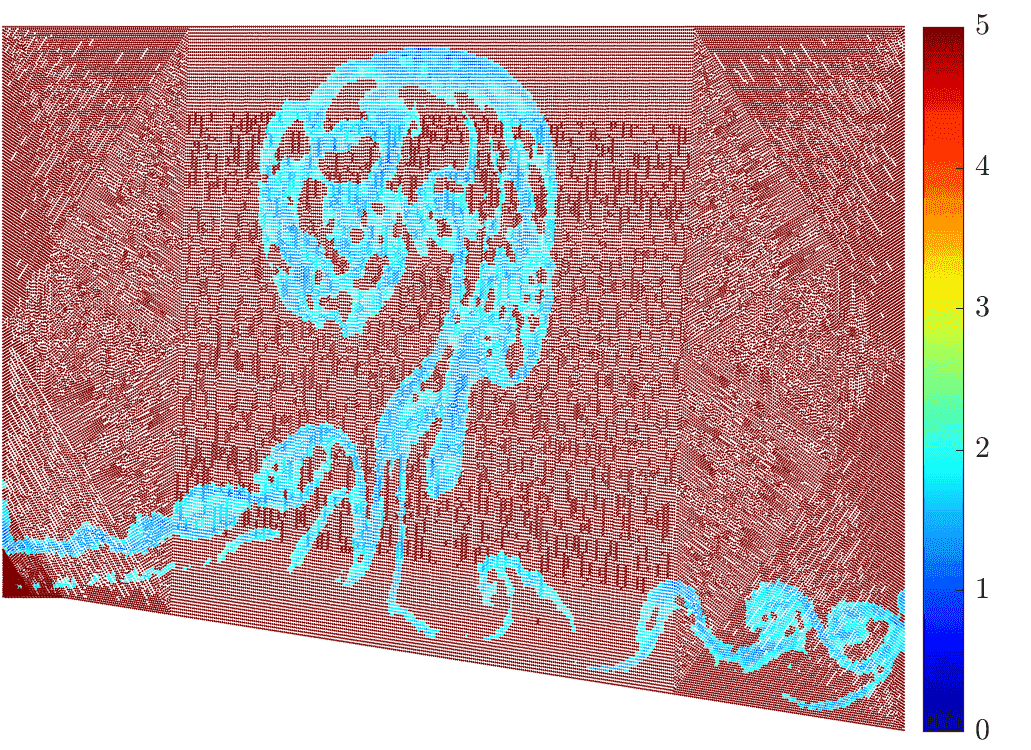}}
\caption{
\textbf{(Example 4)} Turbulent-flow dataset on $104,\!444$ scattered nodes. (b) Native-norm with $21,\!847$ outliers identified by IQR-based screening.
}
\label{fig:turbulent_screening}
\end{figure}

Figure~\ref{fig:turbulent_screening} illustrates the screening and subsequent regularity analysis.
The figure displays (a) the PDE solution field, (b) the native‑norm, and (c) the data‑driven regularity distribution following adaptive stencil refinement.
Irregular spatial sampling leads to elevated native norms near regions of strong gradients,
while smooth areas remain clearly separated.
After adaptive refinement, the data-driven regularity becomes spatially consistent across the domain.

\begin{figure}[t]
\centering
\subfloat[$|\alpha|=2$ (21,812)]{\includegraphics[width=0.31\textwidth]{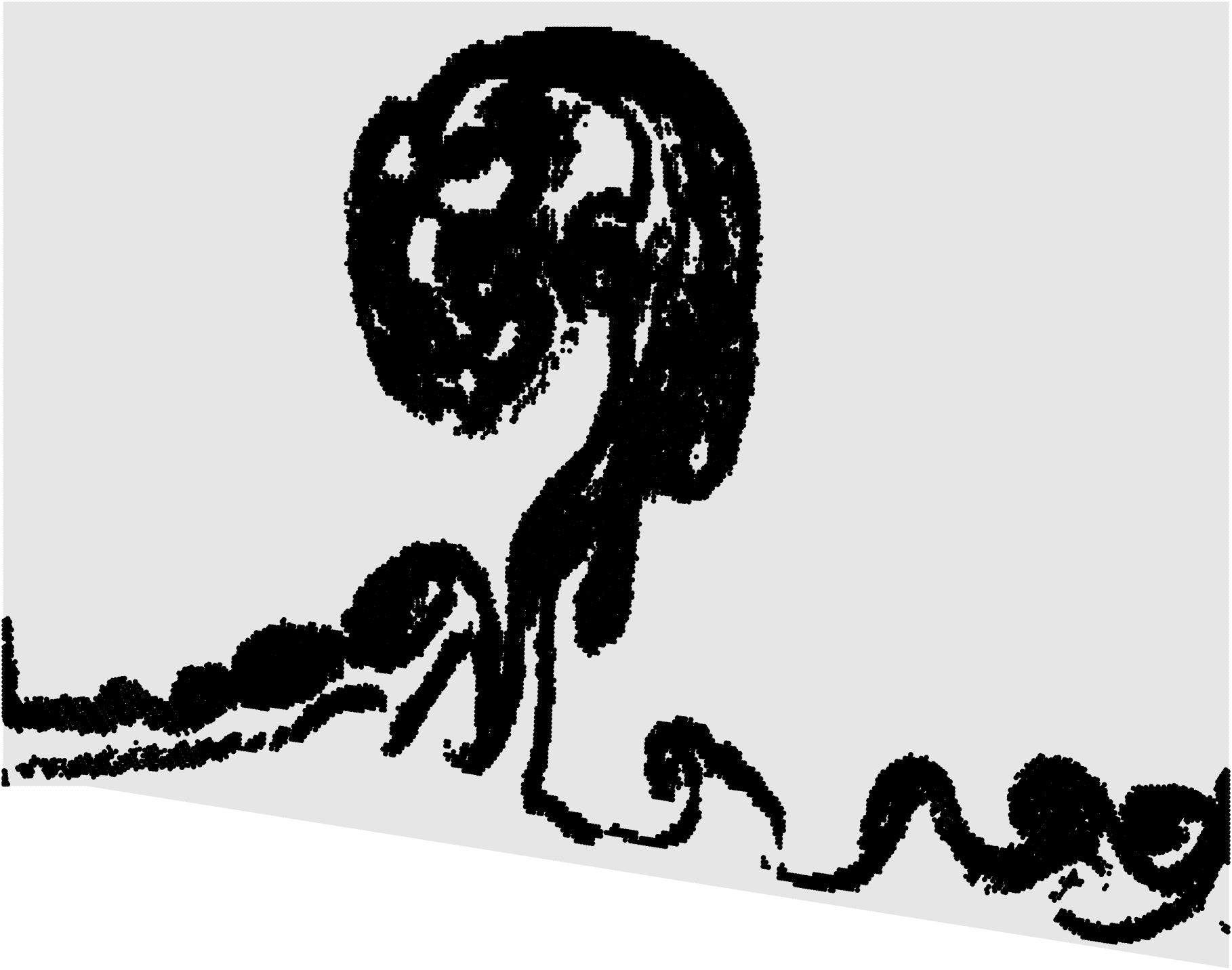}}
\hspace{1mm}
\subfloat[$|\alpha|=1$ (18,237)]{\includegraphics[width=0.31\textwidth]{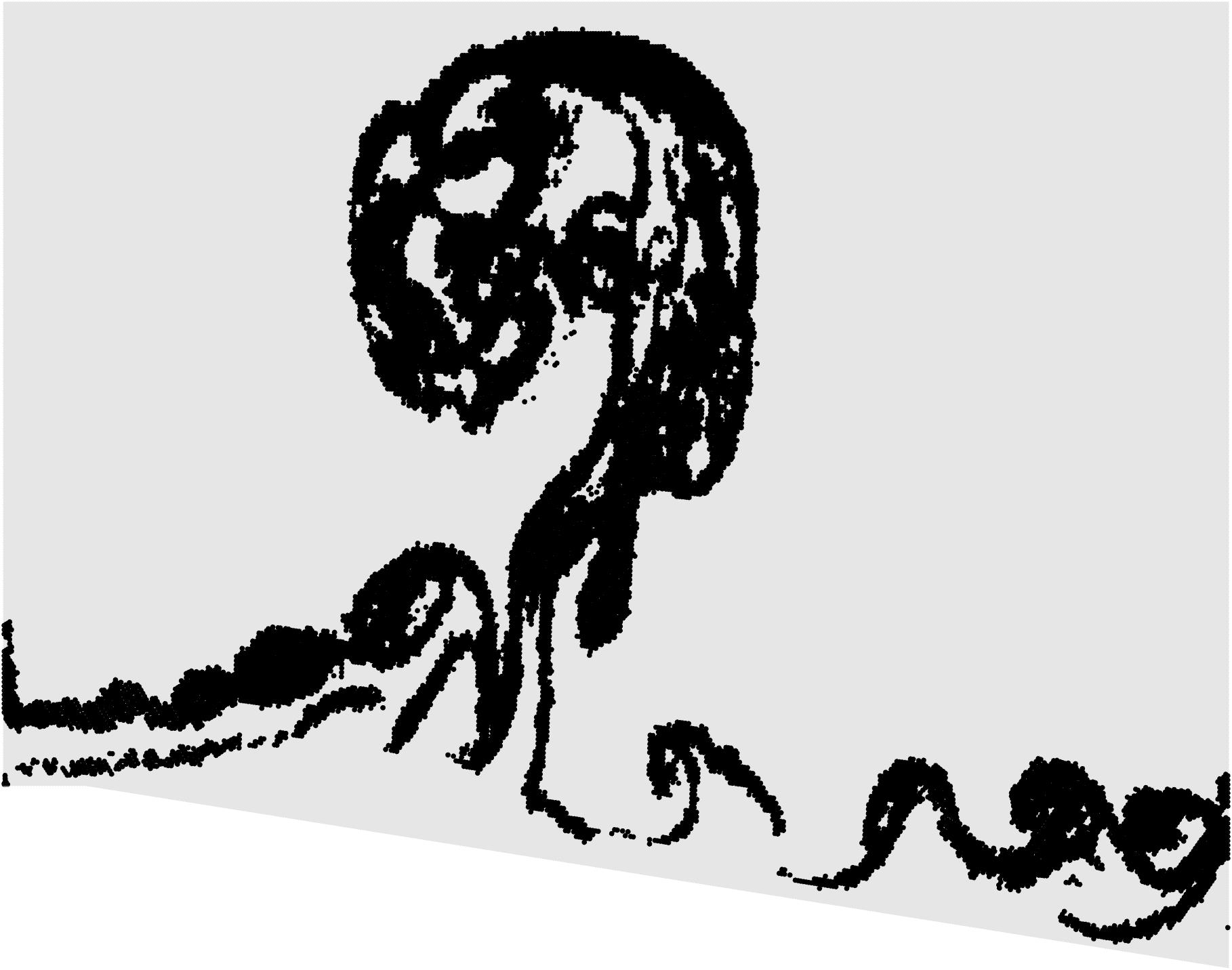}}
\hspace{1mm}
\subfloat[$|\alpha|=0$ (4,385)]{\includegraphics[width=0.31\textwidth]{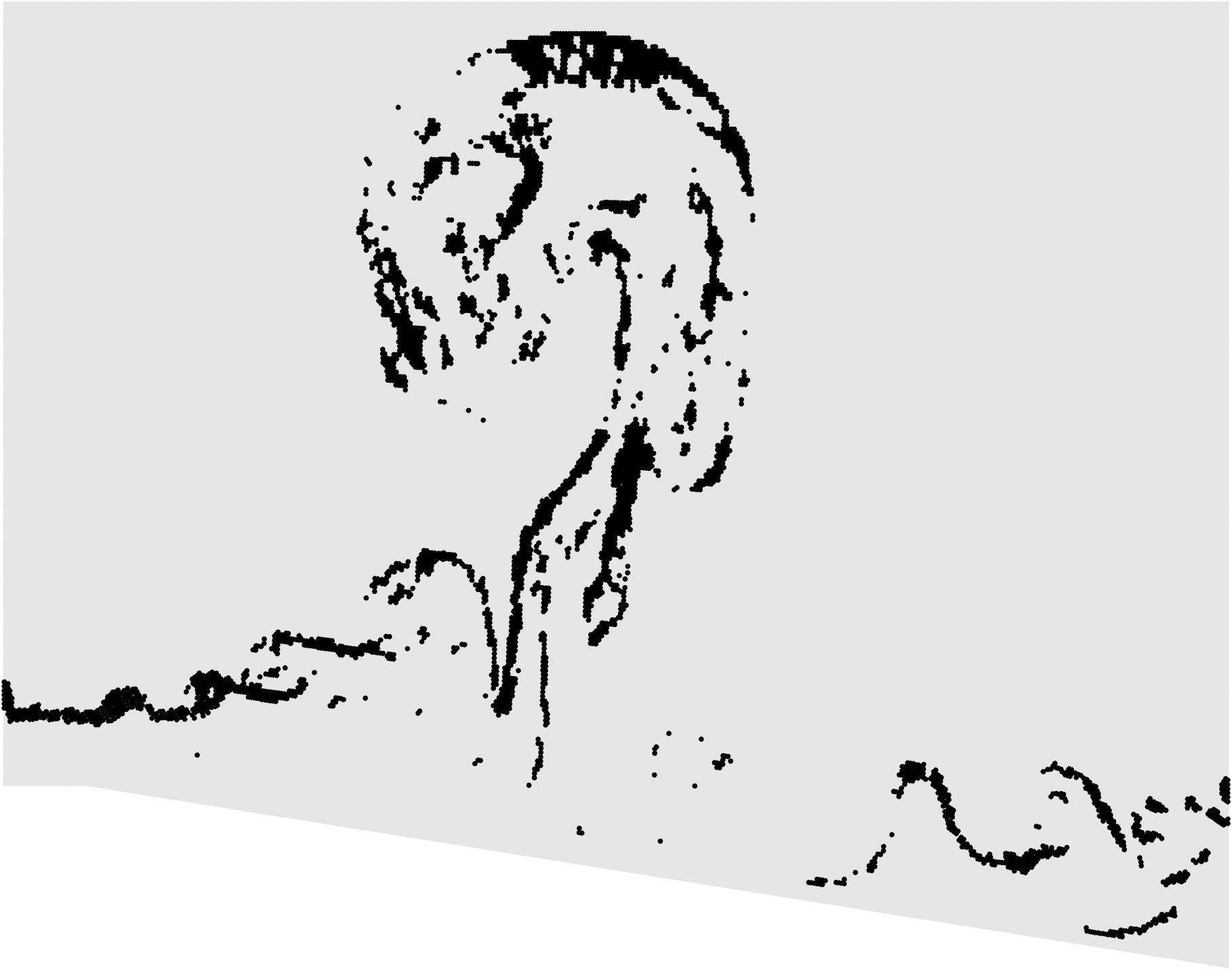}}
\caption{
\textbf{(Example 4)} Flagged points under the regularity-limited criterion
$\tilde{s}(\Omega_z)<|\alpha|+d/2$ for derivative orders $|\alpha|=2$,
$|\alpha|=1$ , and $|\alpha|=0$,
using adaptive RBF-FD stencils ($n=20$).
}
\label{fig:turbpdesingularity_detection}
\end{figure}

Using the discrete regularity criterion $\tilde{s}(\Omega_z)<|\alpha|+d/2$, flagged points are classified by the
highest derivative order locally supported.
Figure~\ref{fig:turbpdesingularity_detection}(a)–(c) show the distributions for
$|\alpha|=2$, $1$, and $0$, respectively.
The identified regions align with coherent physical structures where flow gradients are steep,
while noise and spurious detections in smooth areas are largely suppressed.

\begin{figure}[t]
\centering
\subfloat[FEM reference $D_x f$ with \\ $P1$ element]{\includegraphics[width=0.32\textwidth]{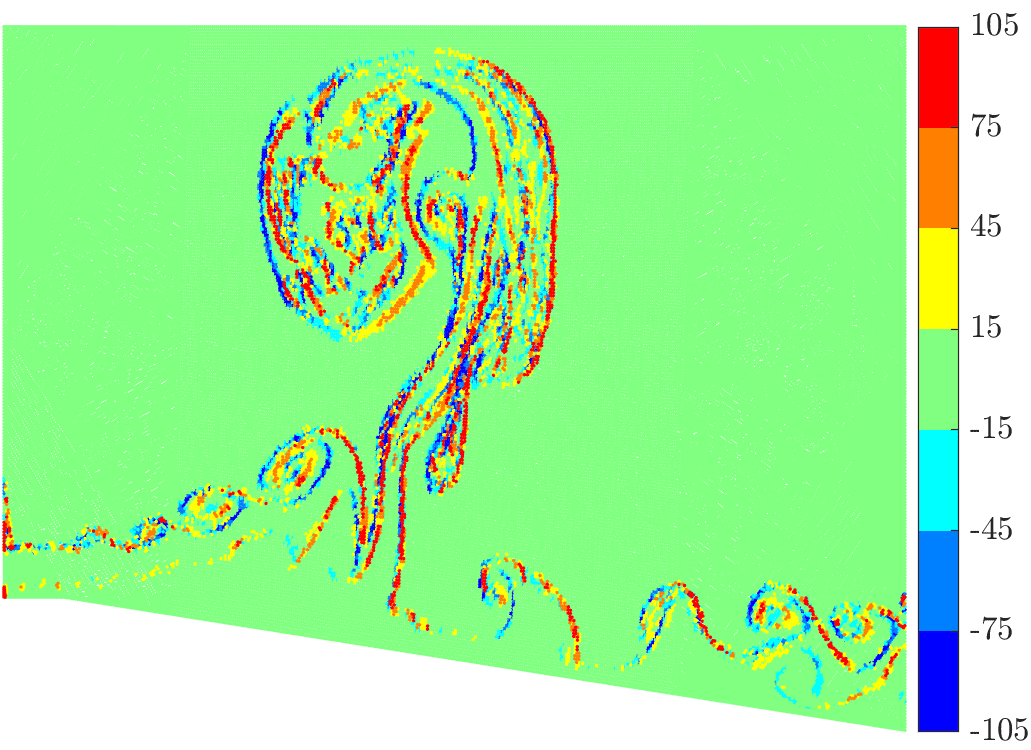}}
\hspace{1mm}
\subfloat[RBF-FD approximation \\ using regularity-aware stencils]{\includegraphics[width=0.32\textwidth]{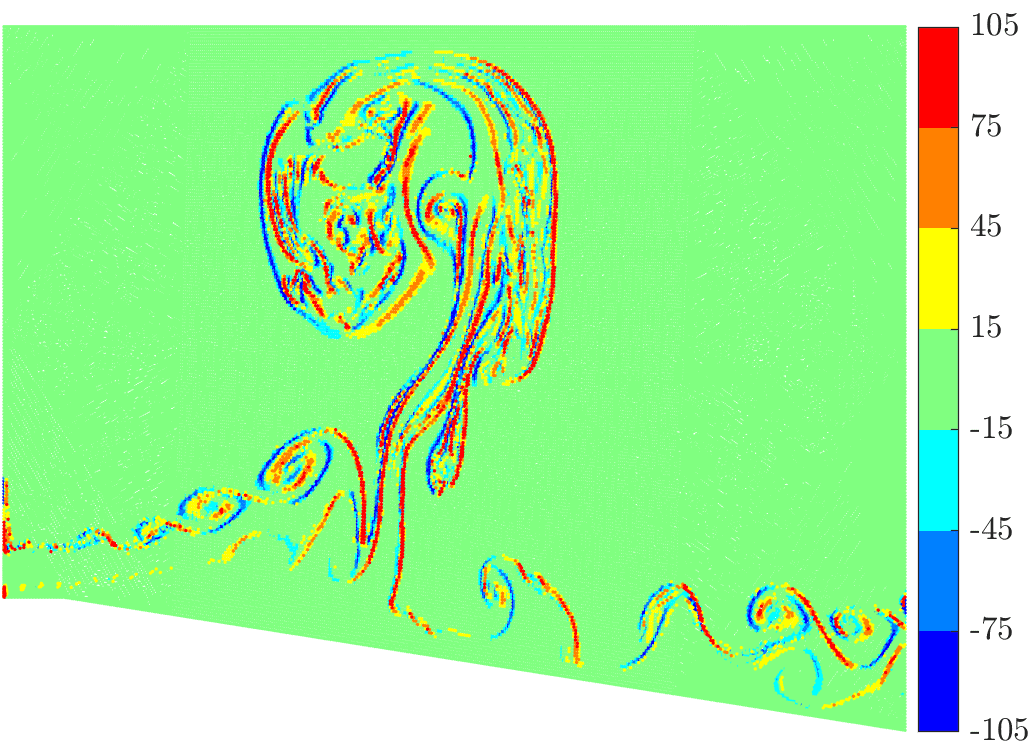}}
\hspace{1mm}
\subfloat[pointwise  relative absolute error]{\includegraphics[width=0.32\textwidth]{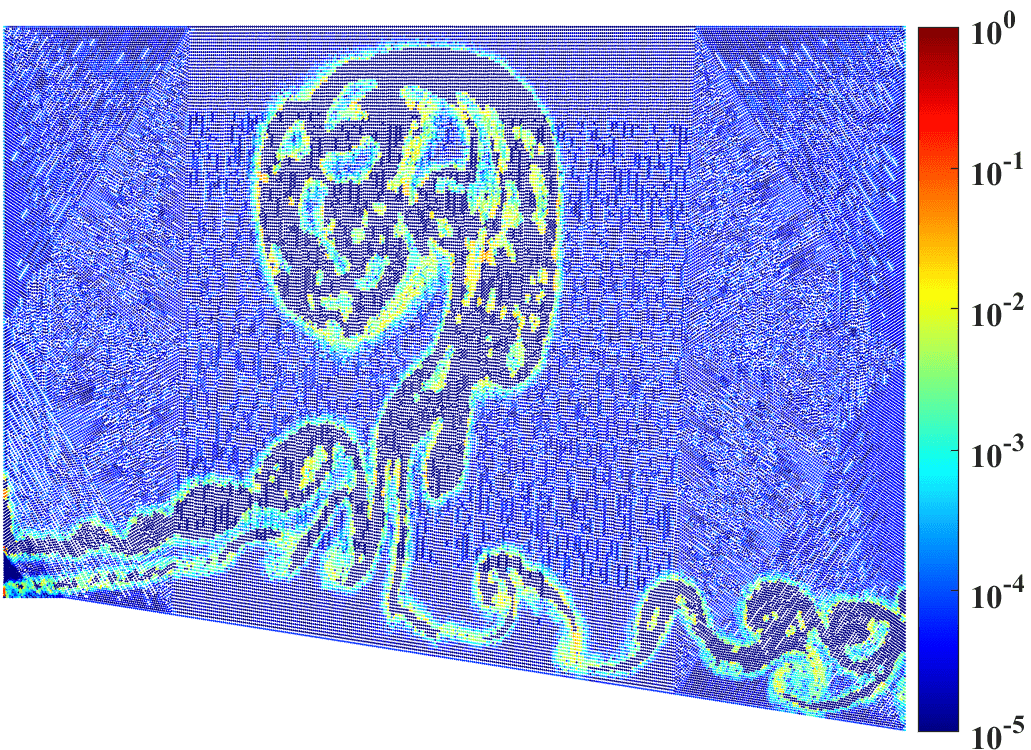}}
\caption{
\textbf{(Example 4)} Derivative validation for the turbulent dataset.
}
\label{fig:turbpdederiv_approx}
\end{figure}

Figure~\ref{fig:turbpdederiv_approx} compares the FEM reference $D_xf$, the RBF-FD approximation using
regularity-aware stencils, and the pointwise relative absolute error.
The main discrepancies occur near low-regularity regions and boundaries.

Near boundaries or strongly irregular regions, high-order finite-difference operators
can become ill-conditioned.
Following the reliability map, we employ low-order stencils in such regions to maintain stability.
More sophisticated treatments, such as adaptive viscosity~\cite{bracco2025positive},
least-squares RBF-FD~\cite{tominec2021least}, or refinement by null rules~\cite{bracco2023discontinuity},
represent natural extensions.
This example demonstrates that the proposed regularity mapping enables
stable and accurate differentiation even in non-uniform, physically complex datasets.


\begin{figure}[t]
\centering
\begin{overpic}[width=0.32\textwidth]{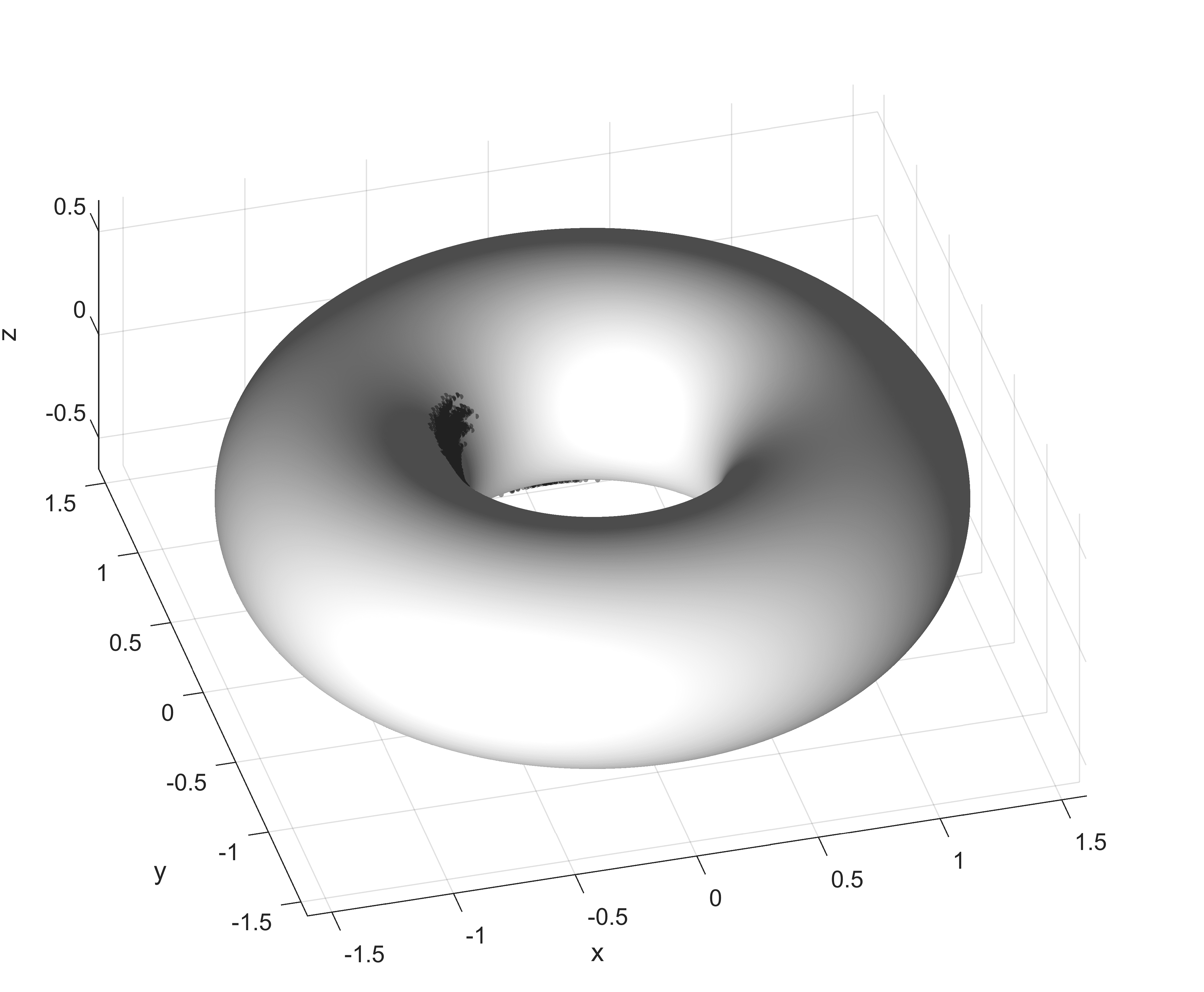}
\put(38,-8){\scriptsize(a) $t\approx 1.0$}
\put(-10,30){\tiny\rotatebox{90}{detection}}
\end{overpic}
\begin{overpic}[width=0.32\textwidth]{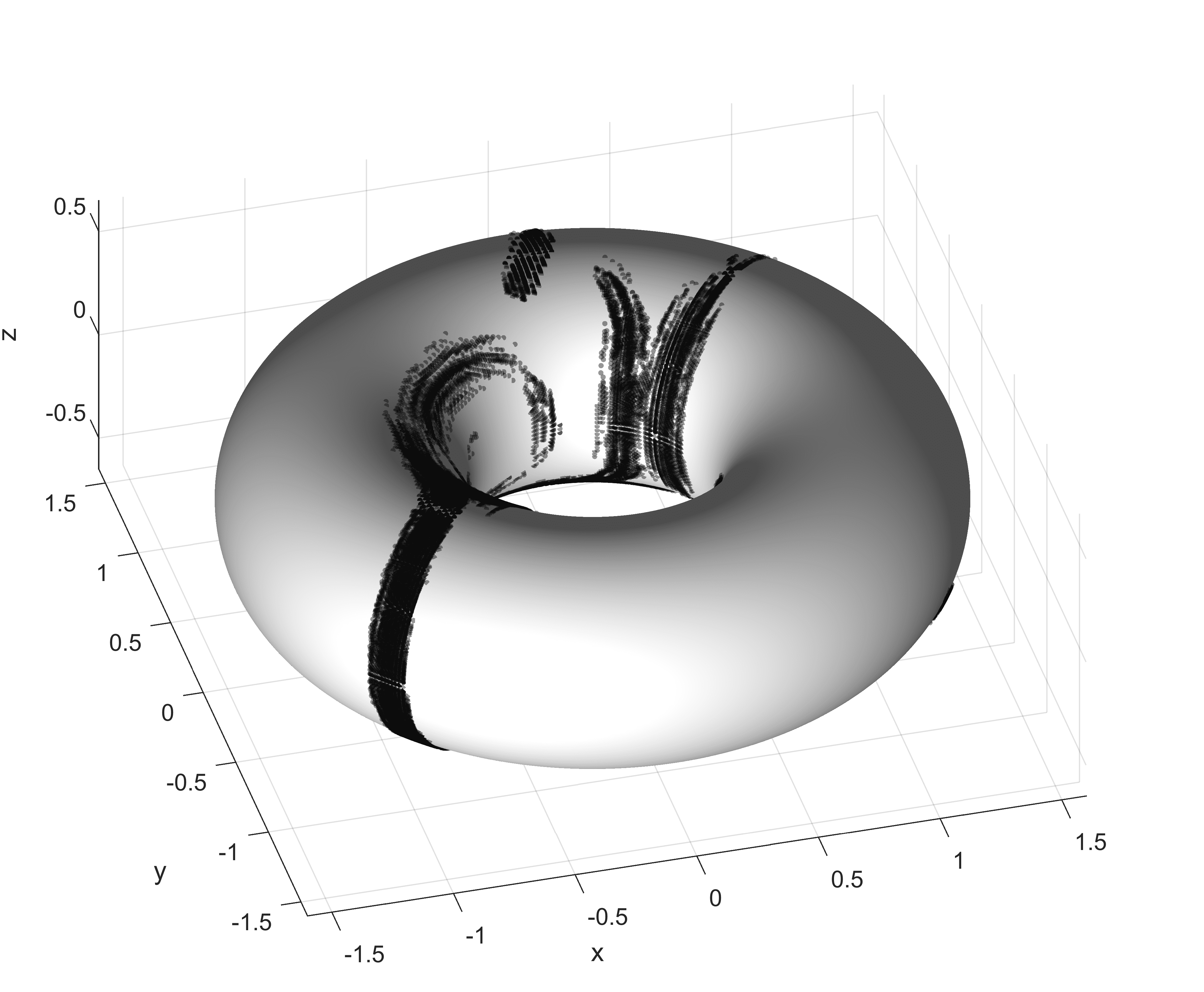}
\put(38,-8){\scriptsize(b) $t\approx 5.0$}
\end{overpic}
\begin{overpic}[width=0.32\textwidth]{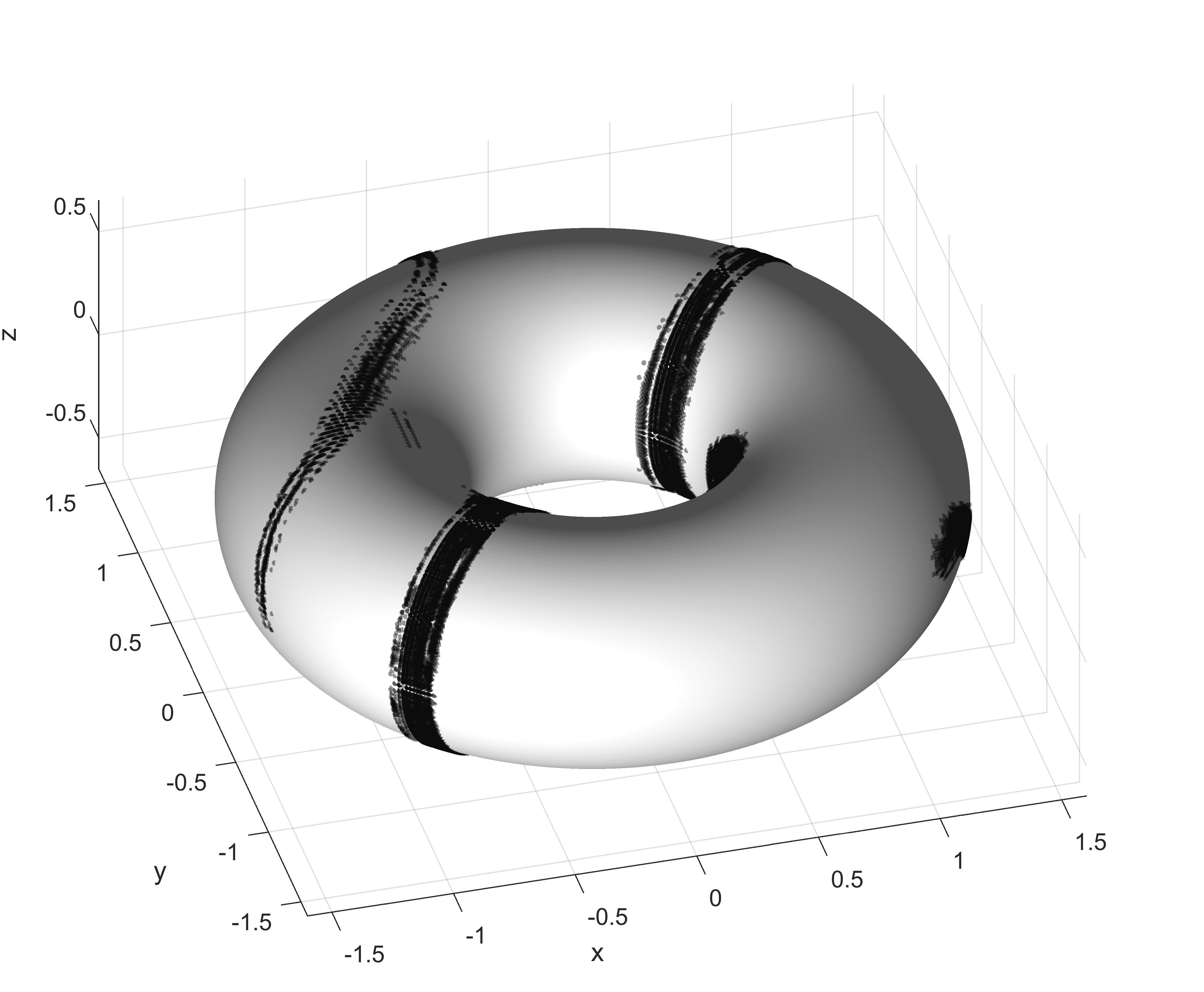}
\put(38,-8){\scriptsize(c) $t\approx 10.0$}
\end{overpic}
\\
\begin{overpic}[width=0.32\textwidth]{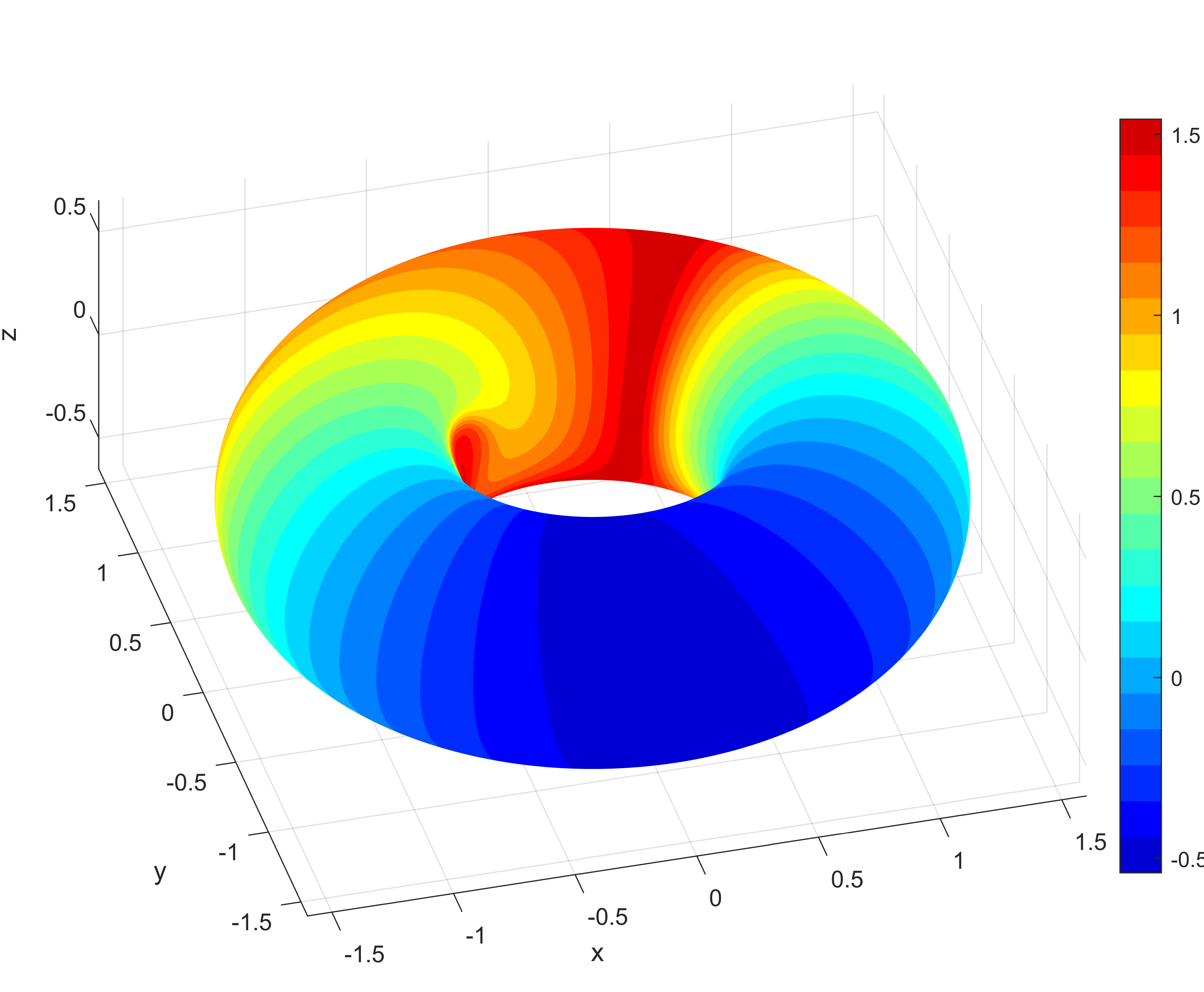}
\put(38,-8){\scriptsize(d) $t\approx 1.0$}
\put(-10,30){\tiny\rotatebox{90}{$u$ on $\mathbb{T}^2 $}}
\end{overpic}
\begin{overpic}[width=0.32\textwidth]{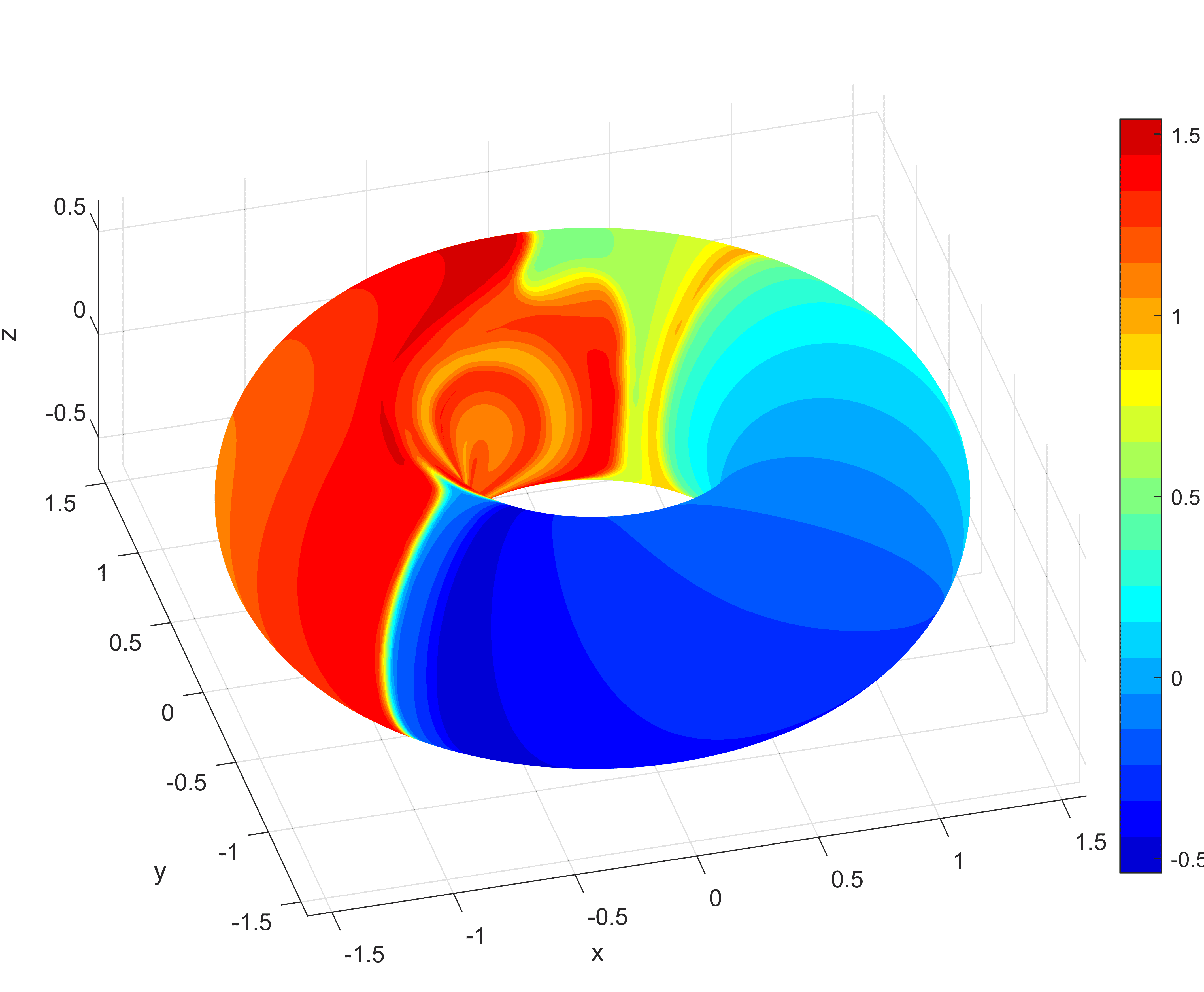}
\put(38,-8){\scriptsize(e) $t\approx 5.0$}
\end{overpic}
\begin{overpic}[width=0.32\textwidth]{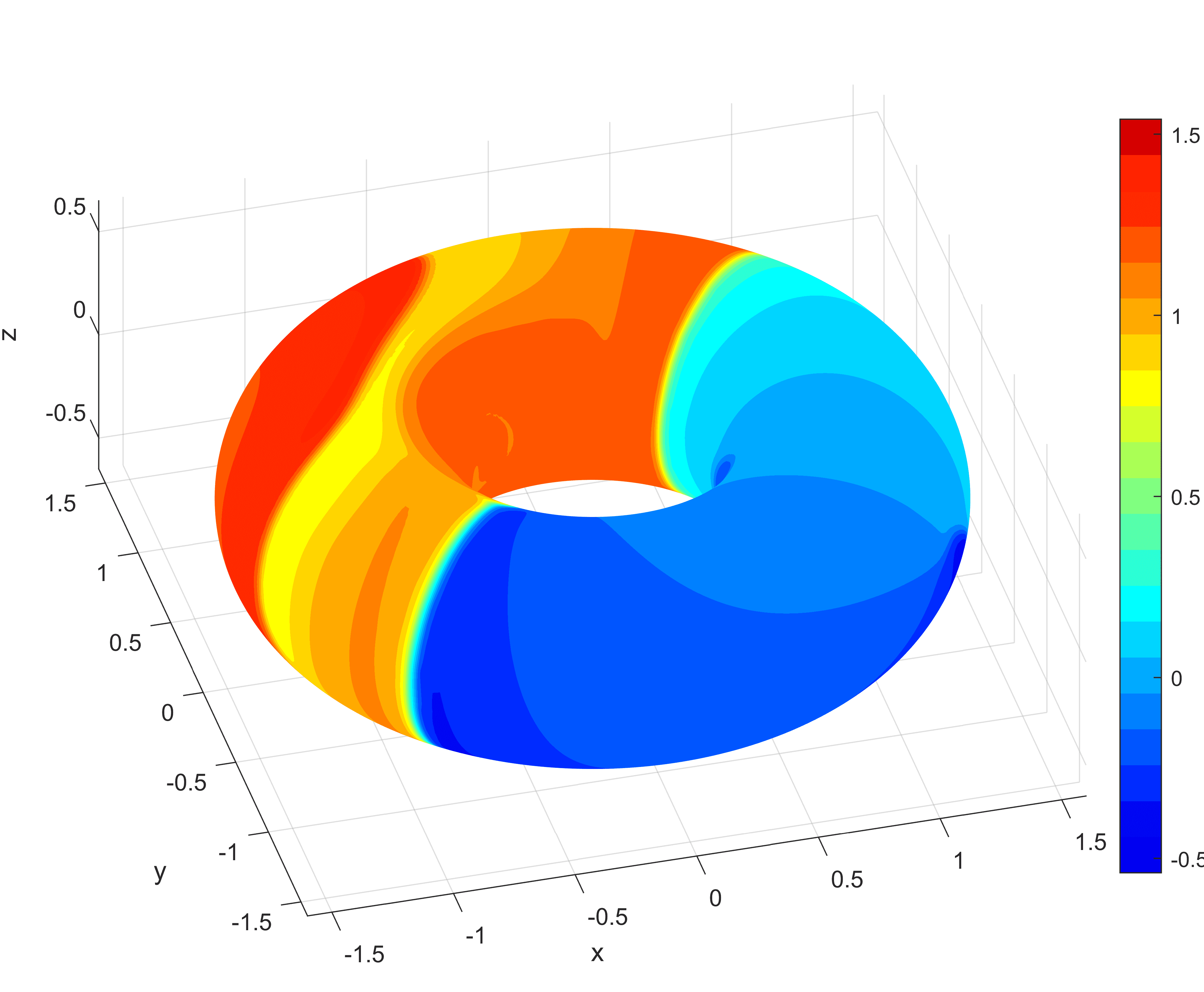}
\put(38,-8){\scriptsize(f) $t\approx 10.0$}
\end{overpic}
\caption{
\textbf{(Example 5) Regularity detection for the Burgers equation on a torus.}
From left to right, the snapshots correspond to $t\approx1.0$, $t\approx5.0$,
and $t\approx10.0$. The first row shows the points detected as locally
low-regularity regions on the torus, and the second row shows the numerical
solution $u$ on the surface $\mathbb{T}^2 $. The initial data are smooth; hence the
detected low-regularity bands are generated dynamically by nonlinear steepening
during the PDE evolution. The marked points form the adaptive mask used in
the closest-point extension step.
}
\label{fig:burgers_torus_regularity_detection}
\end{figure}

\begin{figure}[t]
\centering
\begin{overpic}[width=0.32\textwidth]{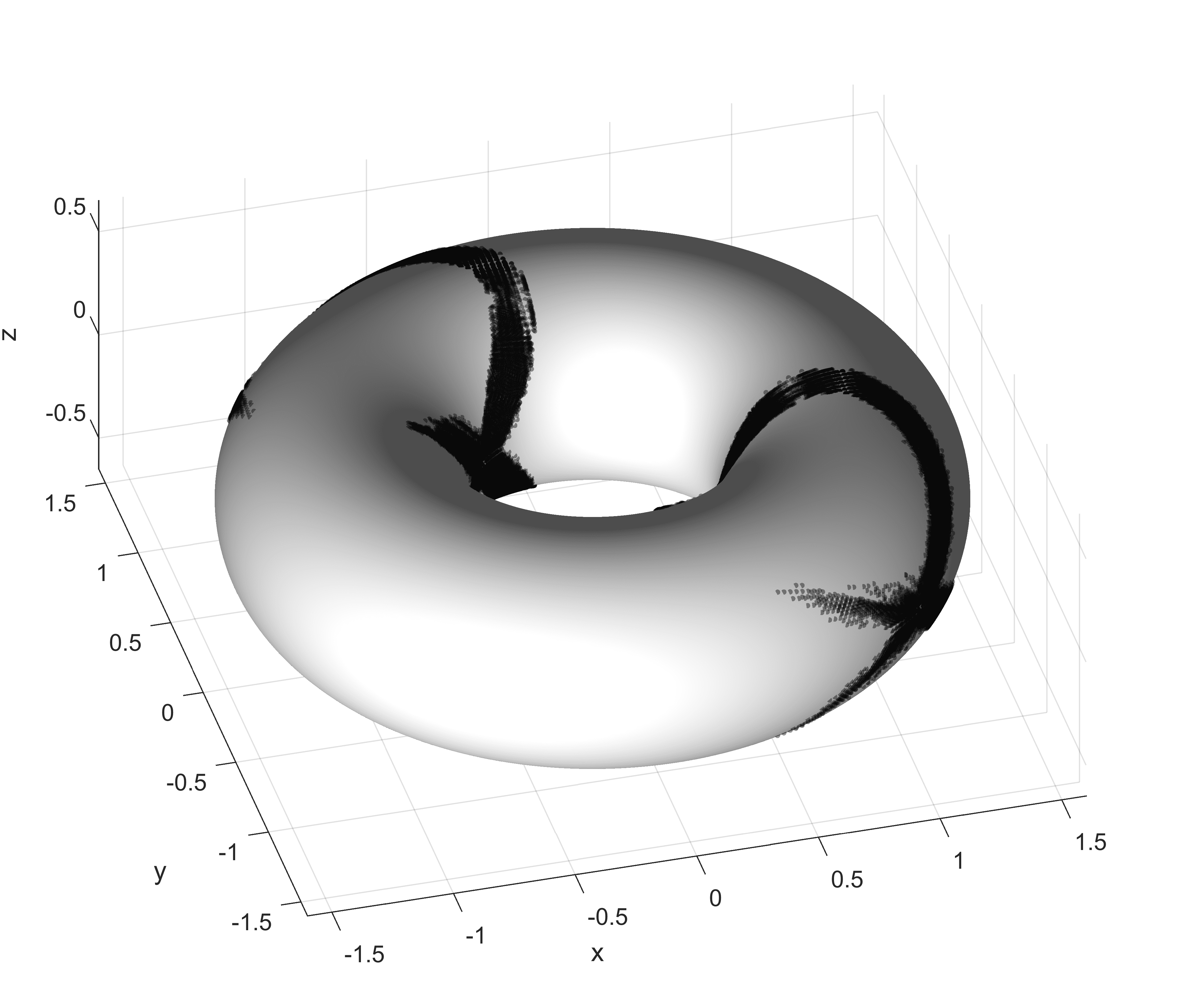}
\put(38,-8){\scriptsize(a) $t\approx 0.1$}
\put(-10,30){\tiny\rotatebox{90}{detection}}
\end{overpic}
\begin{overpic}[width=0.32\textwidth]{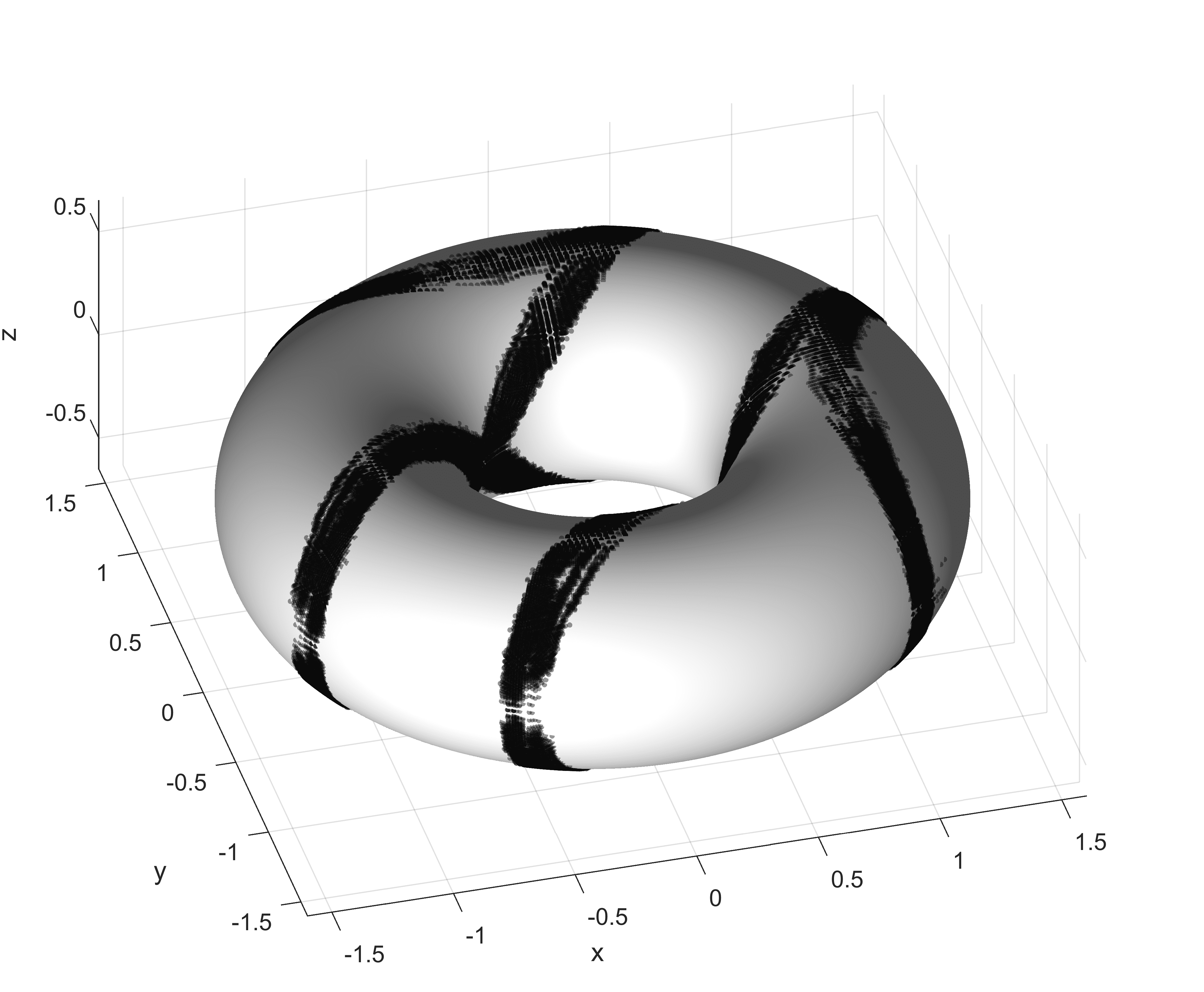}
\put(38,-8){\scriptsize(b) $t\approx 0.5$}
\end{overpic}
\begin{overpic}[width=0.32\textwidth]{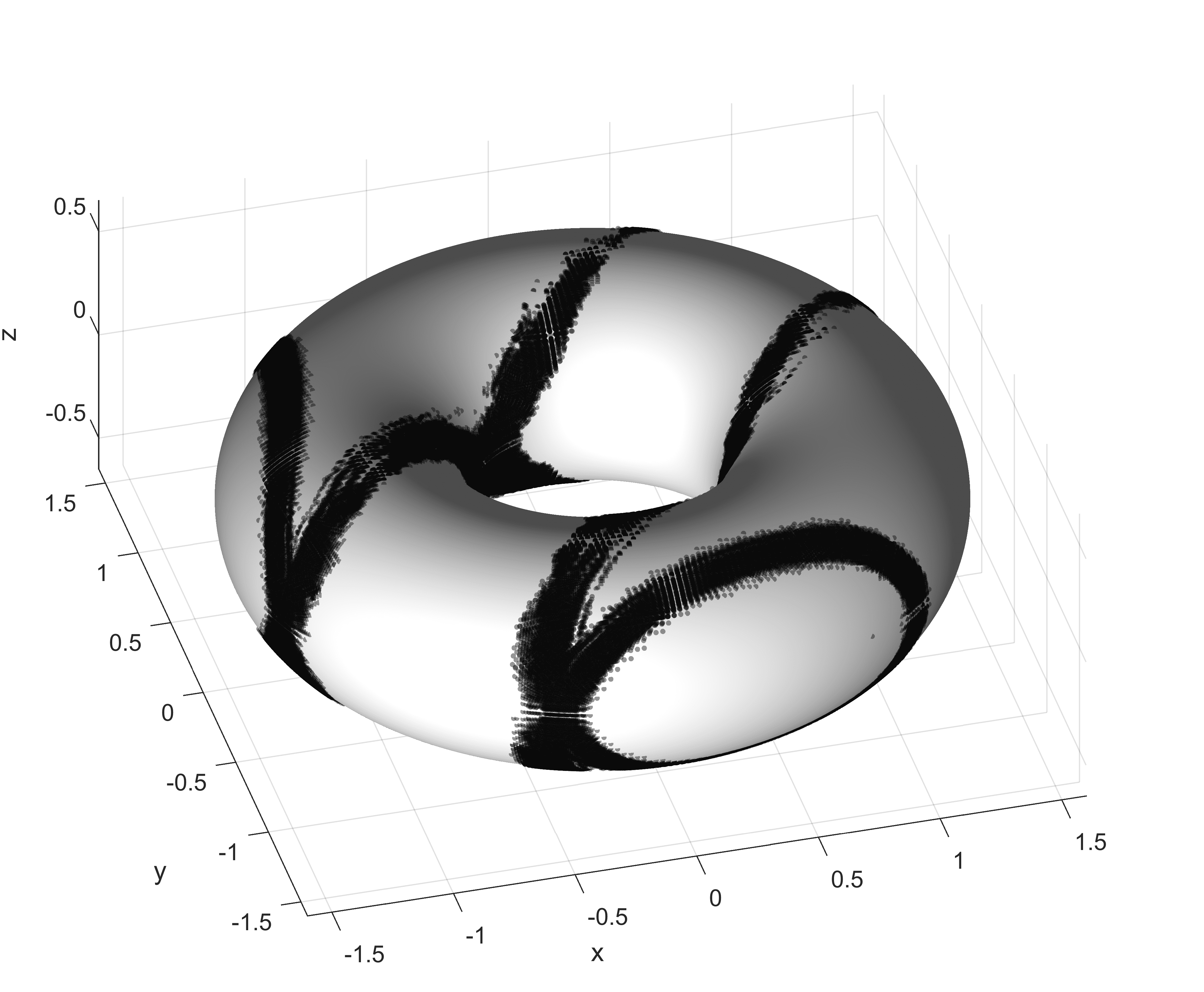}
\put(38,-8){\scriptsize(c) $t\approx 1.0$}
\end{overpic}
\\
\begin{overpic}[width=0.32\textwidth]{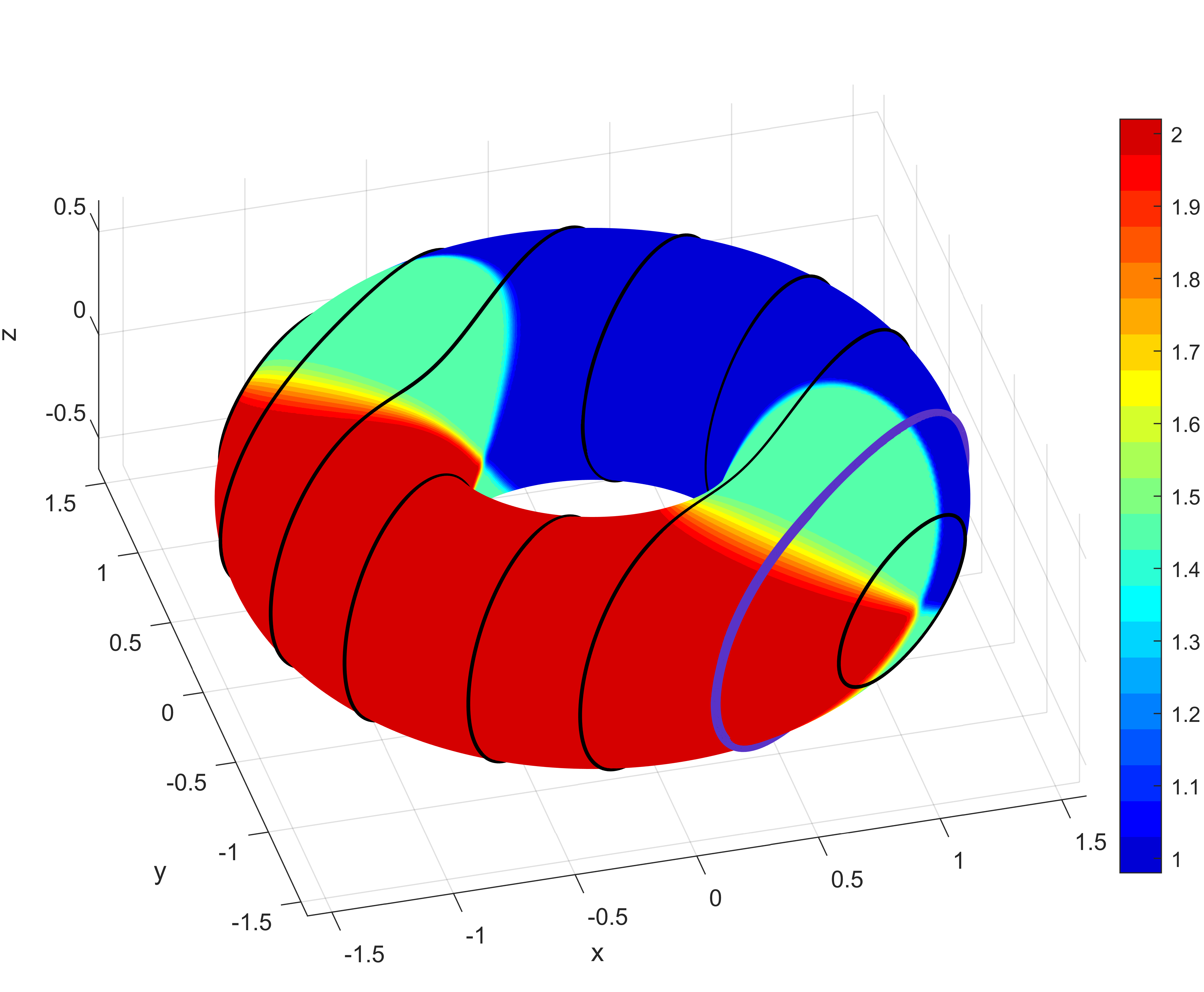}
\put(38,-8){\scriptsize(d) $t\approx 0.1$}
\put(-10,30){\tiny\rotatebox{90}{$h$ on $\mathbb{T}^2 $}}
\end{overpic}
\begin{overpic}[width=0.32\textwidth]{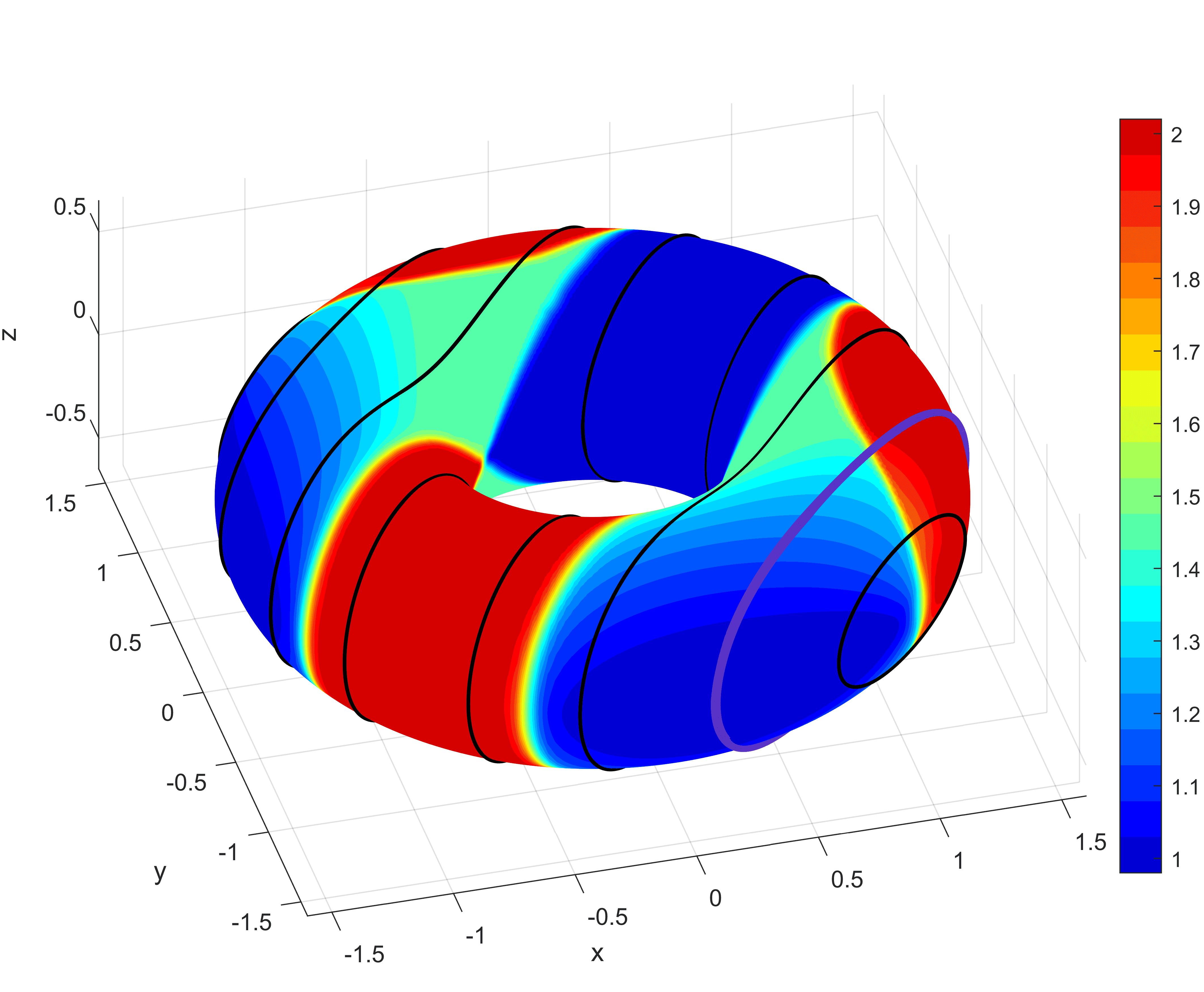}
\put(38,-8){\scriptsize(e) $t\approx 0.5$}
\end{overpic}
\begin{overpic}[width=0.32\textwidth]{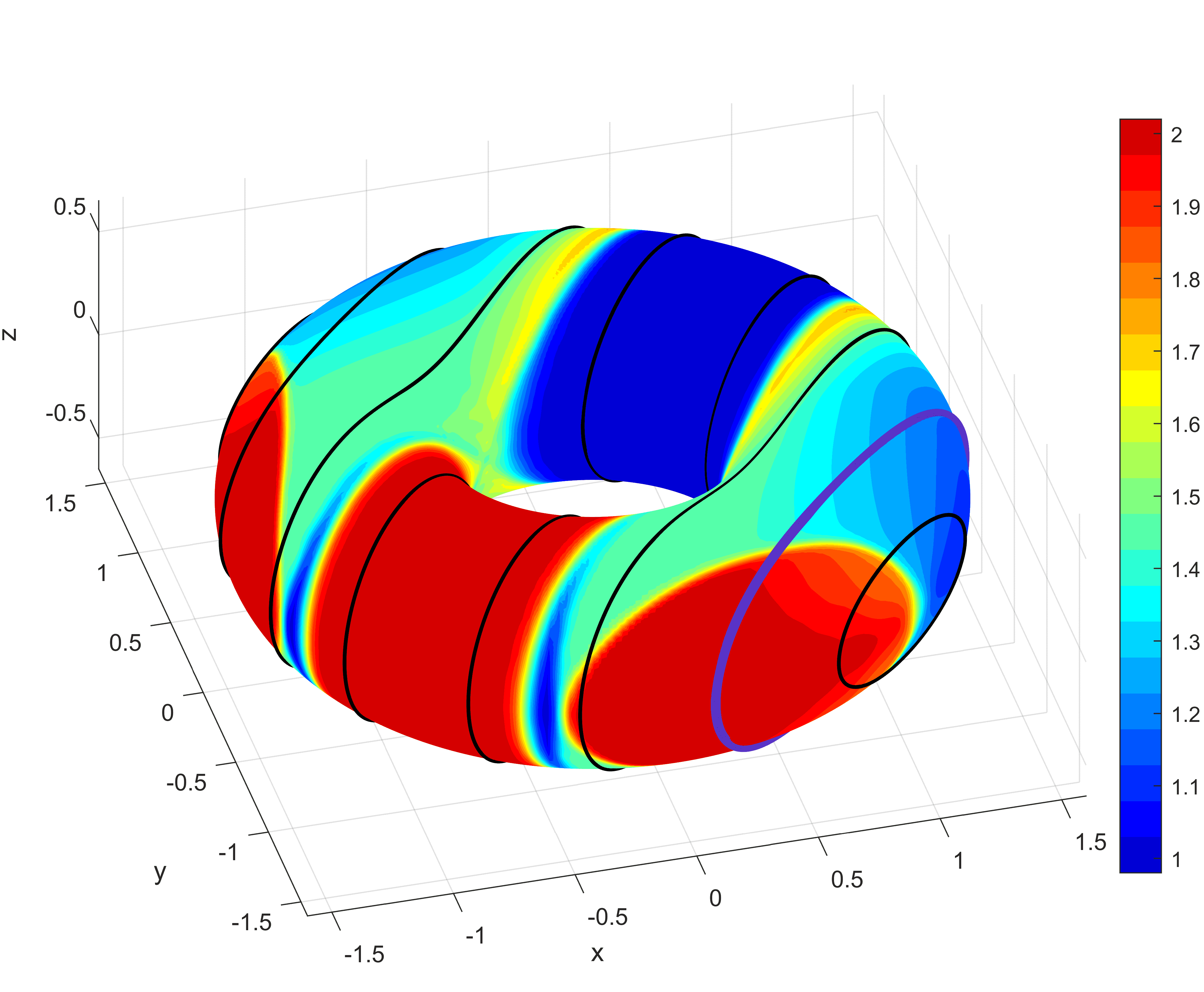}
\put(38,-8){\scriptsize(f) $t\approx 1.0$}
\end{overpic}
\\
\begin{overpic}[width=0.32\textwidth]{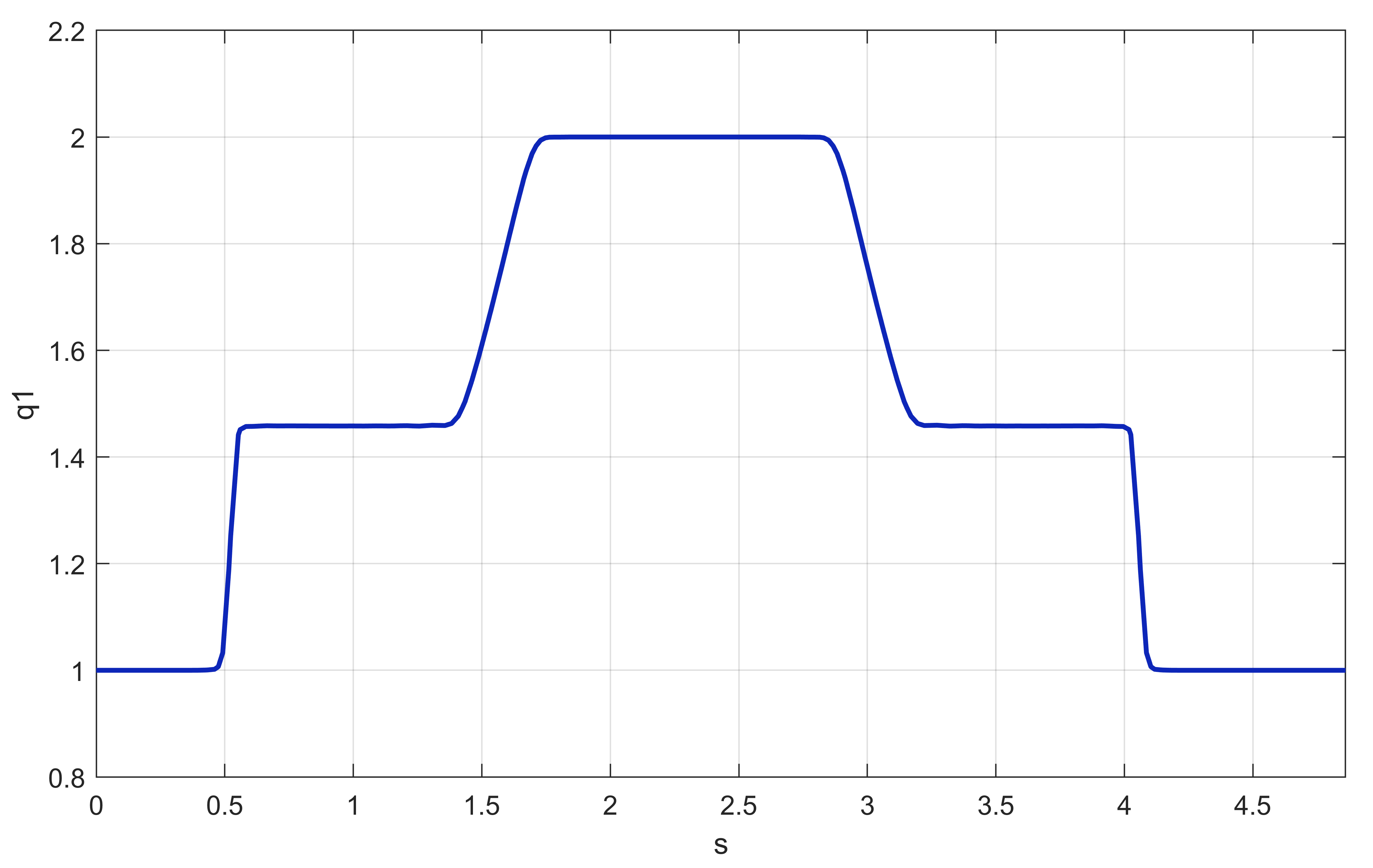}
\put(38,-8){\scriptsize(g) $t\approx 0.1$}
\put(-10,30){\tiny\rotatebox{90}{profile}}
\end{overpic}
\begin{overpic}[width=0.32\textwidth]{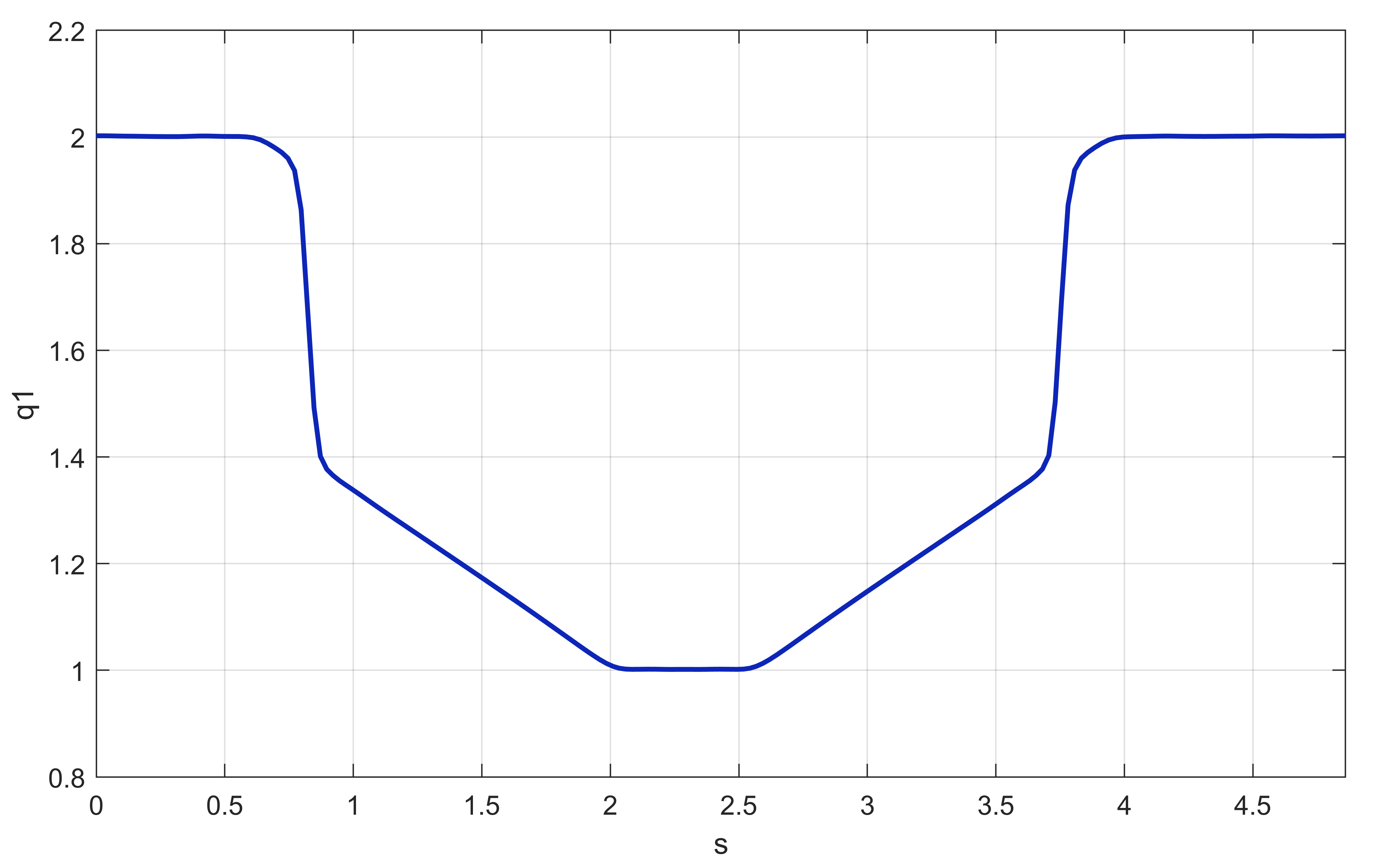}
\put(38,-8){\scriptsize(h) $t\approx 0.5$}
\end{overpic}
\begin{overpic}[width=0.32\textwidth]{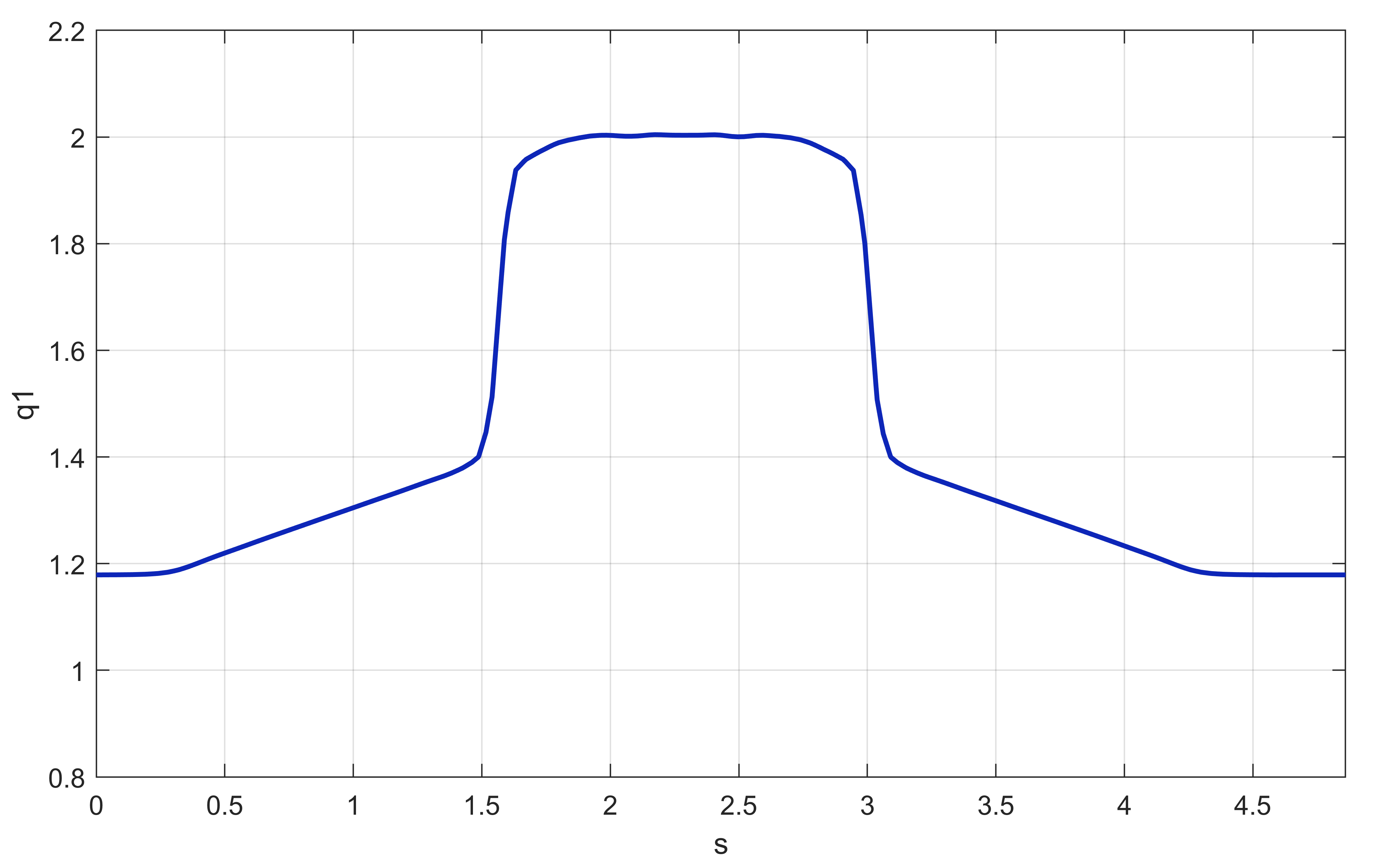}
\put(38,-8){\scriptsize(i) $t\approx 1.0$}
\end{overpic}
\caption{
\textbf{(Example 5) Regularity detection for the reduced shallow-water-type system on a torus.}
From left to right, the snapshots correspond to $t\approx0.1$, $t\approx0.5$,
and $t\approx1.0$. The first row shows the points detected as locally
low-regularity regions on the torus, the second row shows the numerical
solution $h$ on the surface $\mathbb{T}^2 $, and the third row shows the
corresponding profile of $h$ along the selected foliation curve. The marked
points are the low-regularity mask produced by the proposed detector and used
in the adaptive closest-point extension step.
}
\label{fig:swe_torus_regularity_detection}
\end{figure}

\subsection{Example 5: PDE-Coupled Regularity Screening for Surface Conservation Laws}
\label{subsec:surface-conservation-laws}

We finally test the proposed regularity detector within a time-dependent surface-PDE computation. We use the closest-point WENO framework of~\cite{Alex2026}, in which the surface conservation law is evolved on a Cartesian embedding grid and the solution is re-extended through the closest-point map. In the present test, this re-extension is adaptively weakened near detected low-regularity regions, while the standard extension is retained elsewhere. The low-regularity mask is generated by the detector developed in this paper.

The WENO flux reconstruction itself is unchanged; the proposed detector is used only to determine the regions where the closest-point extension is modified. We consider two tests on a torus: a scalar Burgers equation with smooth initial data and a reduced shallow-water-type system with a discontinuous initial height. The governing equations and numerical parameters are given in Appendix~\ref{appendix:surface-pde-models}. For the Burgers equation, the mask is computed from $u$, while for the shallow-water-type test it is computed from the height variable $h$. In both cases, local neighborhoods contain $n_X=50$ closest-point samples.

Figure~\ref{fig:burgers_torus_regularity_detection} shows the scalar Burgers
test. The initial data are smooth, so the flagged regions are not inherited from
an initial jump. Instead, as the nonlinear evolution steepens the solution,
narrow low-regularity bands form on the torus and are detected at later times.
This demonstrates that the detector can identify dynamically generated steep
fronts.

Figure~\ref{fig:swe_torus_regularity_detection} shows the reduced
shallow-water-type test. Since the initial height contains a discontinuity in
the angular variable, detected low-regularity points are already present at
early times and remain concentrated near propagating fronts as they wrap around
the torus. The one-dimensional profiles in the last row show that the flagged
points lie near sharp transitions of the height variable \(h\), while smooth
portions of the profile are mostly unflagged.



These examples illustrate that the proposed detector can be used dynamically within a surface-PDE solver to provide the low-regularity mask required by the adaptive closest-point extension.

\section{Conclusion}
\label{sec:conclusion}

We have presented a kernel-based framework for quantitative estimation of local Sobolev regularity
on scattered data. The approach infers an effective regularity order
$\tilde{s}(\Omega_z)$ from the growth of localized native norms as the kernel smoothness parameter~$m$
is increased. The resulting Sobolev-scale profiles provide a direct data-driven measure
of spatially varying regularity, without requiring any prior knowledge of the underlying function or geometry.
A simple transition-based criterion yields $\tilde{s}(\Omega_z)$ efficiently from discrete norm sweeps, while a
band-limited surrogate analysis explains the onset of native-norm amplification once the kernel
smoothness exceeds the locally resolvable regularity.
To enhance spatial selectivity and suppress artifacts due to interface-crossing stencils, an adaptive stencil-shift refinement is incorporated. For large or non-uniform point sets, the accelerated tail-based screening strategy identifies candidate low-regularity neighborhoods before more expensive refinement steps. The numerical comparison with the full norm profile shows that this screening is conservative and less sharply localized, but substantially cheaper. 

Numerical experiments on synthetic benchmarks, turbulent-flow data, and surface conservation-law computations demonstrate the proposed detector in several settings. The benchmark tests show that the method resolves spatially heterogeneous regularity and remains robust under severe non-uniform sampling. The RBF-FD experiments show how the regularity map can guide admissible differentiation stencils. The surface conservation-law example further shows that the detector can be applied to time-dependent PDE solutions and can identify evolving low-regularity regions on the surface.

In summary, the proposed regularity-profiling framework establishes a data-driven connection between kernel native norms, local Sobolev behavior, and stencil reliability. It provides a practical foundation for regularity-aware scattered data approximation, differentiation, and the analysis of numerical PDE data.

\appendix
\section{1D Test Function in  Section \ref{subsec:1d-shape}}
\label{appendix:testfun1d}
The test function $f(x)$ used in Section~\ref{sec:kernel-detector} and Section~\ref{sec:Numerical experiments} is defined on $[-1,1]$ as
\begin{equation}
f(x)=
\begin{cases}
0.8, & -1 \le x < -0.8,\\
1.2, & -0.8 \le x < -0.6,\\
1.2 + 100(x+0.6)^2, & -0.6 \le x < -0.4,\\
5\lvert 5x+1\rvert + 4, & -0.4 \le x < 0,\\
6 + 3\sin(24\pi x), & 0 \le x < \tfrac12,\\
2 + \sin(6\pi x), & \tfrac12 \le x \le 1.
\end{cases}
\label{eq:testfun1d}
\end{equation}

\section{2D Test Function in Section \ref{subsec:example3}}
\label{appendix:testfun}
The test function $f(x,y)$ used in Section~\ref{sec:Numerical experiments}
is defined on $[0,6]^2$.
Let $p_1=(1.5,4.5)$, $p_2=(4.5,2.0)$, $p_3=(3.5,3.5)$ and
$r_i=\sqrt{(x-p_{ix})^2+(y-p_{iy})^2}$ for $i=1,2,3$.
Define
\[
f_1(x,y)=-0.01\sin(4\pi x)\cos(4\pi y)+\exp({-2r_2})+0.2(x+y-2.75),
\]
and set
\begin{equation}
f(x,y)=
\begin{cases}
f_1(x,y), & |x-2.25|+|y-2.25|\le1.75,\\
0.4\!\left(1-{|x-\sin y-4.5|}/{0.3}\right), & |x-\sin y-4.5|\le0.3,\\
1-2r_1, & r_1\le0.5,\\
\dfrac{1}{2}\,{\sin(20r_3)}/({1+3r_3}), & r_3\le1,\\
0.25, & 0.5\le x,y\le5.5,\\
0, & \text{otherwise.}
\end{cases}
\label{eq:testfun2d}
\end{equation}

\section{Surface Conservation-Law Models for Example~\ref{subsec:surface-conservation-laws}}
\label{appendix:surface-pde-models}
Let $\mathbb{T}^2 \subset \mathbb{R}^3$ be a torus with major radius $R=1$ and minor radius $r=0.5$. For $\vec{x} = (x,y,z)^\top \in \mathbb{T}^2$, denote the outward unit normal by $\vec{n}(\vec{x})$ and the azimuthal angle by $\theta(\vec{x}) = \operatorname{atan2}(y,x)$. The tangential vector field $\vec{T}(\vec{x}) = \vec{n}(\vec{x}) \times \vec{V}$ is generated by the constant ambient vector $\vec{V} = (1,-1,0)^\top$.

\subsection{Burgers equation}
\begin{equation}
\partial_t u + \nabla_{\mathbb{T}^2}\!\cdot\!\bigl(\tfrac{u^2}{2} \vec{T}\bigr) = 0,
\qquad t \in (0,12],
\end{equation}
with initial condition
\begin{equation}
u(\vec{x},0) = \sin\theta(\vec{x}) + 0.5.
\end{equation}

\subsection{Reduced shallow-water system}
The reduced shallow-water equations on $\mathbb{T}^2$ are
\begin{equation}
\partial_t \begin{pmatrix} h \\ u \end{pmatrix} 
+ \nabla_{\mathbb{T}^2} \cdot \left( \begin{pmatrix} hu \\ \frac{u^2}{2} + gh \end{pmatrix} \otimes \vec{T} \right) = 0,
\end{equation}
with gravitational constant $g = 9.812$ and $t \in (0,4]$. The initial condition is a fluid at rest,
$u(\vec{x},0) = 0$, with a piecewise-constant depth profile:
\begin{equation}
h(\vec{x},0) =
\begin{cases}
2, & \theta(\vec{x}) < -\pi/4 \;\text{or}\; \theta(\vec{x}) > 3\pi/4,\\
1, & \text{otherwise}.
\end{cases}
\end{equation}

\subsection{Discretization parameters}
The computations use a Cartesian embedding grid with uniform spacing
$\Delta x=0.025$. The spatial terms are discretized by the fifth-order
cp-WENO scheme with global Lax--Friedrichs flux splitting and ENO-based
closest-point extension. Time integration uses the third-order TVD
Runge--Kutta method with CFL number $0.45$. For the regularity plots, the
detector is applied to the computed closest-point solution values using local
neighborhoods of $n_X=50$ samples.

\bibliographystyle{siamplain}
\bibliography{name}

@article{bracco2023discontinuity,
  title={Discontinuity detection by null rules for adaptive surface reconstruction},
  author={Bracco, Cesare and Calabr{\`o}, Francesco and Giannelli, Carlotta},
  journal={J. Sci. Comput.},
  volume={97},
  number={2},
  pages={37},
  year={2023}
}

@article{Kempf+Wendland-Highapprwithkern:23,
author = {Kempf, Rüdiger and Wendland, Holger},
year = {2023},
month = {07},
pages = {1-35},
title = {High-dimensional approximation with kernel-based multilevel methods on sparse grids},
volume = {154},
journal = {Numer. Math.}
}

@article{Wendland-MultanalSobospac:10,
author = {Wendland, Holger},
year = {2010},
month = {09},
pages = {493-517},
title = {Multiscale analysis in {Sobolev} spaces on bounded domains},
volume = {116},
journal = {Numer. Math.}
}

@ARTICLE{HartenEngquistOsherChakravarthy1987ENO,
  author       = {Harten, A. and Engquist, B. and Osher, S. and Chakravarthy, S. R.},
  title        = {Uniformly High Order Accurate Essentially Non-Oscillatory Schemes, {III}},
  journal      = {J. Comput. Phys.},
  volume       = {71},
  number       = {2},
  pages        = {231--303},
  year         = {1987}
}

@ARTICLE{JiangShu1996WENO,
  author       = {Jiang, G.-S. and Shu, C.-W.},
  title        = {Efficient Implementation of Weighted {ENO} Schemes},
  journal      = {J. Comput. Phys.},
  volume       = {126},
  number       = {1},
  pages        = {202--228},
  year         = {1996}
}

@ARTICLE{HuShu1999WENOUnstructured,
  author       = {Hu, C. and Shu, C.-W.},
  title        = {Weighted Essentially Non-Oscillatory Schemes on Triangular Meshes},
  journal      = {J. Comput. Phys.},
  volume       = {150},
  number       = {1},
  pages        = {97--127},
  year         = {1999}
}

@ARTICLE{Avesani2025Smoothness,
  author       = {Avesani, S. and Giacchi, G. and Multerer, M.},
  title        = {Multiresolution local smoothness detection in non-uniformly sampled multivariate signals},
  journal      = {arXiv:2507.13480},
  year         = {2025},
}

@ARTICLE{Harbrecht2022Samplets,
  author       = {Harbrecht, H. and Multerer, M.},
  title        = {Samplets: Construction and scattered data compression},
  journal      = {J. Comput. Phys.},
  volume       = {471},
  pages        = {111616},
  year         = {2022}
}

@ARTICLE{Opfer-Multkern:06,
  author       = {Opfer, R.},
  title        = {Multiscale kernels},
  journal      = { Adv. Comput. Math. },
  volume       = {25},
  number       = {4},
  pages        = {357--380},
  year         = {2006},
}

@ARTICLE{Avesani2024Multiscale,
  author       = {Avesani, S. and Kempf, R. and Multerer, M. and Wendland, H.},
  title        = {Multiscale scattered data analysis in samplet coordinates},
  journal      = {arXiv:2409.14791},
  year         = {2025}
}

@book{wendland2004scattered,
title={{Scattered data approximation}},
author={Wendland, Holger},
volume={17},
year={2004},
publisher={Cambridge University Press}
}

@article{narcowich2006sobolev,
  title={{Sobolev} error estimates and a {Bernstein} inequality for scattered data interpolation via radial basis functions},
  author={Narcowich, Francis J and Ward, Joseph D and Wendland, Holger},
  journal={Constr. Approx.},
  volume={24},
  number={2},
  pages={175--186},
  year={2006}
}

@BOOK{Mallat1999Wavelet,
  author       = {Mallat, S.},
  title        = {A Wavelet Tour of Signal Processing},
  publisher    = {Academic Press},
  address      = {San Diego},
  edition      = {2nd},
  year         = {1999}
}

@ARTICLE{canny1986computational,
  title        = {A computational approach to edge detection},
  author       = {Canny, J. F.},
  journal      = {IEEE Trans. Pattern Anal. Mach. Intell.},
  volume       = {PAMI-8},
  number       = {6},
  pages        = {679--698},
  year         = {1986}
}

@article{jung2009iterative,
title={{An iterative adaptive multiquadric radial basis function method for the detection of local jump discontinuities}},
author={Jung, Jae-Hun and Durante, Vincent R.},
journal={Appl. Numer. Math. },
volume={59},
number={7},
pages={1449--1466},
year={2009},
publisher={Elsevier}
}

@article{jung2011iterative,
title={{Iterative adaptive RBF methods for detection of edges in two-dimensional functions}},
author={Jung, Jae-Hun and Gottlieb, Sigal and Kim, Saeja Oh},
journal={Appl. Numer. Math. },
volume={61},
number={1},
pages={77--91},
year={2011},
publisher={Elsevier}
}

@article{shrivakshan2012comparison,
title={{A comparison of various edge detection techniques used in image processing}},
author={Shrivakshan, G. T. and Chandramouli, Chandrasekar},
journal={Int. J. Comput. Sci. Issues},
volume={9},
number={5},
pages={269},
year={2012}
}

@article{hecht2012new,
  title={{New development in FreeFem++}},
  author={Hecht, Fr{\'e}d{\'e}ric},
  journal={J. Numer. Math.},
  volume={20},
  number={3-4},
  pages={1--14},
  year={2012}
}

@article{demarchi2020shape,
title={{Shape-driven interpolation with discontinuous kernels: Error analysis, edge extraction, and applications in magnetic particle imaging}},
author={De Marchi, Stefano and Erb, Wolfgang and Marchetti, Francesco and Perracchione, Emma and Rossini, Milvia},
journal={SIAM J. Sci. Comput.},
volume={42},
number={2},
pages={B472--B491},
year={2020}
}

@article{tominec2021least,
  title={A least squares radial basis function finite difference method with improved stability properties},
  author={Tominec, Igor and Larsson, Elisabeth and Heryudono, Alfa},
  journal={SIAM J. Sci. Comput.},
  volume={43},
  number={2},
  pages={A1441--A1471},
  year={2021}
}

@article{de2020jumping,
  title={{Jumping with variably scaled discontinuous kernels (VSDKs)}},
  author={De Marchi, Stefano and Marchetti, Francesco and Perracchione, Emma},
  journal={BIT Num. Math.},
  volume={60},
  number={2},
  pages={441--463},
  year={2020}
}

@article{bracco2025positive,
  title={A positive meshless finite difference scheme for scalar conservation laws with adaptive artificial viscosity driven by fault detection},
  author={Bracco, Cesare and Davydov, Oleg and Giannelli, Carlotta and Sestini, Alessandra},
  journal={Comput. Math. Appl.},
  volume={190},
  pages={103--121},
  year={2025}
}

@article{karvonen2025scale,
   title={Scale estimation and rate-unbiasedness for {Gaussian} processes under smoothness misspecification}, 
      author={Karvonen, Toni and  Bachoc, François},
      year={2025},
      journal={arXiv preprint arXiv: 2110.02810}
}

@article{Alex2026,
  title        = {Foliation structures and global flow dynamics of scalar hyperbolic conservation laws on manifolds: I. geometry-compatible fluxes and numerical validation on sphere and torus},
  author       = {Wang, Bao–Shan and Chu, Alex Shiu Lun and Ling, Leevan and Don, Wai Sun},
  journal      = {Commun. Nonlinear Sci. Numer. Simul.},
  volume       = {152},
  pages        = {109415},
  year         = {2026},
}

@article{ArchibaldGelbYoon2005,
author = {Archibald, Rick and Gelb, Anne and Yoon, Jungho},
title = {Polynomial Fitting for Edge Detection in Irregularly Sampled Signals and Images},
journal = {SIAM J. on Numer. Analy.},
volume = {43},
number = {1},
pages = {259-279},
year = {2005},
}

\end{document}